\documentclass[twocolumn,final]{IEEEtran}
\usepackage{graphicx}
\usepackage{hyperref} 
\usepackage{amsmath}
\usepackage{amssymb}
\usepackage{theorem}
\usepackage{euscript}
\usepackage{xcolor}
\definecolor{color2b}{rgb}{0.6,0.8,1.0}
\usepackage{epsfig}
\usepackage{mathtools}
\usepackage[usenames,dvipsnames]{pstricks}
\usepackage{pstricks-add}
\usepackage{pgf,tikz}
\usepackage{pst-grad}
\usepackage{pst-plot} 
\usepackage{pst-node} 
\usepackage{pst-eucl} 
\usepackage{pst-coil} 
\usepackage{mathrsfs} 
\usepackage{pgf,tikz}
\usepackage{caption} 
\definecolor{labelkey}{rgb}{0,0.08,0.45}
\definecolor{refkey}{rgb}{0,0.6,0.0}
\definecolor{dblue}{HTML}{0455BF}
\definecolor{dgreen}{HTML}{02724A}
\definecolor{Dblue}{HTML}{8602DC}
\definecolor{dred}{HTML}{D90404}
\usepackage{upref}
\usepackage{hyperref}
\hypersetup{colorlinks=true,linktocpage=true,linkcolor=dgreen,%
citecolor=dblue,urlcolor=dred}
\renewcommand{\leq}{\ensuremath{\leqslant}}
\renewcommand{\geq}{\ensuremath{\geqslant}}

\newcommand{\menge}[2]{\big\{{#1}~\big |~{#2}\big\}} 
 
\newcommand{\emp}{\ensuremath{\varnothing}}
\newcommand{\scal}[2]{{\left\langle{{#1}\mid{#2}}\right\rangle}}
\newcommand{\sscal}[2]{{\big\langle{{#1}\mid{#2}}\big\rangle}}
\newcommand{\HH}{\mathcal H}
\newcommand{\GG}{\mathcal G}

\newcommand{\HHH}{\ensuremath{\boldsymbol{\mathcal H}}}
\newcommand{\GGG}{\ensuremath{\boldsymbol{\mathcal G}}}
\newcommand{\KKK}{\ensuremath{\boldsymbol{\mathcal K}}}
\newcommand{\zeroun}{\ensuremath{\left]0,1\right[}} 
\newcommand{\rzeroun}{\ensuremath{\left]0,1\right]}} 
\newcommand{\exi}{\ensuremath{\exists\,}}
\newcommand{\nabli}[1]{\nabla_{\!{#1}}\,}

\newcommand{\NN}{\ensuremath{\mathbb N}}
\newcommand{\RR}{\ensuremath{\mathbb R}}
\newcommand{\RP}{\ensuremath{\left[0,+\infty\right[}}

\newcommand{\RPP}{\ensuremath{\left]0,+\infty\right[}}
\newcommand{\RX}{\ensuremath{\left]-\infty,+\infty\right]}}
\newcommand{\RXX}{\ensuremath{\left[-\infty,+\infty\right]}}
\newcommand{\Argmin}{\ensuremath{\mathrm{Argmin}\,}}
\newcommand{\reli}{\ensuremath{\operatorname{ri}}}
\newcommand{\epi}{\ensuremath{\mathrm{epi}\,}}
\newcommand{\dom}{\ensuremath{\mathrm{dom}\,}}
\newcommand{\intdom}{\ensuremath{\mathrm{int\,dom}\,}}
\newcommand{\ran}{\ensuremath{\mathrm{range}\,}}
\newcommand{\prox}{\ensuremath{\mathrm{prox}}}
\newcommand{\best}{\ensuremath{\mathrm{best}}}
\newcommand{\proj}{\ensuremath{\mathrm{proj}}}

\newcommand{\gra}{\ensuremath{\mathrm{gra}\,}}
\newcommand{\Sum}{\ensuremath{\displaystyle\sum}}

\newcommand{\Fix}{\ensuremath{\text{\rm Fix}\,}}
\newcommand{\Id}{\ensuremath{\mathrm{Id}}}
\newcommand{\ID}{\ensuremath{\boldsymbol{\mathrm{Id}}}}

\newcommand{\pinf}{\ensuremath{+\infty}}

\newcommand{\EC}[2]{{\mathsf E}(#1\! \mid\! #2)}

\newcommand{\XX}{\ensuremath{\EuScript{X}}}
\newcommand{\as}{\ensuremath{\text{a.~\!s.}}}
\newcommand{\EE}{\ensuremath{\mathsf E}}
\newcommand{\PP}{\ensuremath{\mathsf{Prob}\,}}
\newcommand{\minimize}[2]{\ensuremath{\underset{\substack{{#1}}}%
{\text{\rm minimize}}\;\;#2 }}
\newcommand{\cart}{\ensuremath{\raisebox{-0.5mm}%
{\mbox{\LARGE{$\times$}}}}}
\newcommand{\Argmind}[2]{\ensuremath{\underset{\substack{{#1}}}%
{\text{\rm Argmin}}\;\;#2 }}

\newtheorem{theorem}{Theorem}

\newtheorem{proposition}[theorem]{Proposition}
\theoremstyle{plain}{\theorembodyfont{\rmfamily}%
\newtheorem{assumption}[theorem]{Assumption}}
\theoremstyle{plain}{\theorembodyfont{\rmfamily}%
\newtheorem{definition}[theorem]{Definition}}
\theoremstyle{plain}{\theorembodyfont{\rmfamily}%
\newtheorem{example}[theorem]{Example}}
\theoremstyle{plain}{\theorembodyfont{\rmfamily}%
\newtheorem{problem}[theorem]{Problem}}
\theoremstyle{plain}{\theorembodyfont{\rmfamily}%
\newtheorem{remark}[theorem]{Remark}}
\theoremstyle{plain}{\theorembodyfont{\rmfamily}%
\newtheorem{algorithm}[theorem]{Algorithm}}
\renewcommand\theenumi{\roman{enumi})}
\renewcommand\theenumii{\alph{enumii})}

\newcounter{numperceuse}

\definecolor{dred}{rgb}{0.80,0.0,0.00}
\definecolor{nido}{rgb}{0.60,0.0,0.60}
\newrgbcolor{dgreen}{0.00 0.35 0.00}
\newrgbcolor{dbrown}{0.65 0.40 0.35}

\begin{document}

\title{Fixed Point Strategies in Data Science}
\author{Patrick L. Combettes, {\em Fellow, IEEE}, 
and Jean-Christophe Pesquet, {\em Fellow, IEEE}
\thanks{The work of P. L. Combettes was supported by the National
Science Foundation under grant CCF-1715671. The work of J.-C.
Pesquet was supported by Institut Universitaire de France and the
ANR Chair in Artificial Intelligence BRIGEABLE.}
\thanks{P. L. Combettes (corresponding author) is with 
North Carolina State University,
Department of Mathematics, Raleigh, NC 27695-8205, USA
and J.-C. Pesquet is with Universit\'e Paris-Saclay, 
Inria, CentraleSup\'elec, 
Centre de Vision Num\'erique, 91190 Gif sur Yvette,
France. E-mail: 
\texttt{plc@math.ncsu.edu}, \texttt{jean-christophe@pesquet.eu}.}
}
\markboth{IEEE Transactions on Signal Processing}%
{Combettes and Pesquet: Fixed Point Methods}
\maketitle

\begin{abstract}
The goal of this paper is to promote the use of fixed point
strategies in data science by showing that they provide a
simplifying and unifying framework to model, analyze, and solve a
great variety of problems. They are seen to constitute a natural
environment to explain the behavior of advanced convex optimization
methods as well as of recent nonlinear methods in data science
which are formulated in terms of paradigms that go beyond
minimization concepts and involve constructs such as Nash
equilibria or monotone inclusions. We review the pertinent tools
of fixed point theory and describe the main state-of-the-art
algorithms for provably convergent fixed point construction. We
also incorporate additional ingredients such as stochasticity,
block-implementations, and non-Euclidean metrics, which provide
further enhancements. Applications to signal and image processing,
machine learning, statistics, neural networks, and inverse problems
are discussed. 
\end{abstract}

\begin{IEEEkeywords}
Convex optimization,
fixed point,
game theory,
monotone inclusion,
image recovery,
inverse problems,
machine learning,
neural networks,
nonexpansive operator,
signal processing.
\end{IEEEkeywords}

\section{Introduction}

Attempts to apply mathematical methods to the extraction of
information from data can be traced back to the work of Boscovich
\cite{Bosc57}, Gauss \cite{Gaus09}, Laplace \cite{Lapl89}, and
Legendre \cite{Lege05}. Thus, in connection with the problem of
estimating parameters from noisy observations, Boscovich and
Laplace invented the least-deviations data fitting method, while
Legendre and Gauss invented the least-squares data fitting method.
On the algorithmic side, the gradient method was invented by Cauchy
\cite{Cauc47} to solve a data fitting problem in astronomy, and
more or less heuristic methods have been used from then on. The
early work involving provably convergent numerical solutions
methods was focused mostly on quadratic minimization problems or
linear programming techniques, e.g.,
\cite{Artz79,Herm73,Hunt70,Twom65,Wagn59}. Nowadays, general convex
optimization methods have penetrated virtually all branches of data
science 
\cite{Bach12,Byrn14,Cham16,Aiep96,Banf11,Glow16,Wrig12,Theo20}. 
In fact, the optimization and data science communities have never
been closer, which greatly facilitates technology transfers towards
applications. Reciprocally, many of the recent advances in convex
optimization algorithms have been motivated by data processing
problems in signal recovery, inverse problems, or machine learning.
At the same time, the design and the convergence analysis of some
of the most potent splitting methods in highly structured or
large-scale optimization are based on concepts that are not
found in the traditional optimization toolbox but reach deeper into
nonlinear analysis. Furthermore, an increasing number of problem
formulations go beyond optimization in the sense that their
solutions are not optimal in the classical sense of minimizing a
function but, rather, satisfy more general notions of equilibrium.
Among the formulations that fall outside of the realm of standard
minimization methods, let us mention variational inequality and
monotone inclusion models, game theoretic approaches, neural
network structures, and plug-and-play methods. 

Given the abundance of activity described above and the
increasingly complex formulations of some data processing problems
and their solution methods, it is essential to identify general
structures and principles in order to simplify and clarify the
state of the art. It is the objective of the present paper to
promote the viewpoint that fixed point theory constitutes an ideal
technology towards this goal. Besides its unifying nature, the
fixed point framework offers several advantages. On the algorithmic
front, it leads to powerful convergence principles that demystify
the design and the asymptotic analysis of iterative methods.
Furthermore, fixed point methods can be implemented using
stochastic perturbations, as well as block-coordinate or
block-iterative strategies which reduce the computational load and
memory requirements of the iterations. 

Historically, one of the first uses of fixed point theory in signal
recovery is found in the bandlimited reconstruction method of
\cite{Land61}, which is based on the iterative Banach-Picard
contraction process
\begin{equation}
\label{e:emile}
x_{n+1}=Tx_n,
\end{equation}
where the operator $T$ has Lipschitz constant $\delta<1$. The
importance of dealing with the more general class of nonexpansive
operators, i.e., those with Lipschitz constant $\delta=1$,
was emphasized by Youla in \cite{Youl78} and \cite{Youl82}; 
see also \cite{Scha81,Tomv81,Wile78}. Since then, many problems in
data science have been modeled and solved using nonexpansive
operator theory; see for instance 
\cite{Fixe11,Byrn14,Aiep96,Smds20,Smms05,Daub04,Mari90,%
Pott93,Star87,Theo11}.

The outline of the paper is as follows. In order to make
the paper as self-contained as possible, we present in
Section~\ref{sec:2} the essential tools and results from nonlinear
analysis on which fixed point approaches are grounded. These
include notions of convex analysis, monotone operator theory, and
averaged operator theory. Section~\ref{sec:3} provides an overview
of basic fixed point principles and methods. Section~\ref{sec:4}
addresses the broad class of monotone inclusion problems and their
fixed point modeling. Using the tools of Section~\ref{sec:3},
various splitting strategies are described, as well as
block-iterative and block-coordinate algorithms.
Section~\ref{sec:5} discusses applications of splitting methods to
a large panel of techniques for solving structured convex
optimization problems. Moving beyond traditional optimization,
algorithms for Nash equilibria are investigated in
Section~\ref{sec:6}. Section~\ref{sec:7} shows how fixed point
strategies can be applied to four additional categories of data
science problems that have no underlying minimization
interpretation. Some brief conclusions are drawn in
Section~\ref{sec:8}. For simplicity, we have adopted a Euclidean
space setting. However, most results remain valid in general
Hilbert spaces up to technical adjustments.

\section{Notation and mathematical foundations}
\label{sec:2}
We review the basic tools and principles from nonlinear analysis
that will be used throughout the paper. Unless otherwise stated, 
the material of this section can be found in \cite{Livre1}; for 
convex analysis see also \cite{Rock70}. 

\subsection{Notation}
\label{sec:not}
Throughout, $\HH$, $\GG$, $(\HH_i)_{1\leq i\leq m}$,
and $(\GG_k)_{1\leq k\leq q}$ are Euclidean spaces. 
We denote by $2^\HH$ the collection of all subsets of $\HH$ and 
by $\HHH=\HH_1\times\cdots\times\HH_m$ 
and $\GGG=\GG_1\times\cdots\times\GG_q$ the standard Euclidean 
product spaces. A generic point in $\HHH$ is denoted by 
$\boldsymbol{x}=(x_i)_{1\leq i\leq m}$. 
The scalar product of a Euclidean space is denoted by
$\scal{\cdot}{\cdot}$ and the associated norm by $\|\cdot\|$.
The adjoint of a linear operator $L$ is denoted by $L^*$.
Let $C$ be a subset of $\HH$. Then the 
\emph{distance function} to $C$ 
is $d_C\colon x\mapsto\inf_{y\in C}\|x-y\|$ and 
the \emph{relative interior} of $C$, denoted
by $\reli C$, is its interior relative to its affine hull.

\subsection{Convex analysis}
\label{sec:21}
The central notion in convex analysis is that of a convex set: a
subset $C$ of $\HH$ is \emph{convex} if it contains all the line
segments with end points in the set, that is, 
\begin{equation}
\label{e:convex1}
(\forall x\in C)(\forall y\in C)(\forall\alpha\in\zeroun\,)\quad
\alpha x+(1-\alpha)y\in C.
\end{equation}
The projection theorem is one of the most important results
of convex analysis. 

\begin{theorem}[projection theorem]
\label{t:11}
Let $C$ be a nonempty closed convex subset of $\HH$ and let
$x\in\HH$. Then there exists a unique point $\proj_Cx\in C$, called
the \emph{projection} of $x$ onto $C$, such that
$\|x-\proj_Cx\|=d_C(x)$.
In addition, for every $p\in\HH$, 
\begin{equation}
\label{e:proj}
p=\proj_Cx\quad\Leftrightarrow\quad
\begin{cases}
p\in C\\
(\forall y\in C)\;\scal{y-p}{x-p}\leq 0.
\end{cases}
\end{equation}
\end{theorem}

Convexity for functions is inherited from convexity for sets as
follows. Consider a function $f\colon\HH\to\RX$. 
Then $f$ is \emph{convex} if its \emph{epigraph}
\begin{equation}
\label{e:epi}
\epi f=\menge{(x,\xi)\in\HH\times\RR}{f(x)\leq\xi}
\end{equation}
is a convex set. This is equivalent to requiring that
\begin{multline}
\label{e:convex}
(\forall x\in\HH)(\forall y\in\HH)(\forall\alpha\in\zeroun\,)
\quad\\
f\big(\alpha x+(1-\alpha) y\big)\leq\alpha f(x)+(1-\alpha)f(y).
\end{multline}
If $\epi f$ is closed, then $f$ is 
\emph{lower semicontinuous} in the sense that, for every sequence 
$(x_n)_{n\in\NN}$ in $\HH$ and $x\in\HH$, 
\begin{equation}
\label{e:lsc}
x_n\to x\quad\Rightarrow\quad
f(x)\leq\varliminf f(x_n).
\end{equation}
Finally, we say that $f\colon\HH\to\RX$ is \emph{proper} if 
$\epi f\neq\emp$, which is equivalent to 
\begin{equation}
\label{e:dom}
\dom f=\menge{x\in\HH}{f(x)<\pinf}\neq\emp.
\end{equation}
The class of functions $f\colon\HH\to\RX$ which are proper, lower
semicontinuous, and convex is denoted by $\Gamma_0(\HH)$. The
following result is due to Moreau \cite{Mor62b}.

\begin{theorem}[proximation theorem]
\label{t:12}
Let $f\in\Gamma_0(\HH)$ and let $x\in\HH$. Then there exists a
unique point $\prox_fx\in\HH$, called the \emph{proximal point} of
$x$ relative to $f$, such that
\begin{multline}
\label{e:prox}
f\big(\prox_fx\big)+\frac12\|x-\prox_fx\|^2=\\
\min_{y\in\HH}\bigg(f(y)+\frac{1}{2}\|x-y\|^2\bigg).
\end{multline}
In addition, for every $p\in\HH$, 
\begin{multline}
\label{e:moreau1}
p=\prox_fx\;\Leftrightarrow\\
(\forall y\in\HH)\;\scal{y-p}{x-p}+f(p)\leq f(y).
\end{multline}
\end{theorem}

\begin{figure}
\scalebox{0.63} 
{
\begin{pspicture}(-2.9,-3.5)(4.0,4.3)
\psline[linewidth=0.04cm,arrowsize=0.09cm 4.0,%
arrowlength=1.4,arrowinset=0.4]{->}(-0.0,-3.0)(-0.0,3.2)
\psline[linewidth=0.04cm,arrowsize=0.09cm 4.0,%
arrowlength=1.4,arrowinset=0.4]{->}(-2.7,-0.0)(3.0,-0.0)
\rput(3.3,0.0){\large $\HH$}
\rput(0.0,3.5){\large $\RR$}
\rput(-1.0,0.0){\large $\boldsymbol{|}$}
\rput(-1.0,0.5){\large $-1$}
\rput(0.0,-1.0){\large $\boldsymbol{-}$}
\rput(0.6,-1.0){\large $-1$}
\rput(1.0,0.0){\large $\boldsymbol{|}$}
\rput(1.0,0.5){\large $1$}
\psplot[linewidth=0.05cm,linestyle=solid,algebraic,%
linecolor=dbrown]{-2.6}{-0.99}{-2*x-3}
\psline[linewidth=0.05cm,linecolor=dbrown](-1.0,-1.0)(1.01,0.0)
\psplot[linewidth=0.05cm,linestyle=solid,algebraic,%
linecolor=dbrown]{1.0}{2.6}{0.5*x^2-0.5}
\end{pspicture}

\begin{pspicture}(-2.9,-3.5)(4.0,4.3)
\psline[linewidth=0.04cm,arrowsize=0.09cm 4.0,%
arrowlength=1.4,arrowinset=0.4]{->}(-0.0,-3.0)(-0.0,3.2)
\psline[linewidth=0.04cm,arrowsize=0.09cm 4.0,%
arrowlength=1.4,arrowinset=0.4]{->}(-2.7,-0.0)(3.0,-0.0)
\rput(3.3,0.0){\large $\HH$}
\rput(0.0,3.5){\large $\HH$}
\psline[linewidth=0.05cm,linecolor=blue](-2.6,-2.0)(-1.0,-1.97)
\psline[linewidth=0.05cm,linecolor=blue](-1.0,-1.995)(-1.0,0.525)
\psline[linewidth=0.05cm,linecolor=blue](-1.0,0.5)(1.024,0.5)
\psline[linewidth=0.05cm,linecolor=blue](1.0,0.5)(1.0,1.01)
\psplot[linewidth=0.05cm,linestyle=solid,algebraic,linecolor=blue]%
{0.995}{2.6}{x}
\rput(0.0,-2.0){\large $\boldsymbol{-}$}
\rput(0.6,-2.0){\large $-2$}
\rput(0.0,1.0){\large $\boldsymbol{-}$}
\rput(0.4,1.0){\large $1$}
\rput(1.0,0.0){\large $\boldsymbol{|}$}
\rput(1.0,-0.5){\large $1$}
\end{pspicture} 
}
\caption{Left: Graph of a function defined on $\HH=\RR$. 
Right: Graph of its subdifferential.}
\label{fig:7}
\end{figure}
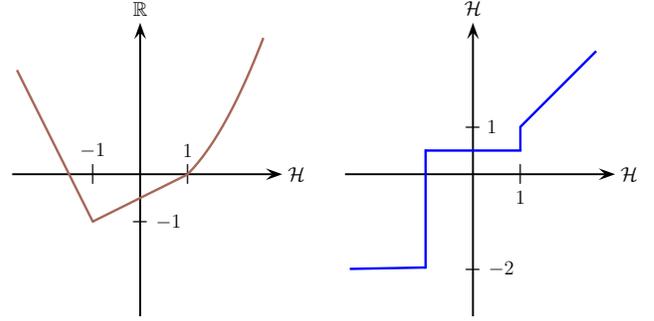

The above theorem defines an operator $\prox_f$ called the
\emph{proximity operator} of $f$ (see \cite{Banf11} for a
tutorial, and \cite[Chapter~24]{Livre1} and \cite{Comb18} for a
detailed account with various properties). 
\begin{figure}
\scalebox{0.53} 
{
\begin{pspicture}(-6.0,-2.0)(11.5,11.1) 
\psplot[linewidth=0.06cm,linestyle=solid,algebraic,%
linecolor=dbrown]{-2.32}{6.3}{abs(x/2-1)+(x/2-1)^2+3}
\psplot[linewidth=0.05cm,linestyle=solid,algebraic,%
linecolor=red]{-2.3}{6.7}{1.2*x+1.22}
\psplot[linewidth=0.05cm,linestyle=solid,algebraic,%
linecolor=blue]{-3.0}{6.6}{x/2+2}
\psline[linewidth=0.05cm,arrowsize=3.0mm,%
linestyle=solid,linecolor=dgreen]{->}%
(3.10,3.85)(3.10,4.9)
\psline[linecolor=black,linewidth=0.045cm,linestyle=dashed]{-}%
(2.00,0.0)(2.00,3)
\psline[linecolor=black,linewidth=0.045cm,linestyle=dashed]{-}%
(-1,3.0)(2.00,3)
\psline[linewidth=0.05cm,arrowsize=0.11cm 4.0,%
arrowlength=1.4,arrowinset=0.4]{->}(-4,0)(7.2,0)
\psline[linewidth=0.05cm,arrowsize=0.11cm 4.0,%
arrowlength=1.4,arrowinset=0.4]{->}(-1,-1)(-1,10.2)
\rput(5.2,5.5){\LARGE \color{dbrown}$\gra f$}
\rput(-2.7,1.5){\LARGE \color{blue}$\gra{{m_{x,w}}}$}
\rput(7.7,8.6){\LARGE \color{red}$\gra{{\scal{\cdot}{u}}}$}
\rput(2.00,-0.4){\LARGE $x$}
\rput(1.70,8.2){\LARGE \color{dbrown}$\epi f$}
\rput(-1.7,3.0){\LARGE \color{black}$f(x)$}
\rput(-1,10.5){\LARGE $\RR$}
\rput(7.5,0){\LARGE $\HH$}
\rput(2.3,5.2){\LARGE \color{dgreen}{${f^*(u)}$}}
\end{pspicture} 
}
\caption{The graph of a function $f\in\Gamma_0(\HH)$ is shown in
brown. The area above the graph is the closed convex set 
$\epi f$ of \eqref{e:epi}. 
Let $u\in\HH$ and let the red line be the graph of the linear
function $\scal{\cdot}{u}$. In view of \eqref{e:conj}, the
value of $f^*(u)$ (in green) is the maximum signed difference
between the red line and the brown line. 
Now fix $x\in\HH$ and $w\in\partial f(x)$. Additionally, by
\eqref{e:subdiff}, the affine function 
$m_{x,w}\colon y\mapsto\scal{y-x}{w}+f(x)$ satisfies 
$m_{x,w}\leq f$ and it coincides with $f$ at $x$. 
Its graph is represented in blue. Every subgradient $w$ gives such
an affine minorant.
}
\label{fig:8}
\end{figure}
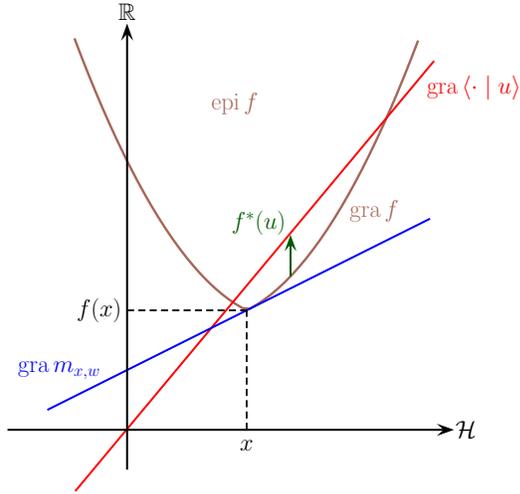
Now let $C$ be a nonempty closed convex subset of $\HH$. Then its
\emph{indicator function} $\iota_C$, defined by
\begin{equation}
\label{e:iota}
\iota_C\colon\HH\to\RX\colon x\mapsto 
\begin{cases}
0,&\text{if}\;\:x\in C;\\
\pinf,&\text{if}\;\:x\notin C,
\end{cases}
\end{equation}
lies in $\Gamma_0(\HH)$ and it follows from \eqref{e:proj} and
\eqref{e:moreau1} that
\begin{equation}
\label{e:98}
\prox_{\iota_C}=\proj_C. 
\end{equation}
This shows that Theorem~\ref{t:12} generalizes Theorem~\ref{t:11}.
Let us now introduce basic convex analytical tools (see
Fig.~\ref{fig:8}). The \emph{conjugate} of $f\colon\HH\to\RX$ is
\begin{equation}
\label{e:conj}
f^*\colon\HH\to\RXX\colon u
\mapsto\sup_{x\in\HH}\big(\scal{x}{u}-f(x)\big).
\end{equation}
The \emph{subdifferential} of a proper function $f\colon\HH\to\RX$
is the set-valued operator
$\partial f\colon\HH\to 2^{\HH}$ which maps a point $x\in\HH$ to
the set (see Fig.~\ref{fig:7})
\begin{equation}
\label{e:subdiff}
\partial f(x)\!=\!
\menge{u\in\HH}
{(\forall y\in\HH)\;\scal{y-x}{u}+f(x)\leq f(y)}.
\end{equation}
A vector in $\partial f(x)$ is a \emph{subgradient} of $f$ at $x$.
If $C$ is a nonempty closed convex subset of $\HH$,
$N_C=\partial\iota_C$ is the \emph{normal cone} operator of 
$C$, that is, for every $x\in\HH$,
\begin{equation} 
\label{e:normalcone}
N_Cx=
\begin{cases}
\menge{u\in\HH}{\!(\forall y\in C)\:\scal{y-x}{u}\leq 0},
&\text{if}\;\;x\in C;\\
\emp,&\text{otherwise.}
\end{cases}
\end{equation}
Let us denote by $\Argmin f$ the set of minimizers of a 
function $f\colon\HH\to\RX$ (the notation 
$\Argmin_{x\in\HH}f(x)$ will also be used).
The most fundamental result in optimization is actually the
following immediate consequence of \eqref{e:subdiff}. 

\begin{theorem}[Fermat's rule]
\label{t:4}
Let $f\colon\HH\to\RX$ be a proper function. Then 
$\Argmin f=\menge{x\in\HH}{0\in\partial f(x)}$.
\end{theorem}

\begin{theorem}[Moreau]
\label{t:6}
Let $f\in\Gamma_0(\HH)$. Then $f^*\in\Gamma_0(\HH)$, 
$f^{**}=f$, and $\prox_f+\prox_{f^*}=\Id$.
\end{theorem}

A function $f\in\Gamma_0(\HH)$ is \emph{differentiable} at 
$x\in\dom f$ if there exists a vector $\nabla f(x)\in\HH$, 
called the \emph{gradient} of $f$ at $x$, such that
\begin{equation}
\label{e:grad}
(\forall y\in\HH)\quad
\lim_{\alpha\downarrow 0}\dfrac{f(x+\alpha y)-f(x)}{\alpha}=
\scal{y}{\nabla f(x)}.
\end{equation}

\begin{example}
\label{ex:jjm7}
Let $C$ be a nonempty closed convex subset of $\HH$. Then 
$\nabla d_C^2/2=\Id-\proj_C$.
\end{example}

\begin{proposition}
\label{p:11}
Let $f\in\Gamma_0(\HH)$, let $x\in\dom f$, and suppose that $f$
is differentiable at $x$. Then $\partial f(x)=\{\nabla f(x)\}$.
\end{proposition}

We close this section by examining fundamental properties of a
canonical convex minimization problem. 

\begin{proposition}
\label{p:17}
Let $f\in\Gamma_0(\HH)$, let $g\in\Gamma_0(\GG)$, and let
$L\colon\HH\to\GG$ be linear. Suppose that
$L(\dom f)\cap\dom g\neq\emp$ and set $S=\Argmin(f+g\circ L)$.
Then the following hold: 
\begin{enumerate}
\item
\label{p:17i}
Suppose that $\lim_{\|x\|\to\pinf}f(x)+g(Lx)=\pinf$. Then 
$S\neq\emp$.
\item
\label{p:17ii}
Suppose that $\reli(L(\dom f))\cap\reli(\dom g)\neq\emp$. Then 
\begin{align*}
\label{e:17ii}
S
&=\menge{x\in\HH}{0\in\partial f(x)+L^*\big(\partial g(Lx)\big)}\\
&=\menge{x\in\HH}{(\exi v\in\partial g(Lx))\;-L^*v\in\partial
f(x)}.
\end{align*}
\end{enumerate}
\end{proposition}

\subsection{Nonexpansive operators}
\label{sec:23}
We introduce the main classes of operators pertinent to our
discussion. First, we need to define the notion of a relaxation for
an operator. 

\begin{definition}
\label{d:relax1}
Let $T\colon\HH\to\HH$ and let $\lambda\in\RPP$. Then the operator
$R=\Id+\lambda(T-\Id)$ is a \emph{relaxation} of $T$. If
$\lambda\leq 1$, then $R$ is an \emph{underrelaxation} of $T$ and,
if $\lambda\geq 1$, $R$ is an \emph{overrelaxation} of $T$; in
particular, if $\lambda=2$, $R$ is the \emph{reflection} of $T$.
\end{definition}

\begin{definition}
\label{d:relax2}
Let $\alpha\in\rzeroun$. An $\alpha$-relaxation sequence is a
sequence $(\lambda_n)_{n\in\NN}$ in $\left]0,1/\alpha\right[$ such
that $\sum_{n\in\NN}\lambda_n(1-\alpha\lambda_n)=\pinf$. 
\end{definition}

\begin{example}
\label{ex:relax}
Let $\alpha\in\rzeroun$ and let $(\lambda_n)_{n\in\NN}$ be a
sequence in $\RPP$. Then $(\lambda_n)_{n\in\NN}$ is an
$\alpha$-relaxation sequence in each of the following cases:
\begin{enumerate}
\item
\label{ex:relaxi}
$\alpha<1$ and $(\forall n\in\NN)$ $\lambda_n=1$.
\item
\label{ex:relaxii}
$(\forall n\in\NN)$ $\lambda_n=\lambda\in\left]0,1/\alpha\right[$.
\item
\label{ex:relaxiii}
$\inf_{n\in\NN}\lambda_n>0$ and
$\sup_{n\in\NN}\lambda_n<1/\alpha$.
\item
\label{ex:relaxiv}
There exists $\varepsilon\in\zeroun$ such that 
$(\forall n\in\NN)$ $\varepsilon/\sqrt{n+1}\leq\lambda_n\leq
1/\alpha-\varepsilon/\sqrt{n+1}$.
\end{enumerate}
\end{example}

An operator $T\colon\HH\to\HH$ is 
\emph{Lipschitzian} with constant $\delta\in\RPP$ if 
\begin{equation}
\label{e:100}
(\forall x\in\HH)(\forall y\in\HH)\quad\|Tx-Ty\|\leq\delta\|x-y\|.
\end{equation}
If $\delta<1$ above, then $T$ is a 
\emph{Banach contraction} (also called a \emph{strict
contraction}). If $\delta=1$, that is, 
\begin{equation}
\label{e:103}
(\forall x\in\HH)(\forall y\in\HH)\quad\|Tx-Ty\|\leq\|x-y\|,
\end{equation}
then $T$ is \emph{nonexpansive}.
On the other hand, $T$ is \emph{cocoercive} with constant 
$\beta\in\RPP$ if
\begin{multline}
\label{e:101}
(\forall x\in\HH)(\forall y\in\HH)\quad
\scal{x-y}{Tx-Ty}\geq\\
\beta\|Tx-Ty\|^2.
\end{multline}
If $\beta=1$ in \eqref{e:101}, then
$T$ is \emph{firmly nonexpansive}. Alternatively, $T$ is 
firmly nonexpansive if 
\begin{multline}
\label{e:107}
(\forall x\in\HH)(\forall y\in\HH)\;\;
\|Tx-Ty\|^2\leq\|x-y\|^2\\
-\|(\Id-T)x-(\Id-T)y\|^2.
\end{multline}
Equivalently, $T$ is firmly nonexpansive if the reflection
\begin{equation}
\label{e:104}
\Id+2(T-\Id)\;\:\text{is nonexpansive.}
\end{equation}
More generally, let $\alpha\in\rzeroun$. Then $T$ is
$\alpha$-\emph{averaged} if the overrelaxation
\begin{equation}
\label{e:105-}
\Id+\alpha^{-1}(T-\Id)\;\:\text{is nonexpansive}
\end{equation}
or, equivalently, if there exists a nonexpansive operator
$Q\colon\HH\to\HH$ such that $T$ can be written as
the underrelaxation
\begin{equation}
\label{e:105}
T=\Id+\alpha(Q-\Id).
\end{equation}
An alternative characterization of $\alpha$-averagedness is
\begin{multline}
\label{e:106}
(\forall x\in\HH)(\forall y\in\HH)\;\;
\|Tx-Ty\|^2\leq\|x-y\|^2\\
-{\dfrac{1-\alpha}{\alpha}}\|(\Id-T)x-(\Id-T)y\|^2.
\end{multline}
Averaged operators will be the most important class of nonlinear 
operators we use in this paper. They were introduced in 
\cite{Bail78} and their central role in many nonlinear analysis 
algorithms was pointed out in \cite{Opti04}, with further
refinements in \cite{Jmaa15,Huan20}. Note that 
\begin{align}
\label{e:110}
T\;\text{is firmly nonexpansive}
&\Leftrightarrow\:\Id-T\;\text{is firmly nonexpansive}
\nonumber\\
&\Leftrightarrow\:T\;\text{is $1/2$-averaged}
\nonumber\\
&\Leftrightarrow\:T\;\text{is $1$-cocoercive}.
\end{align}
Here is an immediate consequence of \eqref{e:moreau1} and 
\eqref{e:110}.

\begin{example}
\label{ex:1}
Let $f\in\Gamma_0(\HH)$. Then 
$\prox_f$ and $\Id-\prox_f$ are firmly nonexpansive. 
In particular, if $C$ is a nonempty closed convex subset of $\HH$,
then \eqref{e:98} implies that $\proj_C$ and $\Id-\proj_C$ are
firmly nonexpansive.
\end{example}

\begin{figure}
\begin{center}
\scalebox{0.82} 
{
\begin{pspicture}(-0.9,-1.2)(9.5,7.1)
\rput(4.0,1.0){\large\color{nido} projection operators}
\rput(4.0,1.9){\large\color{nido} proximity operators}
\rput(4.0,3.2){\large\color{nido} firmly nonexpansive }
\rput(4.0,2.7){\large\color{nido} operators/resolvents }
\rput(2.7,4.2){\large\color{dred} $\alpha$-averaged operators,
$\alpha<1$}
\rput(1.7,5.2){\large\color{dred} nonexpansive operators}
\rput(5.0,6.2){\large Lipschitzian operators}
\rput(8.0,3.0){\large\color{blue}cocoercive}
\rput(7.9,2.5){\large\color{blue}operators}
\psframe[linecolor=nido,linewidth=0.04,dimen=outer]
(2.0,0.5)(6.0,1.5)
\psframe[linecolor=nido,linewidth=0.04,dimen=outer]
(1.8,0.3)(6.2,2.4)
\psframe[linecolor=nido,linewidth=0.04,dimen=outer]
(1.6,0.1)(6.4,3.6)
\psframe[linecolor=dred,linewidth=0.04,dimen=outer]
(-0.4,-0.3)(6.6,4.8)
\psframe[linecolor=dred,linewidth=0.04,dimen=outer]
(-0.6,-0.5)(6.8,5.8)
\psframe[linecolor=blue,linewidth=0.04,dimen=outer]
(1.4,-0.1)(9.2,3.8)
\psframe[linewidth=0.04,dimen=outer]
(-0.8,-0.7)(9.4,6.8)
\end{pspicture} 
}
\end{center}
\vskip -4mm
\caption{Classes of nonlinear operators.}
\label{fig:1}
\end{figure}
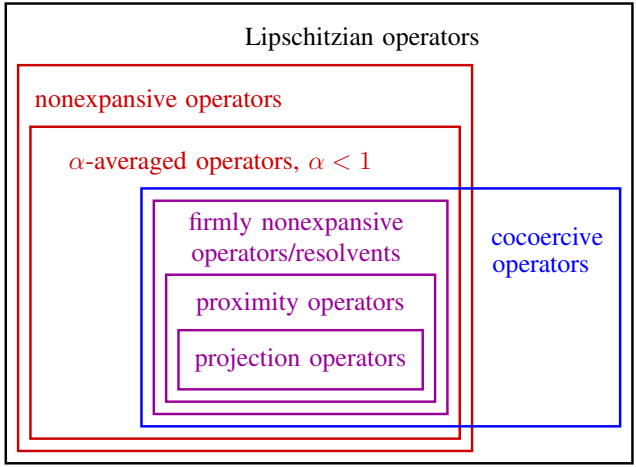

The relationships between the different types of nonlinear
operators discussed so far are depicted in Fig.~\ref{fig:1}.
The next propositions provide further connections between them.

\begin{proposition}
\label{p:19b}
Let $\delta\in\zeroun$, let $T\colon\HH\to\HH$ be 
$\delta$-Lipschitzian, and set $\alpha=(\delta+1)/2$. 
Then $T$ is $\alpha$-averaged. 
\end{proposition}

\begin{proposition}
\label{p:2004-5}
Let $T\colon\HH\to\HH$, let $\beta\in\RPP$, 
and let $\gamma\in\left]0,2\beta\right[$. Then $T$ is
$\beta$-cocoercive if and only if
$\Id-\gamma T$ is $\gamma/(2\beta)$-averaged. 
\end{proposition}

It follows from the Cauchy-Schwarz inequality that a 
$\beta$-cocoercive operator is $\beta^{-1}$-Lipschitzian. 
In the case of gradients of convex functions, the converse is 
also true.

\begin{proposition}[Baillon-Haddad]
\label{p:BH}
Let $f\colon\HH\to\RR$ be a differentiable convex function 
such that $\nabla f$ is $\beta^{-1}$-Lipschitzian for some
$\beta\in\RPP$. Then $\nabla f$ is $\beta$-cocoercive. 
\end{proposition}

We now describe operations that preserve averagedness and
cocoercivity.

\begin{proposition}
\label{p:av4}
Let $T\colon\HH\to\HH$, let $\alpha\in\zeroun$, and let
$\lambda\in\left]0,1/\alpha\right[$. Then $T$ is $\alpha$-averaged
if and only if $(1-\lambda)\Id+\lambda T$ is
$\lambda\alpha$-averaged.
\end{proposition}

\begin{proposition}
\label{p:av2}
For every $i\in\{1,\ldots,m\}$, let $\alpha_i\in\zeroun$, let
$\omega_i\in\rzeroun$, and let $T_i\colon\HH\to\HH$ be 
$\alpha_i$-averaged. Suppose that $\sum_{i=1}^m\omega_i=1$
and set $\alpha=\sum_{i=1}^m\omega_i\alpha_i$. 
Then $\sum_{i=1}^m\omega_iT_i$ is $\alpha$-averaged.
\end{proposition}

\begin{example}
\label{ex:2}
For every $i\in\{1,\ldots,m\}$, let $\omega_i\in\rzeroun$ and 
let $T_i\colon\HH\to\HH$ be firmly nonexpansive. Suppose that 
$\sum_{i=1}^m\omega_i=1$. Then $\sum_{i=1}^m\omega_iT_i$ is 
firmly nonexpansive.
\end{example}

\begin{proposition}
\label{p:roma}
For every $i\in\{1,\ldots,m\}$, let $\alpha_i\in\zeroun$ and let 
$T_i\colon\HH\to\HH$ be $\alpha_i$-averaged. Set 
\begin{equation}
\label{e:nyc2014-09-25}
T=T_1\circ\cdots\circ T_m\quad\text{and}\quad
\alpha=\dfrac{1}{1+\dfrac{1}
{\Sum_{i=1}^m\dfrac{\alpha_i}{1-\alpha_i}}}.
\end{equation}
Then $T$ is $\alpha$-averaged.
\end{proposition}

\begin{example}
\label{ex:roma}
Let $\alpha_1\in\zeroun$, 
let $\alpha_2\in\zeroun$, let $T_1\colon\HH\to\HH$ be 
$\alpha_1$-averaged, and let $T_2\colon\HH\to\HH$ be 
$\alpha_2$-averaged. Set
\begin{equation}
T=T_1\circ T_2\quad\text{and}\quad
\alpha=\frac{\alpha_1+\alpha_2-2\alpha_1\alpha_2}
{1-\alpha_1\alpha_2}.
\end{equation}
Then $T$ is $\alpha$-averaged.
\end{example}

\begin{proposition}[\cite{Davi17}]
\label{p:3magic}
Let $T_1\colon\HH\to\HH$ and $T_2\colon\HH\to\HH$ be firmly
nonexpansive, let $\alpha_3\in\zeroun$, and let 
$T_3\colon\HH\to\HH$ be $\alpha_3$-averaged. Set
$\alpha=1/(2-\alpha_3)$ and 
\begin{equation}
T=T_1\circ(T_2-\Id+T_3\circ T_2)+\Id-T_2.
\end{equation}
Then $T$ is $\alpha$-averaged.
\end{proposition}

\begin{proposition} 
\label{p:coco}
For every $k\in\{1,\ldots,q\}$, let 
$0\neq L_k\colon\HH\to\GG_k$ be linear,
let $\beta_k\in\RPP$, and let $T_k\colon\GG_k\to\GG_k$ be 
$\beta_k$-cocoercive. Set
\begin{equation}
\label{e:a12}
T=\sum_{k=1}^qL_k^*\circ T_k\circ L_k\quad\text{and}\quad
\beta=\dfrac{1}{{\Sum_{k=1}^q}\dfrac{\|L_k\|^2}{\beta_k}}.
\end{equation}
Then the following hold:
\begin{enumerate}
\item
\label{p:cocoi}
$T$ is $\beta$-cocoercive {\rm\cite{Livre1}}.
\item
\label{p:cocoii}
Suppose that $\sum_{k=1}^q\|L_k\|^2\leq 1$ and that the operators
$(T_k)_{1\leq k\leq q}$ are firmly nonexpansive. Then $T$ is
firmly nonexpansive {\rm\cite{Livre1}}.
\item
\label{p:cocoiii}
Suppose that $\sum_{k=1}^q\|L_k\|^2\leq 1$ and that 
$(T_k)_{1\leq k\leq q}$ are proximity operators. Then $T$ is a
proximity operator {\rm\cite{Comb18}}.
\end{enumerate}
\end{proposition}

\begin{remark}
\label{r:2018}
The statement of Proposition~\ref{p:coco}\ref{p:cocoiii} can be
made more precise \cite{Comb18}. To wit, for every
$k\in\{1,\ldots,q\}$, let $\omega_k\in\RPP$, let
$0\neq L_k\colon\HH\to\GG_k$ be linear, let
$g_k\in\Gamma_0(\GG_k)$, and let
$h_k\colon v\mapsto\inf_{w\in\GG_k}(g_k^*(w)+\|v-w\|^2/2)$ be
the \emph{Moreau envelope} of $g_k^*$. Then, if 
$\sum_{k=1}^q\omega_k\|L_k\|^2\leq 1$, we have 
\begin{multline}
\sum_{k=1}^q\omega_k\big(L_k^*\circ\prox_{g_k}\circ L_k\big)
=\prox_{f},\quad\text{where}\\
f=\Bigg(\sum_{k=1}^q\omega_kh_k
\circ L_k\Bigg)^*-\dfrac{\|\cdot\|^2}{2}.
\end{multline}
\end{remark}

\bigskip
Let $T\colon\HH\to\HH$ and let
\begin{equation}
\Fix T=\menge{x\in\HH}{Tx=x}
\end{equation}
be its set of \emph{fixed points}. If $T$ is a Banach contraction,
then it admits a unique fixed point. However, if $T$ is merely
nonexpansive, the situation is quite different. Indeed, a
nonexpansive operator may have no fixed point (take $T\colon
x\mapsto x+z$, with $z\neq 0$), exactly one (take $T=-\Id$), or
infinitely many (take $T=\Id$). Even those operators which are
firmly nonexpansive can fail to have fixed points.

\begin{example}
\label{ex:f8}
$T\colon\RR\to\RR\colon x\mapsto({x+\sqrt{x^2+4}})/2$
is firmly nonexpansive and $\Fix T=\emp$.
\end{example}

\begin{proposition}
\label{p:f1}
Let $T\colon\HH\to\HH$ be nonexpansive. 
Then $\Fix T$ is closed and convex. 
\end{proposition}

\begin{proposition}
\label{p:f2}
Let $(T_i)_{1\leq i\leq m}$ be nonexpansive operators from $\HH$ to
$\HH$, and let $(\omega_i)_{1\leq i\leq m}$ be real numbers in 
$\rzeroun$ such that $\sum_{i=1}^m\omega_i=1$. Suppose that
$\bigcap_{i=1}^m\Fix T_i\neq\emp$. Then 
$\Fix(\sum_{i=1}^m\omega_iT_i)=\bigcap_{i=1}^m\Fix T_i$. 
\end{proposition}

\begin{proposition}
\label{p:f3}
For every $i\in\{1,\ldots,m\}$, let $\alpha_i\in\zeroun$ and let 
$T_i\colon\HH\to\HH$ be $\alpha_i$-averaged. Suppose that
$\bigcap_{i=1}^m\Fix T_i\neq\emp$. Then 
$\Fix(T_1\circ\cdots\circ T_m)=\bigcap_{i=1}^m\Fix T_i$. 
\end{proposition}

\subsection{Monotone operators}
\label{sec:22}

\begin{figure}
\scalebox{0.63} 
{
\begin{pspicture}(-2.9,-3.5)(4.0,4.3)
\psline[linewidth=0.04cm,arrowsize=0.09cm 4.0,%
arrowlength=1.4,arrowinset=0.4]{->}(-0.0,-3.0)(-0.0,3.2)
\psline[linewidth=0.04cm,arrowsize=0.09cm 4.0,%
arrowlength=1.4,arrowinset=0.4]{->}(-2.7,-0.0)(3.0,-0.0)
\rput(3.3,0.0){\large $\HH$}
\rput(0.0,3.5){\large $\HH$}
\psplot[linewidth=0.05cm,linestyle=solid,algebraic,linecolor=blue]%
{-2.4}{-0.3}{(x+1)^2-1}
\psline[linewidth=0.05cm,linecolor=blue](-0.3,-0.2955)(-0.3,0.99)
\psplot[linewidth=0.05cm,linestyle=solid,algebraic,linecolor=blue]%
{0.3}{1.8}{2.5+0.4*(x-1.8)^3}
\psline[linewidth=0.05cm,linecolor=blue](2.1,-2.1)(2.1,2.0)
\psline[linewidth=0.05cm,linecolor=blue](2.075,-2.1)(2.7,-2.1)
\end{pspicture} 

\begin{pspicture}(-2.9,-3.5)(4.0,4.3)
\psline[linewidth=0.04cm,arrowsize=0.09cm 4.0,%
arrowlength=1.4,arrowinset=0.4]{->}(-0.0,-3.0)(-0.0,3.2)
\psline[linewidth=0.04cm,arrowsize=0.09cm 4.0,%
arrowlength=1.4,arrowinset=0.4]{->}(-2.7,-0.0)(3.0,-0.0)
\rput(3.3,0.0){\large $\HH$}
\rput(0.0,3.5){\large $\HH$}
\psplot[linewidth=0.05cm,linestyle=solid,algebraic,linecolor=nido]%
{-1.00}{-0.51}{sqrt(x+1)-1}
\psplot[linewidth=0.05cm,linestyle=solid,algebraic,linecolor=nido]%
{-1.00}{0.96}{-sqrt(x+1)-1}
\psline[linewidth=0.05cm,linecolor=nido](-0.2955,-0.3)(0.99,-0.3)
\psplot[linewidth=0.05cm,linestyle=solid,algebraic,linecolor=nido]%
{1.15}{2.5}{1.8-((2.5-x)/0.4)^(0.333333)}
\psline[linewidth=0.05cm,linecolor=nido](-2.1,2.1)(2.0,2.1)
\psline[linewidth=0.05cm,linecolor=nido](-2.08,2.1)(-2.08,2.7)
\end{pspicture} 
}
\caption{Left: Graph of a (nonmonotone) set-valued operator. 
Right: Graph of its inverse.}
\label{fig:5}
\end{figure}
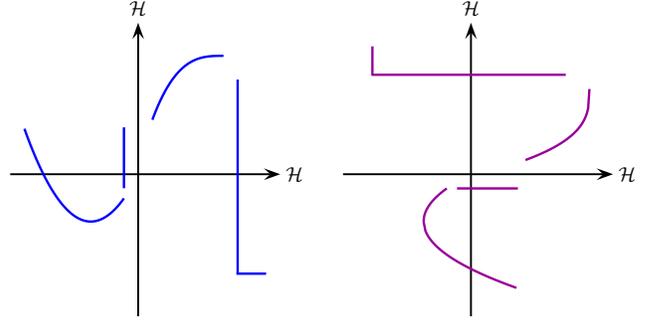

Let $A\colon\HH\to 2^{\HH}$ be a set-valued operator. Then $A$ is
described by its graph
\begin{equation}
\label{e:ge1}
\gra A=\menge{(x,u)\in\HH\times\HH}{u\in Ax},
\end{equation}
and its inverse $A^{-1}$, defined by the relation
\begin{equation}
\label{e:ge-1}
(\forall (x,u)\in\HH\times\HH)\quad x\in A^{-1}u\quad
\Leftrightarrow\quad u\in Ax,
\end{equation}
always exists (see Fig.~\ref{fig:5}). 
The operator $A$ is \emph{monotone} if 
\begin{multline}
\label{e:ge2}
\big(\forall (x,u)\in\gra A\big)\big(\forall (y,v)\in\gra A\big)\\
\scal{x-y}{u-v}\geq 0,
\end{multline}
in which case $A^{-1}$ is also monotone.

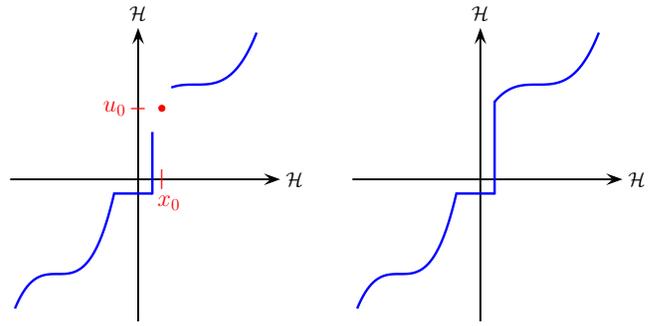
\begin{figure}
\scalebox{0.63} 
{
\begin{pspicture}(-2.9,-3.5)(4.2,4.3)
\psline[linewidth=0.04cm,arrowsize=0.09cm 4.0,%
arrowlength=1.4,arrowinset=0.4]{->}(-0.0,-3.0)(-0.0,3.2)
\psline[linewidth=0.04cm,arrowsize=0.09cm 4.0,%
arrowlength=1.4,arrowinset=0.4]{->}(-2.7,-0.0)(3.0,-0.0)
\rput(3.3,0.0){\large $\HH$}
\rput(0.0,3.5){\large $\HH$}
\rput(0.65,-0.5){\color{red}\Large $x_0$}
\rput(0.5,0.0){\color{red}\large $\boldsymbol{|}$}
\rput(-0.5,1.5){\color{red}\Large $u_0$}
\rput(-0.0,1.5){\color{red}\large $\boldsymbol{-}$}
\psplot[linewidth=0.05cm,linestyle=solid,algebraic,linecolor=blue]%
{-2.6}{-0.51}{(x+1.7)^3-2}
\psline[linewidth=0.05cm,linecolor=blue](-0.53,-0.3)(0.3,-0.3)
\psplot[linewidth=0.05cm,linestyle=solid,algebraic,linecolor=blue]%
{0.7}{2.5}{2.0+0.5*(x-1.2)^3}
\psline[linewidth=0.05cm,linecolor=blue](0.3,-0.32)(0.3,1.0)
\psdot[linewidth=0.04cm,linecolor=red](0.5,1.50)
\end{pspicture} 

\begin{pspicture}(-2.9,-3.5)(4.2,4.3)
\psline[linewidth=0.04cm,arrowsize=0.09cm 4.0,%
arrowlength=1.4,arrowinset=0.4]{->}(-0.0,-3.0)(-0.0,3.2)
\psline[linewidth=0.04cm,arrowsize=0.09cm 4.0,%
arrowlength=1.4,arrowinset=0.4]{->}(-2.7,-0.0)(3.0,-0.0)
\rput(3.3,0.0){\large $\HH$}
\rput(0.0,3.5){\large $\HH$}
\psplot[linewidth=0.05cm,linestyle=solid,algebraic,linecolor=blue]%
{-2.6}{-0.51}{(x+1.7)^3-2}
\psline[linewidth=0.05cm,linecolor=blue](-0.53,-0.3)(0.3,-0.3)
\psplot[linewidth=0.05cm,linestyle=solid,algebraic,linecolor=blue]%
{0.295}{2.5}{2.0+0.5*(x-1.2)^3}
\psline[linewidth=0.05cm,linecolor=blue](0.3,-0.32)(0.3,1.65)
\end{pspicture} 
}

\caption{Left: Graph of a monotone operator which is not maximally
monotone: we can add the point $(x_0,u_0)$ to its graph and still
get a monotone graph. Right: Graph of a maximally monotone
operator: adding any point to this graph destroys its
monotonicity.
}
\label{fig:6}
\end{figure}

\begin{example}
Let $f\colon\HH\to\RX$ be a proper function, let
$(x,u)\in\gra\partial f$, and let $(y,v)\in\gra\partial f$. Then 
\eqref{e:subdiff} yields
\begin{equation}
\begin{cases}
\scal{x-y}{u}+f(y)\geq f(x)\\
\scal{y-x}{v}+f(x)\geq f(y).
\end{cases}
\end{equation}
Adding these inequality yields $\scal{x-y}{u-v}\geq 0$, which 
shows that $\partial f$ is monotone.
\end{example}

A natural question is whether the operator obtained by adding a
point to the graph of a monotone operator $A\colon\HH\to 2^{\HH}$
is still monotone. If it is not, then $A$ is said to be 
\emph{maximally monotone}. Thus, $A$ is maximally monotone if, 
for every $(x,u)\in\HH\times\HH$, 
\begin{equation} 
\label{e:maxmon2}
(x,u)\in\gra A\;\Leftrightarrow\;(\forall (y,v)\in\gra A)\;\;
\scal{x-y}{u-v}\geq 0.
\end{equation}
These notions are illustrated in Fig.~\ref{fig:6}.
Let us provide some basic examples of maximally monotone
operators, starting with the subdifferential of \eqref{e:subdiff}
(see Fig.~\ref{fig:7}).

\begin{example}[Moreau]
\label{ex:mono1}
Let $f\in\Gamma_0(\HH)$. Then $\partial f$ is maximally monotone
and $(\partial f)^{-1}=\partial f^*$. 
\end{example}

\begin{example} 
\label{ex:mono5}
Let $T\colon\HH\to\HH$ be monotone and continuous. Then $T$ is
maximally monotone. In particular, if $T$ is cocoercive, it is 
maximally monotone. 
\end{example}

\begin{example}
\label{ex:mono3}
Let $T\colon\HH\to\HH$ be nonexpansive. Then $\Id-T$ is
maximally monotone. 
\end{example}

\begin{example}
\label{ex:mono4}
Let $T\colon\HH\to\HH$ be linear (hence continuous) and 
\emph{positive} in the sense
that $(\forall x\in\HH)$ $\scal{x}{Tx}\geq 0$.
Then $T$ is maximally monotone. In particular, if $T$ is
\emph{skew}, i.e., $T^*=-T$, then it is maximally monotone. 
\end{example}

Given $A\colon\HH\to 2^{\HH}$, the \emph{resolvent} of $A$ is 
the operator $J_A=(\Id+A)^{-1}$, that is,
\begin{equation}
\label{e:gm}
(\forall(x,p)\in\HH\times\HH)\quad p\in J_{\!A}x
\quad\Leftrightarrow\quad x-p\in Ap.
\end{equation}
In addition, the \emph{reflected resolvent} of $A$ is 
\begin{equation}
\label{e:RA}
R_A=2J_A-\Id.
\end{equation}
A profound result which connects monotonicity and
nonexpansiveness is Minty's theorem \cite{Mint62}. It implies that
if, $A\colon\HH\to 2^{\HH}$ is maximally monotone, then $J_A$ is
single-valued, defined everywhere on $\HH$, and firmly
nonexpansive.

\begin{theorem}[Minty]
\label{t:minty}
Let $T\colon\HH\to\HH$. 
Then $T$ is firmly nonexpansive if and only if it 
is the resolvent of a maximally monotone operator 
$A\colon\HH\to 2^{\HH}$.
\end{theorem}

\begin{example}
\label{ex:mono2}
Let $f\in\Gamma_0(\HH)$. Then $J_{\partial f}=\prox_f$.
\end{example}

Let $f$ and $g$ be functions in $\Gamma_0(\HH)$ which satisfy the
constraint qualification 
$\reli(\dom f)\cap\reli(\dom g)\neq\emp$. In view of 
Proposition~\ref{p:17}\ref{p:17ii} and Example~\ref{ex:mono1},
the minimizers of $f+g$ are precisely the solutions to the
inclusion $0\in Ax+Bx$ involving the maximally monotone operators
$A=\partial f$ and $B=\partial g$. Hence, 
it may seem that in minimization problems the 
theory of subdifferentials should suffice to analyze and solve 
problems without invoking general monotone operator theory. 
As discussed in \cite{Comb18}, this is not the case and monotone
operators play an indispensable role in various aspects of 
convex minimization. We give below an illustration of this fact in
the context of Proposition~\ref{p:17}.

\begin{example}[\cite{Siop11}]
\label{ex:m+s}
Given $f\in\Gamma_0(\HH)$, $g\in\Gamma_0(\GG)$, and a linear
operator $L\colon\HH\to\GG$, the objective is to 
\begin{equation}
\label{e:p}
\minimize{x\in\HH}{f(x)+g(Lx)}
\end{equation}
using $f$ and $g$ separately by means of their respective proximity
operators. To this end, let us bring into play the
Fenchel-Rockafellar dual problem 
\begin{equation}
\label{e:d}
\minimize{v\in\GG}{f^*(-L^*v)+g^*(v)}.
\end{equation}
We derive from \cite[Theorem~19.1]{Livre1} that, if
$(x,v)\in\HH\times\GG$ solves the inclusion
\begin{equation}
\label{e:d3}
\begin{bmatrix}
0\\
0
\end{bmatrix}
\in
\underbrace{\begin{bmatrix}
\partial f&0\\
0&\partial g^*\\
\end{bmatrix}}_{\text{subdifferential}}
\begin{bmatrix}
x\\
v
\end{bmatrix}
+
\underbrace{\begin{bmatrix}
0&L^*\\
-L&0\\
\end{bmatrix}}_{\text{skew}}
\begin{bmatrix}
x\\
v
\end{bmatrix},
\end{equation}
then $x$ solves \eqref{e:p} and $v$ solves \eqref{e:d}.
Now introduce the variable $\boldsymbol{z}=(x,v)$, the function 
$\Gamma_0(\HH\times\GG)\ni\boldsymbol{h}
\colon\boldsymbol{z}\mapsto f(x)+g^*(v)$, the
operator $\boldsymbol{A}=\partial\boldsymbol{h}$, and the skew 
operator $\boldsymbol{B}\colon\boldsymbol{z}\mapsto(L^*v,-Lx)$.
Then it follows from Examples~\ref{ex:mono1} and \ref{ex:mono4}
that \eqref{e:d3} can be written as the maximally monotone 
inclusion
$\boldsymbol{0}\in\boldsymbol{Az}+\boldsymbol{Bz}$, which does not
correspond to a minimization problem since $\boldsymbol{B}$ is not
a gradient \cite[Proposition~2.58]{Livre1}. As a result, genuine
monotone operator splitting methods were employed in \cite{Siop11} 
to solve \eqref{e:d3} and, thereby, \eqref{e:p} and \eqref{e:d}. 
Applications of this framework can be found in image restoration
\cite{Ocon14} and in empirical mode decomposition \cite{Pust14}.
\end{example}

\begin{example}
The primal-dual pair \eqref{e:p}--\eqref{e:d} can be exploited in
various ways; see for instance 
\cite{Cham11,Icip14,Svva10,Komo15}. 
A simple illustration is found in sparse signal recovery and
machine learning, where one often aims at solving \eqref{e:p} by
choosing $g$ to be a norm $|||\cdot|||$
\cite{Argy12,Bach12,Nume19,Dono03,McDo16}. Now let 
$|||\cdot|||_*\colon\GG\to\RR\colon v\mapsto
\sup_{|||y|||\leq 1}\scal{y}{v}$ be the dual norm and let 
$B_*=\menge{v\in\GG}{|||v|||_*\leq 1}$ be the associated unit 
ball. Then \eqref{e:d} is the constrained optimization problem 
\begin{equation} 
\minimize{v\in B_*}{f^*(-L^*v)}.
\end{equation} 
This dual formulation underlies several investigations, e.g.,
\cite{Elgh12,Ndia17}.
\end{example}

\section{Fixed point algorithms}
\label{sec:3}

We review the main fixed point construction algorithms.

\subsection{Basic iteration schemes}
First, we recall that finding the fixed point of a
Banach contraction is relatively straightforward via the standard
Banach-Picard iteration scheme \eqref{e:emile}.

\begin{theorem}[\cite{Livre1}]
\label{t:banach}
Let $\delta\in\zeroun$, let $T\colon\HH\to\HH$ be 
$\delta$-Lipschitzian, and let $x_0\in\HH$. Set
\begin{equation}
\label{e:banach}
(\forall n\in\NN)\quad x_{n+1}=Tx_n.
\end{equation}
Then $T$ has a unique fixed point $\overline{x}$ and
$x_n\to\overline{x}$. More precisely, $(\forall n\in\NN)$
$\|x_n-\overline{x}\|\leq\delta^n\|x_0-\overline{x}\|$.
\end{theorem}

If $T$ is merely nonexpansive (i.e., $\delta=1$) with 
$\Fix T\neq\emp$, Theorem~\ref{t:banach} fails. For instance, let
$T\neq\Id$ be a rotation in the Euclidean plane. Then it is
nonexpansive with $\Fix T=\{0\}$ but the sequence $(x_n)_{n\in\NN}$
constructed by the successive approximation process
\eqref{e:banach} does not converge. Such scenarios can be handled
via the following result. 

\begin{theorem}[\cite{Livre1}]
\label{t:1}
Let $\alpha\in\rzeroun$, let $T\colon\HH\to\HH$ be an
$\alpha$-averaged operator such that 
$\Fix T\neq\emp$, let $(\lambda_n)_{n\in\NN}$ be an
$\alpha$-relaxation sequence. Set
\begin{equation}
\label{e:groetsch1}
(\forall n\in\NN)\quad x_{n+1}=x_n+\lambda_n\big(Tx_n-x_n\big).
\end{equation}
Then $(x_n)_{n\in\NN}$ converges to a point in $\Fix T$. 
\end{theorem}

\begin{remark}
\label{r:1}
In connection with Theorems~\ref{t:banach} and \ref{t:1}, let us
make the following observations.
\begin{enumerate}
\item
\label{r:1i}
If $\alpha<1$ in Theorem~\ref{t:1}, choosing $\lambda_n=1$
in \eqref{e:groetsch1} (see 
Example~\ref{ex:relax}\ref{ex:relaxi}) yields \eqref{e:banach}.
\item
\label{r:1ii}
In contrast with Theorem~\ref{t:banach}, the convergence in 
Theorem~\ref{t:1} is not linear in general \cite{Baus09,Borw15}.
\item
\label{r:1iii}
When $\alpha=1$, \eqref{e:groetsch1} is known as the 
\emph{Krasnosel'ski\u\i-Mann iteration}.
\end{enumerate}
\end{remark}

Next, we present a more flexible fixed point theorem which involves
iteration-dependent composite averaged operators.

\begin{theorem}[\cite{Jmaa15}]
\label{t:2}
Let $\varepsilon\in\left]0,1/2\right[$ and let 
$x_0\in\HH$. For every $n\in\NN$, let 
$\alpha_{1,n}\in\left]0,1/(1+\varepsilon)\right]$, let
$\alpha_{2,n}\in\left]0,1/(1+\varepsilon)\right]$, 
let $T_{1,n}\colon\HH\to\HH$ be $\alpha_{1,n}$-averaged, and
let $T_{2,n}\colon\HH\to\HH$ be $\alpha_{2,n}$-averaged. 
In addition, for every $n\in\NN$, let
\begin{equation}
\lambda_n\in\big[\varepsilon,
{(1-\varepsilon)(1+\varepsilon\alpha_n)}/{\alpha_n}\big],
\end{equation}
where
$\alpha_n=({\alpha_{1,n}+\alpha_{2,n}-2\alpha_{1,n}\alpha_{2,n}})/
(1-\alpha_{1,n}\alpha_{2,n})$, 
and set 
\begin{equation}
\label{e:beforecluster}
x_{n+1}=x_n+\lambda_n\big(T_{1,n}(T_{2,n}x_n)-x_n\big).
\end{equation}
Suppose that 
$S=\bigcap_{n\in\NN}\Fix(T_{1,n}\circ T_{2,n})\neq\emp$.
Then the following hold:
\begin{enumerate}
\item
\label{t:2i}
$(\forall x\in S)$
$\sum_{n\in\NN}\|T_{2,n}x_n-x_n-T_{2,n}x+x\|^2<\pinf$.
\item
\label{t:2ii}
Suppose that a subsequence of $(x_n)_{n\in\NN}$ converges to a
point in $S$. Then $(x_n)_{n\in\NN}$ converges to a point in $S$.
\end{enumerate}
\end{theorem}

\begin{remark}
\label{r:mi}
The assumption in Theorem~\ref{t:2}\ref{t:2ii} holds in particular
when, for every $n\in\NN$, $T_{1,n}=T_1$ and $T_{2,n}=T_2$. 
\end{remark}

Below, we present a variant of Theorem~\ref{t:1} obtained by
considering the composition of $m$ operators. In the case of firmly
nonexpansive operators, this result is due to Martinet
\cite{Mart72}.

\begin{theorem}[\cite{Opti04}]
\label{t:5}
For every $i\in\{1,\ldots,m\}$, let $\alpha_i\in\zeroun$ and let
$T_i\colon\HH\to\HH$ be $\alpha_i$-averaged. Let $x_0\in\HH$,
suppose that $\Fix(T_1\circ\cdots\circ T_m)\neq\emp$, and
iterate
\begin{equation}
\label{e:2004}
\begin{array}{l}
\text{for}\;n=0,1,\ldots\\
\left\lfloor
\begin{array}{ll}
x_{mn+1}&\hskip -3mm={T_m}x_{mn}\\
x_{mn+2}&\hskip -3mm={T_{m-1}}x_{mn+1}\\
&\hskip -2mm\vdots\\
x_{mn+m-1}&\hskip -3mm={T_2}x_{mn+m-2}\\
x_{mn+m}&\hskip -3mm={T_1}x_{mn+m-1}.
\end{array}
\right.\\[2mm]
\end{array}
\end{equation}
Then $(x_{mn})_{n\in\NN}$ converges to a point $\overline{x}_1$ in 
$\Fix(T_1\circ\cdots\circ T_m)$. Now set 
$\overline{x}_m=T_m\overline{x}_1$,
$\overline{x}_{m-1}=T_{m-1}\overline{x}_m$,
\ldots, $\overline{x}_2=T_2\overline{x}_3$. 
Then, for every $i\in\{1,\ldots,m-1\}$,
$(x_{mn+i})_{n\in\NN}$ converges to $\overline{x}_{m+1-i}$.
\end{theorem}

\subsection{Algorithms for fixed point selection}

The algorithms discussed so far construct an unspecified fixed
point of a nonexpansive operator $T\colon\HH\to\HH$. In
some applications, one may be interested in finding a specific
fixed point, for instance one of minimum norm or, more generally,
one that minimizes some quadratic function \cite{Artz79,Sign03}.
One will find in \cite{Sign03} several algorithms to minimize
convex quadratic functions over fixed point sets, as well as 
signal recovery applications. Beyond quadratic selection, 
one may wish to minimize a strictly convex function
$g\in\Gamma_0(\HH)$ over the closed convex set (see
Proposition~\ref{p:f1}) $\Fix T$, i.e.,
\begin{equation}
\label{e:bi}
\minimize{x\in\Fix T} g(x).
\end{equation}
Instances of such formulations can be found in signal interpolation
\cite{Ono15} and machine learning \cite{Naka20}. Algorithms to 
solve \eqref{e:bi} have been proposed in 
\cite{Sico00,Hirs06,Yama01} under various hypotheses. Here is an
example.

\begin{proposition}[\cite{Yama01}]
\label{y:isao}
Let $T\colon\HH\to\HH$ be nonexpansive, let
$g\colon\HH\to\RR$ be strongly convex and differentiable with a
Lipschitzian gradient, let $x_0\in\HH$, and let
$(\alpha_n)_{n\in\NN}$ be a sequence
in $[0,1]$ such that $\alpha_n\to 0$,
$\sum_{n\in\NN}\alpha_n=\pinf$, and
$\sum_{n\in\NN}|\alpha_{n+1}-\alpha_n|<\pinf$. Suppose that
\eqref{e:bi} has a solution and iterate
\begin{equation}
(\forall n\in\NN)\quad x_{n+1}=Tx_n-\alpha_n\nabla g(Tx_n).
\end{equation}
Then $(x_n)_{n\in\NN}$ converges to the solution to \eqref{e:bi}.
\end{proposition}

\subsection{A fixed point method with block operator updates}
\label{sec:upda}

We turn our attention to a composite fixed point problem.

\begin{problem}
\label{prob:21}
Let $(\omega_i)_{1\leq i\leq m}$ be real numbers in $\rzeroun$ such
that $\sum_{i=1}^m\omega_i=1$. For every $i\in\{0,\ldots,m\}$, let
$T_i\colon\HH\to\HH$ be $\alpha_i$-averaged for some
$\alpha_i\in\zeroun$. The task is to find a fixed point of
$T_0\circ\sum_{i=1}^m\omega_iT_i$, assuming that such a point
exists.
\end{problem}

A simple strategy to solve Problem~\ref{prob:21} is to set
$R=\sum_{i=1}^m\omega_iT_i$, observe that $R$ is averaged by
Proposition~\ref{p:av2}, and then use Theorem~\ref{t:2} and
Remark~\ref{r:mi} to find a fixed point of $T_0\circ R$. This,
however, requires the activation of the $m$ operators 
$(T_i)_{1\leq i\leq m}$ to evaluate 
$R$ at each iteration, which is a significant
computational burden when $m$ is sizable. In the degenerate case
when the operators $(T_i)_{0\leq i\leq m}$ have common fixed
points, Problem~\ref{prob:21} amount to finding such a point (see
Propositions~\ref{p:f2} and \ref{p:f3}) and this can be done using
the strategies devised in \cite{Baus96,Nume06,Aiep96,Lopu97} which
require only the activation of blocks of operators at each
iteration. Such approaches fail in our more challenging setting,
which assumes only that
$\Fix(T_0\circ\sum_{i=1}^m\omega_iT_i)\neq\emp$. However, with 
a strategy based on tools from 
mean iteration theory \cite{Siop17}, it is possible
to devise an algorithm which operates by updating only a block of
operators $(T_i)_{i\in I_n}$ at iteration $n$.

\begin{theorem}[\cite{Upda20}]
\label{t:3}
Consider the setting of Problem~\ref{prob:21}. Let $M$ be a
strictly positive integer and let $(I_n)_{n\in\NN}$ be a sequence
of nonempty subsets of $\{1,\ldots,m\}$ such that 
\begin{equation}
\label{e:K}
(\forall n\in\NN)\quad\bigcup_{k=n}^{n+M-1}I_k=\{1,\ldots,m\}.
\end{equation}
Let $x_0\in\HH$, let $(t_{i,-1})_{1\leq i\leq m}\in\HH^m$, and
iterate
\begin{equation}
\label{e:a2}
\begin{array}{l}
\text{for}\;n=0,1,\ldots\\
\left\lfloor
\begin{array}{l}
\text{for every}\;i\in I_n\\
\left\lfloor
\begin{array}{l}
t_{i,n}=T_ix_n
\end{array}
\right.\\
\text{for every}\;i\in\{1,\ldots,m\}\smallsetminus I_n\\
\left\lfloor
\begin{array}{l}
t_{i,n}=t_{i,n-1}\\
\end{array}
\right.\\[1mm]
x_{n+1}=T_0\big(\sum_{i=1}^m\omega_it_{i,n}\big).
\end{array}
\right.\\
\end{array}
\end{equation}
Then the following hold:
\begin{enumerate}
\item 
\label{c:1i}
Let $x$ be a solution to Problem~\ref{prob:21} and let
$i\in\{1,\ldots,m\}$. Then $x_n-T_ix_n\to x-T_ix$.
\item 
\label{c:1ii}
$(x_n)_{n\in\NN}$ converges to a solution to
Problem~\ref{prob:21}.
\item
\label{c:1iv}
Suppose that, for some $i\in\{0,\ldots,m\}$, $T_i$ is a Banach
contraction. Then $(x_n)_{n\in\NN}$ converges linearly to the
unique solution to Problem~\ref{prob:21}.
\end{enumerate}
\end{theorem}

At iteration $n$, $I_n$ is the set of indices of operators to be
activated. The remaining operators are not used and their most
recent evaluations are recycled to form the update $x_{n+1}$.
Condition \eqref{e:K} imposes the mild requirement that each
operator in $(T_i)_{1\leq i\leq m}$ be evaluated at least once over
the course of any $M$ consecutive iterations. The choice of $M$ is
left to the user.

\subsection{Perturbed fixed point methods}

For various modeling or computational reasons, exact evaluations 
of the operators in fixed point algorithms may not be possible.
Such perturbations can be modeled by deterministic additive errors
\cite{Opti04,Lema96,Mart72} but also by stochastic ones
\cite{Comb15,Ermo69}. Here is a stochastically perturbed version of
Theorem~\ref{t:1}, which is a straightforward variant of
\cite[Corollary~2.7]{Comb15}.

\begin{theorem}
\label{t:1stoch}
Let $\alpha\in\rzeroun$, let $T\colon\HH\to\HH$ be an
$\alpha$-averaged operator such that 
$\Fix T\neq\emp$, and let $(\lambda_n)_{n\in\NN}$ be an
$\alpha$-relaxation sequence.
Let $x_0$ and $(e_n)_{n\in\NN}$ 
be $\HH$-valued random variables. Set
\begin{equation}
\label{e:groetsch2}
(\forall n\in\NN)\quad x_{n+1}=x_n+\lambda_n\big(Tx_n+e_n-x_n\big).
\end{equation}
Suppose that 
$\sum_{n\in\NN}\lambda_n\sqrt{\EC{\|e_n\|^2}{\XX_n}}<\pinf$ $\as$,
where $\XX_n$ is the $\sigma$-algebra generated by 
$(x_0,\ldots,x_n)$. Then $(x_n)_{n\in\NN}$ converges $\as$ to a 
$(\Fix{T})$-valued random variable.
\end{theorem}

\subsection{Random block-coordinate fixed point methods}
\label{sec:87}
We have seen in Section~\ref{sec:upda} that the computational cost
per iteration could be reduced in certain fixed point algorithms by
updating only some of the operators involved in the model. In this 
section, we present another approach to reduce the iteration cost
by considering scenarios in which the underlying Euclidean space
$\HHH$ is decomposable in $m$ factors 
$\HHH=\HH_1\times\cdots\times\HH_m$.
In the spirit of the Gauss-Seidel algorithm, one can explore the
possibility of activating only some of the coordinates of certain
operators at each iteration of a fixed point method. The potential
advantages of such a procedure are a reduced computational cost per
iteration, reduced memory requirements, and an increased
implementation flexibility.

In the product space $\HHH$, consider the basic update process 
\begin{equation}
\label{e:basiciter}
\boldsymbol{x}_{n+1}=
\boldsymbol{T}_{\!n}\boldsymbol{x}_n,
\end{equation}
under the assumption that the operator $\boldsymbol{T}_{\!n}$ is
decomposable explicitly as
\begin{equation}
\boldsymbol{T}_{\!n}\colon\HHH\to\HHH\colon
\boldsymbol{x}\mapsto (T_{1,n}\boldsymbol{x},\ldots,
T_{m,n}\boldsymbol{x}),
\end{equation}
with $T_{i,n}\colon\HHH\to\HH_i$.
Updating only some coordinates is performed by
modifying iteration \eqref{e:basiciter} as 
\begin{equation}
\label{e:stocbasiciter}
(\forall i\in\{1,\ldots,m\})\quad
x_{i,n+1}=x_{i,n}+\varepsilon_{i,n} \big(T_{i,n}x_n-x_{i,n}\big),
\end{equation}
where $\varepsilon_{i,n}\in\{0,1\}$ signals the 
activation of the $i$-th coordinate of $\boldsymbol{x}_n$. If
$\varepsilon_{i,n}=1$, the $i$-th component is updated whereas, if
$\varepsilon_{i,n}=0$, it is unchanged. The main difficulty facing
such an approach is that the nonexpansiveness property of an
operator is usually destroyed by coordinate sampling. To remove
this roadblock, a possibility is to make the activation
variables random, which results in a stochastic algorithm for which
almost sure convergence holds \cite{Comb15,Iutz13}. 

\begin{theorem}[\cite{Comb15}]
Let $\alpha\in\rzeroun$, let $\epsilon\in\left]0,1/2\right[$,
and let $\boldsymbol{T}\colon\HHH\to\HHH\colon
\boldsymbol{x}\mapsto(T_{\!i}\,
\boldsymbol{x})_{1\leq i\leq m}$ be an
$\alpha$-averaged operator where
${T}_i\colon\HHH\to\HH_i$.
Let $(\lambda_n)_{n\in\NN}$ be in
$[\epsilon,\alpha^{-1}-\epsilon]$, 
set $D=\{0,1\}^m\smallsetminus \{\boldsymbol{0}\}$,
let $\boldsymbol{x}_0$ be an $\HHH$-valued random variable, and let
$(\boldsymbol{\varepsilon}_n)_{n\in\NN}$ be identically distributed
$D$-valued random variables. Iterate
\begin{equation}
\label{e:2014-02-09}
\begin{array}{l}
\text{for}\;n=0,1,\ldots\\
\left\lfloor
\begin{array}{l}
\text{for}\;i=1,\ldots,m\\
\left\lfloor
\begin{array}{l}
x_{i,n+1}=x_{i,n}+\varepsilon_{i,n}\lambda_n\big(
T_i\boldsymbol{x}_n-x_{i,n}\big).
\end{array} 
\right.
\end{array} 
\right.
\end{array} 
\end{equation}
In addition, assume that the following hold:
\begin{enumerate}
\item
$\Fix\boldsymbol{T}\neq\emp$.
\item
For every $n\in\NN$, $\boldsymbol{\varepsilon}_n$ and
$(\boldsymbol{x}_0,\ldots,\boldsymbol{x}_n)$ are mutually
independent.
\item
$(\forall i\in\{1,\ldots,m\})$ 
$\PP[\varepsilon_{i,0}=1]>0$.
\end{enumerate}
Then
$(\boldsymbol{x}_n)_{n\in\NN}$ converges $\as$ to a
$\Fix\boldsymbol{T}$-valued random variable.
\end{theorem}

Further results in this vein for iterations involving nonstationary
compositions of averaged operators can be found in \cite{Comb15}.
Mean square convergence results are also available under additional
assumptions on the operators $(\boldsymbol{T}_{\!n})_{n\in\NN}$
\cite{Comb19}. 

\section{Fixed point modeling of monotone inclusions}
\label{sec:4}

\subsection{Splitting sums of monotone operators}
\label{sec:A} 

Our first basic model is that of finding a zero of the sum of two
monotone operators. It will be seen to
be central in understanding and solving data science problems in
optimization form (see also Example~\ref{ex:m+s} for a
special case) and beyond.

\begin{problem}
\label{prob:7}
Let $A\colon\HH\to 2^{\HH}$ and $B\colon\HH\to 2^{\HH}$ be
maximally monotone operators. The task is to 
\begin{equation}
\label{e:r29p}
\text{find}\;\;x\in\HH\;\;\text{such that}\;\;0\in Ax+Bx,
\end{equation}
under the assumption that a solution exists. 
\end{problem}

A classical method for solving Problem~\ref{prob:7} is the 
\emph{Douglas-Rachford} algorithm, which was first proposed in
\cite{Lion79} (see also \cite{Ecks92}; the following relaxed 
version is from \cite{Eoop01}).

\begin{proposition}[Douglas-Rachford splitting]
\label{p:DR}
Let $(\lambda_n)_{n\in\NN}$ be a $1/2$-relaxation sequence,
let $\gamma\in\RPP$, and let $y_0\in\HH$. Iterate
\begin{equation}
\label{e:dr4}
\begin{array}{l}
\text{for}\;n=0,1,\ldots\\
\left\lfloor
\begin{array}{l}
x_n=J_{\gamma B}y_n\\
z_n=J_{\gamma A}(2x_n-y_n)\\
y_{n+1}=y_n+\lambda_n(z_n-x_n).
\end{array}
\right.\\[2mm]
\end{array}
\end{equation}
Then $(x_n)_{n\in\NN}$ converges to a
solution to Problem~\ref{prob:7}.
\end{proposition}

The Douglas-Rachford algorithm requires the ability to evaluate two
resolvents at each iteration. However, if one of the operators is
single-valued and Lipschitzian, it is possible apply it explicitly,
hence requiring only one resolvent evaluation per iteration. The
resulting algorithm, proposed by Tseng \cite{Tsen00}, is often
called the \emph{forward-backward-forward} splitting algorithm
since it involves two explicit (forward) steps using $B$ and one
implicit (backward) step using $A$.

\begin{proposition}[Tseng splitting]
\label{p:khk09}
In Problem~\ref{prob:7}, assume that $B$ is $\delta$-Lipschitzian 
for some $\delta\in\RPP$. Let $x_0\in\HH$, let 
$\varepsilon\in\left]0,1/(\delta+1)\right[$,
let $(\gamma_n)_{n\in\NN}$ be in 
$[\varepsilon,(1-\varepsilon)/\delta]$, and iterate
\begin{equation}
\label{e:rio1010}
\begin{array}{l}
\text{for}\;n=0,1,\ldots\\
\left\lfloor
\begin{array}{l}
y_n=x_n-\gamma_n Bx_n\\
z_n=J_{\gamma_n A}y_n\\
r_n=z_n-\gamma_n Bz_n\\
x_{n+1}=x_n-y_n+r_n.
\end{array}
\right.\\[2mm]
\end{array}
\end{equation}
Then $(x_n)_{n\in\NN}$ converges to a solution to 
Problem~\ref{prob:7}.
\end{proposition}

As noted in Section~\ref{sec:23}, if $B$ is cocoercive, then it is
Lipschitzian, and Proposition~\ref{p:khk09} is applicable. However,
in this case it is possible to devise an algorithm which requires
only one application of $B$ per iteration, as opposed
to two in \eqref{e:rio1010}. 
To see this, let $\gamma_n\in\left]0,2\beta\right[$ and $x\in\HH$. 
Then it follows at once from \eqref{e:gm} that $x$ solves
Problem~\ref{prob:7}
$\Leftrightarrow$ $-\gamma_nBx\in\gamma_n Ax$
$\Leftrightarrow$ $(x-\gamma_nBx)-x\in\gamma_n Ax$
$\Leftrightarrow$ $x=J_{\gamma_n A}(x-\gamma_nBx)$
$\Leftrightarrow$ $x\in\Fix(T_{1,n}\circ T_{2,n})$, where
$T_{1,n}=J_{\gamma_n A}$ and $T_{2,n}=\Id-\gamma_nB$. As seen in 
Theorem~\ref{t:minty}, $T_{1,n}$ is $1/2$-averaged. On the other
hand, we derive from Proposition~\ref{p:2004-5} that, if 
$\alpha_{2,n}=\gamma_n/(2\beta)$, then $T_{2,n}$ is 
$\alpha_{2,n}$-averaged. With these considerations, we invoke
Theorem~\ref{t:2} to obtain the following algorithm, which goes
back to \cite{Merc79}.

\begin{proposition}[forward-backward splitting \cite{Jmaa15}]
\label{p:fb13}
Suppose that, in Problem~\ref{prob:7}, $B$ is $\beta$-cocoercive
for some $\beta\in\RPP$. Let
$\varepsilon\in\left]0,\min\{1/2,\beta\}\right[$,
let $x_0\in\HH$, and let $(\gamma_n)_{n\in\NN}$ be in
$\left[\varepsilon,2\beta/(1+\varepsilon)\right]$. Let 
\begin{equation}
(\forall n\in\NN)\quad
\lambda_n\in\big[\varepsilon,(1-\varepsilon)
\big(2+\varepsilon-{\gamma_n}/(2\beta)\big)\big].
\end{equation}
Iterate
\begin{equation}
\label{e:FB1}
\begin{array}{l}
\text{for}\;n=0,1,\ldots\\
\left\lfloor
\begin{array}{l}
u_n=x_n-\gamma_n Bx_n\\
x_{n+1}=x_n+\lambda_n\big(J_{\gamma_n A}u_n-x_n\big).
\end{array}
\right.\\[2mm]
\end{array}
\end{equation}
Then $(x_n)_{n\in\NN}$ converges to a solution to
Problem~\ref{prob:7}.
\end{proposition}

We now turn our attention to a more structured version of 
Problem~\ref{prob:7}, which includes an additional Lipschitzian
monotone operator. 

\begin{problem}
\label{prob:8}
Let $A\colon\HH\to 2^{\HH}$ and $B\colon\HH\to 2^{\HH}$ be
maximally monotone operators, let $\delta\in\RPP$, and let 
$C\colon\HH\to\HH$ be monotone and $\delta$-Lipschitzian. The task
is to 
\begin{equation}
\label{e:r29q}
\text{find}\;\;x\in\HH\;\;\text{such that}\;\;0\in Ax+Bx+Cx,
\end{equation}
under the assumption that a solution exists. 
\end{problem}

The following approach provides also a dual solution.

\begin{proposition}[splitting three operators I \cite{Svva12}]
\label{p:71}
Consider Problem~\ref{prob:8} and let
$\varepsilon\in\left]0,1/(2+\delta)\right[$.
Let $(\gamma_n)_{n\in\NN}$ be in 
$\left[\varepsilon,(1-\varepsilon)/(1+\delta)\right]$,
let $x_0\in\HH$, and let $u_0\in\HH$. Iterate
\begin{equation}
\label{e:3ops1}
\begin{array}{l}
\text{for}\;n=0,1,\ldots\\
\left\lfloor
\begin{array}{l}
y_n=x_n-\gamma_n(Cx_n+u_n)\\
p_n=J_{\gamma_n A}\,y_n\\
q_n=u_n+\gamma_n\big(x_n-J_{B/\gamma_n}(u_n/\gamma_n+x_n)\big)\\
x_{n+1}=x_n-y_n+p_n-\gamma_n(Cp_n+q_n)\\
u_{n+1}=q_n+\gamma_n(p_n-x_n).
\end{array}
\right.\\[2mm]
\end{array}
\end{equation}
Then $(x_n)_{n\in\NN}$ converges to a solution to 
Problem~\ref{prob:8} and $(u_n)_{n\in\NN}$ converges to a 
solution $u$ to the dual problem, i.e.,
$0\in-(A+C)^{-1}(-u)+B^{-1}u$.
\end{proposition}

When $C$ is $\beta$-cocoercive in Problem~\ref{prob:8}, 
we can take $\delta=1/\beta$. 
In this setting, an alternative algorithm 
is obtained as follows. Let us fix $\gamma\in\RPP$ and define 
\begin{align}
\label{e:T3ops}
T=J_{\gamma A}\circ \big(2J_{\gamma B}-\Id-\gamma 
C\circ J_{\gamma B}\big)+\Id-J_{\gamma B}.
\end{align}
By setting $T_1=J_{\gamma A}$, $T_2=J_{\gamma B}$, and
$T_3=\Id-\gamma C$ in Proposition~\ref{p:3magic}, we deduce
from Proposition~\ref{p:2004-5} that, if 
$\gamma\in\left]0,2\beta\right[$ and
$\alpha=2\beta/(4\beta-\gamma)$, then $T$ is
$\alpha$-averaged. Now take $y\in\HH$ and set
$x=J_{\gamma B}y$, hence $y-x\in\gamma Bx$ by \eqref{e:gm}. 
Then $y\in\Fix T$ $\Leftrightarrow$ 
$J_{\gamma A}(2x-y-\gamma Cx)+y-x=y$ $\Leftrightarrow$ 
$J_{\gamma A}(2x-y-\gamma Cx)=x$
$\Leftrightarrow$ $x-y-\gamma Cx\in\gamma Ax$ by \eqref{e:gm}. 
Thus, $0=(x-y)+(y-x)\in\gamma(Ax+Bx+Cx)$, which shows that $x$
solves Problem~\ref{prob:8}. Altogether, since $y$ can be
constructed via Theorem~\ref{t:1}, we obtain the following
convergence result.

\begin{proposition}[splitting three operators II \cite{Davi17}]
\label{p:72}
In Problem~\ref{prob:8}, assume that 
$C$ is $\beta$-cocoercive for some $\beta\in\RPP$. 
Let $\gamma\in\left]0,2\beta\right[$ and set 
$\alpha=2\beta/(4\beta-\gamma)$.
Furthermore, let $(\lambda_n)_{n\in\NN}$ be an
$\alpha$-relaxation sequence and let $y_0\in\HH$. 
Iterate
\begin{equation}
\label{e:3ops2}
\begin{array}{l}
\text{for}\;n=0,1,\ldots\\
\left\lfloor
\begin{array}{l}
x_n=J_{\gamma B}\,y_n\\
r_n=y_n+\gamma C x_n\\
z_n=J_{\gamma A}(2x_n-r_n)\\
y_{n+1}=y_n+\lambda_n(z_n-x_n).
\end{array}
\right.\\[2mm]
\end{array}
\end{equation}
Then $(x_n)_{n\in\NN}$ converges to a solution to 
Problem~\ref{prob:8}.
\end{proposition}

\begin{remark}
\label{re:3op}\
\begin{enumerate}
\item 
Work closely related to Proposition~\ref{p:72} can be found in 
\cite{Brice15,Bric18,Ragu13}. See also \cite{Ragu19}, which
provides further developments and a discussion of
\cite{Brice15,Davi17,Ragu13}.
\item 
Unlike algorithm \eqref{e:3ops1}, \eqref{e:3ops2} imposes constant
proximal parameters and requires the cocoercivity of $C$, but it
involves only one application of $C$ per iteration. An extension
of \eqref{e:3ops2} appears in \cite{Yanm18} in the context of
minimization problems. 
\end{enumerate}
\end{remark}

\subsection{Splitting sums of composite monotone operators}
\label{se:splitmon}
The monotone inclusion problems of Section~\ref{sec:A} are
instantiations of the following formulation, which involves an
arbitrary number of maximally monotone operators and compositions
with linear operators. 

\begin{problem}
\label{prob:88}
Let $\delta\in\RPP$ and let $A\colon\HH\to 2^\HH$ 
be maximally monotone.
For every $k\in\{1,\ldots,q\}$,
let $B_k\colon\GG_k\to 2^{\GG_k}$ be
maximally monotone, let $0\neq L_k\colon\HH\to\GG_k$ be linear,
and let $C_k\colon\GG_k\to\GG_k$ be monotone and
$\delta$-Lipschitzian. The task is to 
\begin{multline}
\label{e:r29qq}
\text{find}\;\;x\in\HH\;\;\text{such that}\\
0\in Ax+\sum_{k=1}^q L_k^*\big((B_k+C_k)(L_kx)\big),
\end{multline}
under the assumption that a solution exists. 
\end{problem}

In the context of Problem~\ref{prob:88}, the principle of a
splitting algorithm is to involve all the operators individually.
In the case of a set-valued operator $A$ or $B_k$, this means using
the associated resolvent, whereas in the case of a single-valued
operator $C_k$ or $L_k$, a direct application can be considered.
An immediate difficulty one faces with \eqref{e:r29qq} is that
it involves many set-valued operators. However, since inclusion is
a binary relation, for reasons discussed in \cite{Siop11,Play13}
and analyzed in more depth in \cite{Ryue20}, it is not possible to
deal with more than two such operators. To circumvent this
fundamental limitation, a strategy is to rephrase
Problem~\ref{prob:88} as a problem involving at most two set-valued
operators in a larger space. This strategy finds its root in convex
feasibility problems \cite{Pier84} and it was first adapted to the
problem of finding a zero of the sum of $m$ operators in
\cite{Gols87,Spin83}. In \cite{Siop11}, it was used to
deal with the presence of linear operators (see in particular
Example~\ref{ex:m+s}), with further developments in
\cite{Botr13,Botr14,Icip14,Svva12,Bang13}. 
In the same spirit, let us reformulate Problem~\ref{prob:88} 
by introducing
\begin{equation}
\label{e:d9}
\begin{cases}
\boldsymbol{L}\colon\HH\to\GGG\colon
x\mapsto (L_1x,\ldots,L_qx)\\
\boldsymbol{B}\colon\GGG\to 2^{\GGG}\colon
(y_k)_{1\leq k\leq q}\mapsto\cart_{\!k=1}^{\!q}B_ky_k\\
\boldsymbol{C}\colon\GGG\to\GGG\colon
(y_k)_{1\leq k\leq q}\mapsto (C_ky_k)_{1\leq k\leq q}\\
\boldsymbol{V}=\ran\boldsymbol{L}.
\end{cases}
\end{equation}
Note that $\boldsymbol{L}$ is linear, $\boldsymbol{B}$ is maximally
monotone, and $\boldsymbol{C}$ is monotone and
$\delta$-Lipschitzian. In addition, the inclusion \eqref{e:r29qq}
can be rewritten more concisely as
\begin{equation}
\label{e:r29qq2}
\text{find}\;\;x\in\HH\;\;\text{such that}\;\;
0\in Ax+\boldsymbol{L}^*\big((\boldsymbol{B}+\boldsymbol{C})
(\boldsymbol{L}x)\big).
\end{equation}
In particular, suppose that $A=0$. Then, upon setting 
$\boldsymbol{y}=\boldsymbol{L}x\in\boldsymbol{V}$, we obtain the
existence of a point 
$\boldsymbol{u}\in (\boldsymbol{B}+\boldsymbol{C})\boldsymbol{y}$
in $\ker\boldsymbol{L}^*=\boldsymbol{V}^\perp$. In other words, 
\begin{equation} 
\label{e:reforprob88}
\boldsymbol{0}\in
N_{\boldsymbol{V}}\boldsymbol{y}+\boldsymbol{B}\boldsymbol{y}+
\boldsymbol{C}\boldsymbol{y}.
\end{equation}
Solving this inclusion is equivalent to solving a problem similar
to Problem~\ref{prob:8}, formulated in $\GGG$. Thus, applying
Proposition~\ref{p:72} to \eqref{e:reforprob88} leads to the
following result.

\begin{proposition}
\label{p:71b}
In Problem~\ref{prob:88}, suppose that $A=0$, that the operators
$(C_k)_{1\leq k\leq q}$ are $\beta$-cocoercive for some
$\beta\in\RPP$, and that $Q=\sum_{k=1}^qL_k^*\circ L_k$ is
invertible. Let $\gamma\in\left]0,2\beta\right[$, set 
$\alpha=2\beta/(4\beta-\gamma)$, and let
$(\lambda_n)_{n\in\NN}$ be an $\alpha$-relaxation sequence.
Further, let $\boldsymbol{y}_0\in\GGG$, set 
$s_0=Q^{-1}\Big(\sum_{k=1}^qL_k^*y_{0,k}\Big)$, 
and iterate
\begin{equation}
\label{e:3opsbis}
\begin{array}{l}
\text{for}\;n=0,1,\ldots\\
\left\lfloor
\begin{array}{l}
\text{for}\;k=1,\ldots,q\\
\left\lfloor
\begin{array}{l}
p_{n,k}=J_{\gamma B_k}\,y_{n,k}\\
\end{array}
\right.\\[2mm]
x_n=Q^{-1}\big(\sum_{k=1}^qL_k^* p_{n,k}\big)\\[1mm]
c_n=Q^{-1}\big(\sum_{k=1}^qL_k^* C_k p_{n,k}\big)\\
z_n=x_n-s_n-\gamma c_n\\
\text{for}\;k=1,\ldots,q\\
\left\lfloor
\begin{array}{l}
y_{n+1,k}=y_{n,k}+\lambda_n (x_n+z_n-p_{n,k})
\end{array}
\right.\\[2mm]
s_{n+1}=s_n+\lambda_n z_n.
\end{array}
\right.\\[2mm]
\end{array}
\end{equation}
Then $(x_n)_{n\in\NN}$ converges to a solution to 
\eqref{e:r29qq}.
\end{proposition}

A strategy for handling Problem~\ref{prob:88} in its general
setting consists of introducing an auxiliary variable
$\boldsymbol{v}\in\boldsymbol{B}(\boldsymbol{L}x)$ in
\eqref{e:r29qq2}, which can then be rewritten as 
\begin{equation}
\label{e:monskew}
\begin{cases}
0\in Ax+\boldsymbol{L}^*\boldsymbol{v}
+\boldsymbol{L}^*\big(\boldsymbol{C}(\boldsymbol{L}x)\big)\\
\boldsymbol{0}\in -\boldsymbol{L} x +\boldsymbol{B}^{-1} 
\boldsymbol{v}.
\end{cases}
\end{equation}
This results in an instantiation of Problem~\ref{prob:7} 
in $\KKK=\HH\times\GGG$ involving the maximally monotone operators
\begin{equation}
\begin{cases}
\boldsymbol{\mathcal{A}}_1\colon\KKK\to 2^{\KKK}\colon
(x,\boldsymbol{v})&\mapsto 
\begin{bmatrix}
A & 0\\
0 & \boldsymbol{B}^{-1}
\end{bmatrix}
\begin{bmatrix}
x\\
\boldsymbol{v}
\end{bmatrix}\\[4mm]
\boldsymbol{\mathcal{B}}_1\colon\KKK\to\KKK
\colon (x, \boldsymbol{v})&\mapsto\begin{bmatrix}
\boldsymbol{L}^*\circ\boldsymbol{C}\circ\boldsymbol{L} 
& \boldsymbol{L}^*\\
-\boldsymbol{L} & 0 \\
\end{bmatrix}
\begin{bmatrix}
x\\
\boldsymbol{v}
\end{bmatrix}.
\end{cases}
\end{equation}
We observe that, in $\KKK$,
$\boldsymbol{\mathcal{B}}_1$ is Lipschitzian with constant
$\chi=\|\boldsymbol{L}\|(1+\delta\|
\boldsymbol{L}\|)$. By applying Proposition~\ref{p:khk09} to
\eqref{e:monskew}, we obtain the following algorithm.

\begin{proposition}[\cite{Svva12}]
\label{prop:MLFB}
Consider Problem~\ref{prob:88}. Set 
\begin{equation}
\chi=\sqrt{\textstyle{\sum_{k=1}^q}\|L_k\|^2}
\Big(1+\delta\sqrt{\textstyle{\sum_{k=1}^q}\|L_k\|^2}\Big).
\end{equation}
Let $x_0\in\HH$, let $\boldsymbol{v}_0\in\GGG$,
let $\varepsilon\in\left]0,1/(\chi+1)\right[$,
let $(\gamma_n)_{n\in\NN}$ be in 
$[\varepsilon,(1-\varepsilon)/\chi]$, and iterate
\begin{equation}
\label{e:rio1010b}
\begin{array}{l}
\text{for}\;n=0,1,\ldots\\
\left\lfloor
\begin{array}{l}
u_n=x_n-\gamma_n 
\sum_{k=1}^q L_k^* (C_k(L_kx_n)+v_{n,k})\\
p_n=J_{\gamma_n A} u_n\\
\text{for}\;k=1,\ldots,q\\
\left\lfloor
\begin{array}{l}
y_{n,k}=v_{n,k}+\gamma_nL_k x_n\\
z_{n,k}=y_{n,k}-\gamma_n J_{\gamma^{-1} B_k}
\big({y_{n,k}}/{\gamma_n}\big)\\ 
s_{n,k}=z_{n,k}+\gamma_nL_k p_n\\
v_{n+1,k}=v_{n,k}-y_{n,k}+s_{n,k}\\
\end{array}
\right.\\[2mm]
r_n=p_n-\gamma_n\sum_{k=1}^q L_k^* 
(C_k(L_kp_n)+z_{n,k})\\ 
x_{n+1}=x_n-u_n+r_n.
\end{array}
\right.\\[2mm]
\end{array}
\end{equation}
Then $(x_n)_{n\in\NN}$ converges to a solution to 
Problem~\ref{prob:88}.
\end{proposition}

An alternative approach consists of reformulating \eqref{e:monskew}
in the form of Problem~\ref{prob:7} with the maximally monotone 
operators 
\begin{equation}
\begin{cases}
\boldsymbol{\mathcal{A}}_2
\colon\KKK\to 2^{\KKK}\colon (x,\boldsymbol{v})
&\mapsto
\begin{bmatrix}
A & \boldsymbol{L}^*\\
-\boldsymbol{L} & \boldsymbol{B}^{-1}\\
\end{bmatrix}
\begin{bmatrix}
x\\
\boldsymbol{v}
\end{bmatrix}\\[4mm]
\boldsymbol{\mathcal{B}}_2\colon\KKK\to\KKK\colon
(x,\boldsymbol{v})
&\mapsto
\begin{bmatrix}
\boldsymbol{L}^*\circ\boldsymbol{C}\circ\boldsymbol{L} & 0\\
0 & 0
\end{bmatrix} 
\begin{bmatrix}
x\\
\boldsymbol{v}
\end{bmatrix}.
\end{cases}
\end{equation}
Instead of working directly with these operators, it may be
judicious to use preconditioned versions
$\boldsymbol{V}\circ\boldsymbol{\mathcal{A}}_2$ and
$\boldsymbol{V}\circ\boldsymbol{\mathcal{B}}_2$, where
$\boldsymbol{V}\colon\KKK\to\KKK$ is a self-adjoint strictly
positive linear operator. If $\KKK$ is renormed with 
\begin{equation}
\|\cdot\|_{\boldsymbol{V}}\colon (x,\boldsymbol{v})\mapsto 
\sqrt{\scal{(x,\boldsymbol{v})}{\boldsymbol{V}^{-1}
(x,\boldsymbol{v})}},
\end{equation}
then 
$\boldsymbol{V}\circ\boldsymbol{\mathcal{A}}_2$ is maximally
monotone in the renormed space and, if $\boldsymbol{C}$ is
cocoercive in $\GGG$, then
$\boldsymbol{V}\circ\boldsymbol{\mathcal{B}}_2$ is cocoercive
in the renormed space. Thus, setting
\begin{equation}
\label{e:metricV}
\boldsymbol{V}=
\begin{bmatrix}
W & 0\\
0 & (\sigma^{-1}\ID-\boldsymbol{L}\circ 
W\circ\boldsymbol{L}^*)^{-1}
\end{bmatrix},
\end{equation}
where $W\colon\HH\to\HH$, and applying Proposition~\ref{p:fb13} in
this context yields the following result (see \cite{Icip14}).

\begin{proposition}
\label{p:fb13b}
Suppose that, in Problem~\ref{prob:88}, $A=0$ and 
$(C_k)_{1\leq k\leq q}$ are $\beta$-cocoercive for some 
$\beta\in\RPP$. 
Let $W\colon\HH\to\HH$ be a self-adjoint strictly positive linear
operator and let $\sigma\in\RPP$ be such that 
$\kappa=\|\boldsymbol{L}\circ W\circ\boldsymbol{L}^*\|
<\min\{1/\sigma,2\beta\}$. Let
$\varepsilon\in\left]0,\min\{1/2,\beta/\kappa\}\right[$, let 
$x_0\in\HH$, and let $\boldsymbol{v}_0\in\GGG$.
For every $n\in\NN$, let 
\begin{equation}
\lambda_n\in\big[\varepsilon,(1-\varepsilon)
\big(2+\varepsilon-{\kappa}/{2\beta}\big)\big].
\end{equation}
Iterate
\begin{equation}
\label{e:6}
\begin{array}{l}
\text{for}\;n=0,1,\ldots\\
\left\lfloor
\begin{array}{l}
\text{for}\;k=1,\ldots,q\\
\left\lfloor
\begin{array}{l}
s_{n,k}=C_k(L_k x_n)\\
\end{array}
\right.\\[2mm]
z_n=x_n-W\big(\sum_{k=1}^q L_k^* (s_{n,k}+v_{n,k})\big)\\
\begin{array}{l}
\text{for}\;k=1,\ldots,q\\
\left\lfloor
\begin{array}{l}
w_{n,k}=v_{n,k}+\sigma L_k z_n\\
y_{n,k}=w_{n,k}-\sigma J_{\sigma^{-1} B_k}
\left({w_{n,k}}/{\sigma}\right)\\
v_{n+1,k}=v_{n,k}+\lambda_n (y_{n,k}-v_{n,k})\\
\end{array}
\right.\\[2mm]
\end{array}\\
u_n=x_n-W\big(\sum_{k=1}^q L_k^*(s_{n,k}+y_{n,k})\big)\\
x_{n+1}=x_n+\lambda_n (u_n-x_n).
\end{array}
\right.\\[2mm]
\end{array}
\end{equation}
Then $(x_n)_{n\in\NN}$ converges to a solution to
Problem~\ref{prob:88}.
\end{proposition}

Other choices of the metric operator $\boldsymbol{V}$ are
possible, which lead to different primal-dual algorithms 
\cite{Optim14,Cond13,Komo15,Bang13}. An
advantage of \eqref{e:rio1010b} and \eqref{e:6}
over \eqref{e:3opsbis} is that the first two do not
require the inversion of linear operators.

\subsection{Block-iterative algorithms}
As will be seen in Problems~\ref{prob:2} and \ref{prob:3}, 
systems of inclusions arise 
in multivariate optimization problems (they will 
also be present in Nash equilibria; see, e.g., \eqref{e:as} 
and \eqref{e:vi1}).
We now focus on general systems of inclusions involving
maximally monotone operators as well as linear operators coupling
the variables. 

\begin{problem}
\label{prob:82}
For every $i\in I=\{1,\ldots,m\}$ and 
$k\in K=\{1,\ldots,q\}$, let $A_i\colon\HH_i\to 2^{\HH_i}$ and 
$B_k\colon\GG_k\to 2^{\GG_k}$ be maximally monotone, 
and let $L_{k,i}\colon\HH_i\to\GG_k$ be linear.
The task is to
\begin{multline}
\label{e:2015-02-26p}
\text{find}\;\;\overline{x}_1\in\HH_1,\ldots,
\overline{x}_m\in\HH_m\;\text{such that}\;
(\forall i\in I)\\
0\in A_i\overline{x}_i+\Sum_{k\in K}L_{k,i}^*
\bigg(B_k\bigg(\Sum_{j\in I}L_{k,j}\overline{x}_j\bigg)\bigg),
\end{multline}
under the assumption that the \emph{Kuhn-Tucker set}
\begin{multline}
\label{e:kt3}
\boldsymbol{Z}=\bigg\{(\overline{\boldsymbol{x}},
\overline{\boldsymbol{v}})
\in\HHH\times\GGG\:\bigg |\: (\forall i\in I)\;
-\sum_{k\in K}L_{k,i}^*\overline{v}_k\in
A_i\overline{x}_i\\
\text{and}\:\;(\forall k\in K)\;\;
\sum_{i\in I}L_{k,i}\overline{x}_i\in B_k^{-1}\overline{v}_k
\bigg\}
\end{multline}
is nonempty.
\end{problem}

We can regard $m$ as the number of coordinates of the solution
vector 
$\overline{\boldsymbol{x}}=(\overline{x}_i)_{1\leq i\leq m}$. 
In large-scale applications, $m$ can be sizable and so can
the number of terms $q$, which is often associated with the number
of observations. We have already discussed in
Sections~\ref{sec:upda} and \ref{sec:87} techniques in which not
all the indices $i$ or $k$ need to be activated at a given
iteration. Below, we describe a block-iterative method proposed in
\cite{MaPr18} which allows for partial activation of both the
families $(A_i)_{1\leq i\leq m}$ and $(B_k)_{1\leq k\leq q}$,
together with individual, iteration-dependent proximal parameters
for each operator. The method displays an unprecedented level of
flexibility and it does not require the inversion of linear
operators or knowledge of their norms.

The principle of the algorithm is as follows. Denote by
$I_n\subset I$ and $K_n\subset K$ the
blocks of indices of operators to be updated at iteration $n$. We
impose the mild condition that there exist $M\in\NN$ such that
each operator index $i$ and $k$ is used at least once within any
$M$ consecutive iterations, i.e., for every $n\in\NN$,
\begin{equation}
\label{e:n24G}
\bigcup_{j=n}^{n+M-1}I_j=\{1,\ldots,m\}
\;\:\text{and}\:
\bigcup_{j=n}^{n+M-1}K_j=\{1,\ldots,q\}.
\end{equation}
For each $i\in I_n$ and $k\in K_n$, we 
select points $(a_{i,n},a^*_{i,n})\in\gra A_i$
and $(b_{k,n},b^*_{k,n})\in\gra B_k$ and use them to construct a
closed half-space $\boldsymbol{H}_n\subset\HHH\times\GGG$ which
contains $\boldsymbol{Z}$. The primal variable
$\boldsymbol{x}_n$ and the dual variable $\boldsymbol{v}_n$ are 
updated as $(\boldsymbol{x}_{n+1},\boldsymbol{v}_{n+1})
=\proj_{\boldsymbol{H}_n}(\boldsymbol{x}_n,\boldsymbol{v}_n)$.
The resulting algorithm can also be implemented with relaxations
and in an asynchronous fashion \cite{MaPr18}. For simplicity, we
present the unrelaxed synchronous version.

\begin{proposition}[\cite{MaPr18}]
\label{p:1}
Consider the setting of Problem~\ref{prob:82}.
Take sequences $(I_n)_{n\in\NN}$ in $I$ and 
$(K_n)_{n\in\NN}$ in $K$ satisfying \eqref{e:n24G}, 
with $I_0=I$ and $K_0=K$.
Let $\varepsilon\in\zeroun$ and, for every $i\in I$ 
and every $k\in K$, let $(\gamma_{i,n})_{n\in\NN}$
and $(\mu_{k,n})_{n\in\NN}$ be sequences in
$[\varepsilon,1/\varepsilon]$. 
Let $\boldsymbol{x}_{0}\in\HHH$, let
$\boldsymbol{v}_{0}\in\GGG$, and iterate
\begin{equation}
\label{e:n03a}
\begin{array}{l}
\text{for}\;n=0,1,\ldots\\
\left\lfloor
\begin{array}{l}
\hskip -2mm
\begin{array}{l}
\text{for every}\;i\in I_n\\
\left\lfloor
\begin{array}{l}
l^*_{i,n}=\sum_{k\in K}L_{k,i}^*v_{k,n}\\
a_{i,n}=J_{\gamma_{i,n} A_i}\big(x_{i,n}-\gamma_{i,n}
l^*_{i,n}\big)\\
a_{i,n}^*=\gamma_{i,n}^{-1}(x_{i,n}-a_{i,n})-l^*_{i,n}\\
\end{array}
\right.\\[1mm]
\text{for every}\;i\in I\smallsetminus I_n\\
\left\lfloor
\begin{array}{l}
(a_{i,n},a_{i,n}^*)=(a_{i,n-1},a_{i,n-1}^*)\\
\end{array}
\right.\\[1mm]
\text{for every}\;k\in K_n\\
\left\lfloor
\begin{array}{l}
l_{k,n}=\sum_{i\in I}L_{k,i}x_{i,n}\\
b_{k,n}=J_{\mu_{k,n}B_k}
\big(l_{k,n}+\mu_{k,n}v_{k,n}\big)\\
b^*_{k,n}=v_{k,n}+\mu_{k,n}^{-1}
(l_{k,n}-b_{k,n})\\
\end{array}
\right.\\[1mm]
\text{for every}\;k\in K\smallsetminus K_n\\
\left\lfloor
\begin{array}{l}
(b_{k,n},b^*_{k,n})=(b_{k,n-1},b^*_{k,n-1})\\
\end{array}
\right.\\[1mm]
\text{for every}\;i\in I\\
\left\lfloor
\begin{array}{l}
t^*_{i,n}=a^*_{i,n}+\sum_{k\in K}L_{k,i}^*b^*_{k,n}\\
\end{array}
\right.\\[1mm]
\text{for every}\;k\in K\\
\left\lfloor
\begin{array}{l}
t_{k,n}=b_{k,n}-\sum_{i\in I}L_{k,i}a_{i,n}\\
\end{array}
\right.\\[1mm]
\tau_n=\sum_{i\in I}\|t_{i,n}^*\|^2+\sum_{k\in K}\|t_{k,n}\|^2\\
\text{if}\;\tau_n>0\\
\left\lfloor
\begin{array}{l}
\theta_n=\dfrac{1}{\tau_n}\,\text{\rm max} 
\big\{0,\sum_{i\in I}\big(\scal{x_{i,n}}{t^*_{i,n}}-
\scal{a_{i,n}}{a^*_{i,n}}\big)\\
\hskip 15mm +\sum_{k\in K}
\big(\sscal{t_{k,n}}{v_{k,n}}-\sscal{b_{k,n}}{b^*_{k,n}}\big)
\big\}\\
\end{array}
\right.\\
\text{else~} \theta_n=0\\
\text{for every}\;i\in I\\
\left\lfloor
\begin{array}{l}
x_{i,n+1}=x_{i,n}-\theta_n t^*_{i,n}\\
\end{array}
\right.\\
\text{for every}\;k\in K\\
\left\lfloor
\begin{array}{l}
v_{k,n+1}=v_{k,n}-\theta_n t_{k,n}.
\end{array}
\right.\\
\end{array}
\end{array}
\right.\\[4mm]
\end{array}
\end{equation}
Then $(\boldsymbol{x}_n)_{n\in\NN}$ converges to a solution to
Problem~\ref{prob:82}.
\end{proposition}

Recent developments on splitting algorithms for
Problem~\ref{prob:82} as well as variants and extensions thereof
can be found in \cite{Jmaa20,Moor21,Gise19,John20,John21}.

\section{Fixed point modeling of minimization problems}
\label{sec:5}

We present key applications of fixed point models in convex
optimization.

\subsection{Convex feasibility problems}
\label{sec:cfp}
The most basic convex optimization problem is the convex
feasibility problem, which asks for compliance with a finite
number of convex constraints the object of interest is known to
satisfy. This approach was formalized by Youla 
\cite{Youl78,Youl82} in signal recovery and it has enjoyed a broad
success \cite{Proc93,Aiep96,Herm09,Star87,Thao21,Trus84}.

\begin{problem}
\label{prob:1}
Let $(C_i)_{1\leq i\leq m}$ be nonempty closed convex subsets
of $\HH$. The task is to 
\begin{equation}
\label{e:cfp1}
\text{find}\;\;x\in\bigcap_{i=1}^mC_i.
\end{equation}
\end{problem}

Suppose that Problem~\ref{prob:1} has a solution and that
each set $C_i$ is modeled as the fixed point set of 
an $\alpha_i$-averaged operator $T_i\colon\HH\to\HH$
some some $\alpha_i\in\zeroun$. Then, applying Theorem~\ref{t:1}
with $T=T_1\circ\cdots\circ T_m$ (which is averaged by 
Proposition~\ref{p:roma}) and $\lambda_n=1$ for every $n\in\NN$, we
obtain that the sequence $(x_n)_{n\in\NN}$ constructed via the
iteration 
\begin{equation}
\label{e:g2}
(\forall n\in\NN)\quad x_{n+1}=(T_1\circ\cdots\circ T_m)x_n
\end{equation}
converges to a fixed point $x$ of $T_1\circ\cdots\circ T_m$.
However, in view of Proposition~\ref{p:f3}, $x$ is a solution to
\eqref{e:cfp1}. In particular, if each $T_i$ is the projection
operator onto $C_i$ (which was seen to be $1/2$-averaged), we
obtain the classical POCS (Projection Onto Convex Sets) algorithm
\cite{Breg65,Erem65}
\begin{equation}
\label{e:g3}
(\forall n\in\NN)\quad x_{n+1}=
\big(\proj_{C_1}\circ\cdots\circ\proj_{C_m}\big)x_n
\end{equation}
popularized in \cite{Youl82} and which goes back to \cite{Kacz37}
in the case of affine hyperplanes.
In this algorithm, the projection operators are used
sequentially. Another basic projection method for solving
\eqref{e:cfp1} is the barycentric projection algorithm
\begin{equation}
\label{e:g4}
(\forall n\in\NN)\quad x_{n+1}=\dfrac{1}{m}\sum_{i=1}^m
\proj_{C_i}x_n,
\end{equation}
which uses the projections simultaneously and goes back to 
\cite{Cimm38} in the case of affine hyperplanes.
Its convergence is proved by applying Theorem~\ref{t:1} to 
$T=m^{-1}\sum_{i=1}^m\proj_{C_i}$ which is 
$1/2$-averaged by Example~\ref{ex:2}. More general fixed point
methods are discussed in \cite{Baus96,Nume06,Imag97,Lopu97}.

\subsection{Split feasibility problems}

The so-called \emph{split feasibility problem} is just a convex
feasibility problem involving a linear operator
\cite{Byrn02,Cens94,Cens05}.

\begin{problem}
\label{prob:9}
Let $C\subset\HH$ and $D\subset\GG$ be closed convex sets and let
$0\neq L\colon\HH\to\GG$ be linear. The task is to
\begin{equation}
\label{e:split1}
\text{find}\;\;x\in C\;\;\text{such that}\;\;Lx\in D,
\end{equation}
under the assumption that a solution exists. 
\end{problem}

In principle, we
can reduce this problem to a 2-set version of \eqref{e:g3} with
$C_1=C$ and $C_2=L^{-1}(D)$. However the projection onto $C_2$ is
usually not tractable, which makes projection algorithms
such as \eqref{e:g3} or \eqref{e:g4} not implementable. To work
around this difficulty, let us define $T_1=\proj_{C}$ and
$T_2=\Id-\gamma G_2$, where $G_2=L^*\circ(\Id-\proj_D)\circ L$ 
and $\gamma\in\RPP$. Then $(\forall x\in\HH)$ $Lx\in D$
$\Leftrightarrow$ $G_2x=0$.\footnote{Set $T=\Id-\proj_D$ and fix
$\overline{x}\in\HH$ such that $L\overline{x}\in D$. Then
$T(L\overline{x})=0$ and thus $G_2\overline{x}=0$. Conversely, take
$x\in\HH$ such that $G_2x=0$. Since $T$ is firmly nonexpansive by
Example~\ref{ex:1}, applying \eqref{e:101} with $\beta=1$ yields
$0=\scal{0}{x-\overline{x}}=
\scal{G_2x-G_2\overline{x}}{x-\overline{x}}=
\scal{L^*(T(Lx)-T(L\overline{x}))}{x-\overline{x}}=
\scal{T(Lx)-T(L\overline{x})}{Lx-L\overline{x}}\geq
\|T(Lx)-T(L\overline{x})\|^2=\|T(Lx)\|^2$. So $T(Lx)=0$ and
therefore $Lx=\proj_D(Lx)\in D$.}
Hence, 
\begin{equation}
\label{e:g7}
\Fix T_1=C\quad\text{and}\quad\Fix T_2=\menge{x\in\HH}{Lx\in D}.
\end{equation}
Furthermore, $T_1$ is $\alpha_1$-averaged with $\alpha_1=1/2$. 
In addition, $\Id-\proj_D$ is firmly nonexpansive by 
\eqref{e:110} and therefore $1$-cocoercive. It follows from
Proposition~\ref{p:coco} that $G_2$ is cocoercive with constant 
$1/\|L\|^2$. Now let $\gamma\in\left]0,2/\|L\|^2\right[$ and set
$\alpha_2=\gamma\|L\|^2/2$. Then 
Proposition~\ref{p:2004-5} asserts that 
$\Id-\gamma G_2$ is $\alpha_2$-averaged. Altogether, we deduce from
Example~\ref{ex:roma} that $T_1\circ T_2$ is $\alpha$-averaged. 
Now let 
$(\lambda_n)_{n\in\NN}$ be an $\alpha$-relaxation sequence.
According to Theorem~\ref{t:1} and Proposition~\ref{p:f3}, 
the sequence produced by the iterations
\begin{align}
\label{e:g9}
&\hskip -1mm(\forall n\in\NN)\quad 
x_{n+1}=x_n+\lambda_n\nonumber\\
&\qquad\quad
\cdot\Big(\proj_C\big(x_n-\gamma L^*(Lx_n-\proj_D(Lx_n))\big)
-x_n\Big)\nonumber\\
&\hskip 23.5mm=x_n+\lambda_n\big(T_1(T_2x_n)-x_n\big)
\end{align}
converges to a point in $\Fix T_1\cap\Fix T_2$, i.e., in view of 
\eqref{e:g7}, to a solution to Problem~\ref{prob:9}. In particular,
if we take $\lambda_n=1$, the update rule in \eqref{e:g9}
becomes
\begin{equation}
\label{e:g6}
x_{n+1}=\proj_C\Big(x_n-\gamma L^*\big(Lx_n-
\proj_D(Lx_n)\big)\Big).
\end{equation}

\subsection{Convex minimization}
\label{sec:convmin}

We deduce from Fermat's rule (Theorem~\ref{t:4}) and 
Proposition~\ref{p:11} the fact that a differentiable 
convex function $f\colon\HH\to\RR$ admits $x\in\HH$ as
a minimizer if and only if $\nabla f(x)=0$. Now let
$\gamma\in\RPP$. Then this property is equivalent to
$x=x-\gamma\nabla f(x)$, which shows that 
\begin{equation}
\label{e:fa}
\Argmin f=\Fix T,\quad\text{where}\quad T=\Id-\gamma\nabla f.
\end{equation}
If we add the assumption that $\nabla f$ is $\delta$-Lipschitzian,
then it is $1/\delta$-cocoercive by Proposition~\ref{p:BH}. Hence,
if $0<\gamma<2/\delta$, it follows from Proposition~\ref{p:2004-5},
that $T$ in \eqref{e:fa} is $\alpha$-averaged with 
$\alpha=\gamma\delta/2$. We then derive from Theorem~\ref{t:1} 
the convergence of the steepest-descent method.

\begin{proposition}[steepest-descent]
\label{p:sd}
Let $f\colon\HH\to\RR$ be a differentiable convex function such
that $\Argmin f\neq\emp$ and $\nabla f$ is $\delta$-Lipschitzian
for some $\delta\in\RPP$. Let $\gamma\in\left]0,2/\delta\right[$,
let $(\lambda_n)_{n\in\NN}$ be a 
$\gamma\delta/2$-relaxation sequence,
and let $x_0\in\HH$. Set
\begin{equation}
\label{e:sd1}
(\forall n\in\NN)\quad x_{n+1}=x_n-\gamma\lambda_n\nabla f(x_n).
\end{equation}
Then $(x_n)_{n\in\NN}$ converges to a point in $\Argmin f$. 
\end{proposition}

Now, let us remove the smoothness assumption by considering
a general function $f\in\Gamma_0(\HH)$.
Then it is clear from \eqref{e:moreau1} that $(\forall x\in\HH)$
$x=\prox_fx$ $\Leftrightarrow$ $(\forall y\in\HH)$ $f(x)\leq f(y)$.
In other words, we obtain the fixed point characterization
\begin{equation}
\label{e:fb}
\Argmin f=\Fix T,\quad\text{where}\quad T=\prox_f.
\end{equation}
In turn, since $\prox_f$ is firmly nonexpansive (see
Example~\ref{ex:1}), we derive at once from Theorem~\ref{t:1} the
convergence of the proximal point algorithm.

\begin{proposition}[proximal point algorithm]
\label{p:ppa}
Let $f\in\Gamma_0(\HH)$ be such that $\Argmin f\neq\emp$.
Let $\gamma\in\RPP$, let $(\lambda_n)_{n\in\NN}$ be a
$1/2$-relaxation sequence, and let $x_0\in\HH$. Set
\begin{equation}
\label{e:sd2}
(\forall n\in\NN)\quad x_{n+1}=
x_n+\lambda_n\big(\prox_{\gamma f}x_n-x_n\big).
\end{equation}
Then $(x_n)_{n\in\NN}$ converges to a point in $\Argmin f$. 
\end{proposition}

\begin{remark}
We can interpret the barycentric projection algorithm \eqref{e:g4}
as an unrelaxed instance of the proximal point algorithm
\eqref{e:sd2} with $\gamma=1$ by applying Remark~\ref{r:2018} with
$q=m$ and, for every $k\in\{1,\ldots,q\}$, $\omega_k=1/q$,
$\GG_k=\HH$, $L_k=\Id$, and $g_k=\iota_{C_k}$.
\end{remark}

A more versatile minimization model is the following instance of
the formulation discussed in Proposition~\ref{p:17}.

\begin{problem}
\label{prob:5}
Let $f\in\Gamma_0(\HH)$ and $g\in\Gamma_0(\HH)$ be such that 
$(\reli\dom f)\cap(\reli\dom g)\neq\emp$
and $\lim_{\|x\|\to\pinf}f(x)+g(x)=\pinf$. The task is to 
\begin{equation}
\label{e:r28p}
\minimize{x\in\HH}{f(x)+g(x)}.
\end{equation}
\end{problem}

It follows from Proposition~\ref{p:17}\ref{p:17i} that
Problem~\ref{prob:5} has a solution and from 
Proposition~\ref{p:17}\ref{p:17ii} that
it is equivalent to Problem~\ref{prob:8}
with $A=\partial f$ and $B=\partial g$. It then remains to invoke
Proposition~\ref{p:DR} and Example~\ref{ex:mono2} to obtain the 
following algorithm, which employs the proximity operators of $f$
and $g$ separately.

\begin{proposition}[Douglas-Rachford splitting]
\label{p:18}
Let $(\lambda_n)_{n\in\NN}$ be a $1/2$-relaxation sequence,
let $\gamma\in\RPP$, and let $y_0\in\HH$. Iterate
\begin{equation}
\label{e:DRv}
\begin{array}{l}
\text{for}\;n=0,1,\ldots\\
\left\lfloor
\begin{array}{l}
x_n=\prox_{\gamma g}y_n\\
z_n=\prox_{\gamma f}(2x_n-y_n)\\
y_{n+1}=y_n+\lambda_n(z_n-x_n).
\end{array}
\right.\\[2mm]
\end{array}
\end{equation}
Then $(x_n)_{n\in\NN}$ converges to a
solution to Problem~\ref{prob:5}.
\end{proposition}

The Douglas-Rachford algorithm was first employed in signal and
image processing in \cite{Comb07} and it has since been applied to
various problems, e.g., \cite{Chiz18,Lind21,Papa14,Stei10,Yuyu17}.
For a recent application to joint scale/regression estimation in
statistical data analysis involving several product space
reformulations, see \cite{Ejst20}. We now present two applications
to matrix optimization problems. Along the same lines, the
Douglas-Rachford algorithm is also used in tensor decomposition
\cite{Gand11}.

\begin{example}
Let $\HH$ be the space of $N\times N$ real symmetric matrices 
equipped with the Frobenius norm. We denote by $\xi_{i,j}$ the
$ij$th component of $X\in\HH$. Let $O\in\HH$.
The {\em graphical lasso problem} \cite{Fried08,Ravi11} is to
\begin{equation}
\label{e:GLASSO}
\minimize{X\in\HH}{f(X)+\ell(X)+\operatorname{trace}(OX)},
\end{equation}
where 
\begin{equation}
f(X)=\chi\sum_{i=1}^N\sum_{j=1}^N|\xi_{i,j}|,\quad
\text{with}\;\chi\in\RP,
\end{equation}
and 
\begin{equation}
\ell(X)=
\begin{cases}
-\ln\det X,&\text{if $X$ is positive definite};\\
\pinf,&\text{otherwise.}
\end{cases}
\end{equation}
Problem \eqref{e:GLASSO} arises in the estimation of a sparse
precision (i.e., inverse covariance) matrix from an observed matrix
$O$ and it has found applications in graph processing. Since
$\ell\in\Gamma_0(\HH)$ is a symmetric function of the eigenvalues 
of its arguments, by \cite[Corollary~24.65]{Livre1}, its proximity
operator at $X$ is obtained by performing an eigendecomposition 
$[U,(\mu_i)_{1\leq i\leq N}]=\operatorname{eig}(X)$
$\Leftrightarrow$
$X=U\operatorname{Diag}(\mu_1,\ldots,\mu_N) U^\top$.
Here, given $\gamma\in\RPP$, \cite[Example~24.66]{Livre1} yields
\begin{equation}
\prox_{\gamma\ell}X=U \operatorname{Diag}
\big((\prox_{-\gamma\ln}\mu_1,\ldots,\prox_{-\gamma\ln}\mu_N)\big) 
U^\top,
\end{equation}
where
$\prox_{-\gamma\ln}\colon\xi\mapsto(\xi+\sqrt{\xi^2+4\gamma})/2$.
Let $(\lambda_n)_{n\in\NN}$ be a $1/2$-relaxation sequence,
let $\gamma\in\RPP$, and let $Y_0\in\HH$. Upon setting
$g=\ell+\scal{\cdot}{O}$, the Douglas-Rachford algorithm of
\eqref{e:DRv} for solving \eqref{e:GLASSO} becomes
\begin{equation}
\label{e:DRGLASSO}
\begin{array}{l}
\text{for}\;n=0,1,\ldots\\
\left\lfloor
\begin{array}{l}
[U_n,(\mu_{i,n})_{1\leq i\leq N}]=
\operatorname{eig}(Y_n-\gamma O)\\
X_n=U_n \operatorname{Diag}\big((\prox_{-\gamma\ln}
\,\mu_{i,n})_{1\leq i\leq N}\big)U_n^\top\\
Z_n=\operatorname{soft}_{\gamma \chi}(2X_n-Y_n)\\
Y_{n+1}=Y_n+\lambda_n(Z_n-X_n),
\end{array}
\right.\\[2mm]
\end{array}
\end{equation}
where $\operatorname{soft}_{\gamma \chi}$ denotes the
soft-thresholding operator on $[-\gamma\chi,\gamma\chi]$ applied
componentwise. Applications of \eqref{e:DRGLASSO} as well as
variants with other choices of $\ell$ and $g$ are discussed in
\cite{Benf20}. 
\end{example}

\begin{example}[robust PCA]
\label{ex:robustPCA}
Let $M$ and $N$ be integers such that $M\geq N>0$, and let $\HH$
be the space of $N\times M$ real matrices equipped with the
Frobenius norm. The robust Principal Component Analysis (PCA) 
problem \cite{Cand11,Vasw18} is to
\begin{equation}
\label{e:robustPCA}
\minimize{\substack{X\in\HH, Y\in\HH\\ X+Y=O}}
{\|Y\|_{\rm nuc}+\chi\|X\|_1},
\end{equation}
where $\|\cdot\|_1$ is the componentwise $\ell_1$-norm,
$\|\cdot\|_{\rm nuc}$ is the nuclear norm, and $\chi\in\RPP$.
Let $X=U\operatorname{Diag}(\sigma_1,\ldots,\sigma_N)V^\top$
be the singular value decomposition of $X\in\HH$.
Then $\|X\|_{\rm nuc}=\sum_{i=1}^N \sigma_i$ and, by 
\cite[Example~24.69]{Livre1},
\begin{equation}
\prox_{\chi\|\cdot\|_{\rm nuc}}X=U\operatorname{Diag}
\big(\operatorname{soft}_{\chi}\sigma_1,\ldots,
\operatorname{soft}_{\chi}\sigma_N\big)V^\top.
\end{equation}
An implementation of the Douglas-Rachford algorithm in the product 
space $\HH\times\HH$ to solve \eqref{e:robustPCA} is detailed in 
\cite[Example~28.6]{Livre1}.
\end{example}

By combining Propositions~\ref{p:fb13}, \ref{p:11}, and \ref{p:BH},
together with Example~\ref{ex:mono2}, we obtain the convergence of
the forward-backward splitting algorithm for minimization. The
broad potential of this algorithm in data science was evidenced in
\cite{Smms05}. Inertial variants are presented in
\cite{Apid20,Atto19,Beck09,Biou07,Cham15,Siop17}.

\begin{proposition}[forward-backward splitting]
\label{p:fb17}
Suppose that, in Problem~\ref{prob:5}, $g$ is differentiable
everywhere and that its gradient is $\delta$-Lipschitzian for 
some $\delta\in\RPP$. Let
$\varepsilon\in\left]0,\min\{1/2,1/\delta\}\right[$, let
$x_0\in\HH$, and let $(\gamma_n)_{n\in\NN}$ be in
$\left[\varepsilon,2/(\delta(1+\varepsilon))\right]$, and let 
\begin{equation}
\label{e:215}
(\forall n\in\NN)\quad
\lambda_n\in\big[\varepsilon,(1-\varepsilon)
\big(2+\varepsilon-{\delta\gamma_n}/{2}\big)\big].
\end{equation}
Iterate
\begin{equation}
\label{e:FB2}
\begin{array}{l}
\text{for}\;n=0,1,\ldots\\
\left\lfloor
\begin{array}{l}
u_n=x_n-\gamma_n\nabla g(x_n)\\
x_{n+1}=x_n+\lambda_n\big(\prox_{\gamma_n f}u_n-x_n\big).
\end{array}
\right.\\[2mm]
\end{array}
\end{equation}
Then $(x_n)_{n\in\NN}$ converges to a solution to 
Problem~\ref{prob:5}.
\end{proposition}

\begin{example}
Let $M$ and $N$ be integers such that $M\geq N>0$, and let $\HH$
be the space of $N\times M$ real-valued matrices 
equipped with the Frobenius norm. The task is to reconstruct 
a low-rank matrix given its projection $O$ onto a vector space
$V\subset\HH$. Let $L=\proj_V$. The problem is formulated as
\begin{equation}
\label{e:matcomp}
\minimize{X\in\HH}{\frac12\|O-LX\|^2+\chi\|X\|_{\rm nuc}},
\end{equation}
where $\chi\in\RPP$. As seen in Example~\ref{ex:robustPCA}, the
proximity operator of the nuclear norm has a closed form
expression. In addition, $g\colon X\mapsto\|O-LX\|^2/2$ is convex
and its gradient $\nabla g\colon X\mapsto L^*(LX-O)=LX-O$ is
nonexpansive. Problem \eqref{e:matcomp} can thus be solved by
algorithm \eqref{e:FB2} where $f=\chi\|\cdot\|_{\rm nuc}$ and
$\delta=1$. A particular case of \eqref{e:matcomp} is the matrix
completion problem \cite{Cand09,Cand10}, where only some components
of the sought matrix are observed. If $\mathbb{K}$ denotes the set
of indices of the unknown matrix components, we have 
$V=\menge{X\in\HH}{(\forall (i,j)\in\mathbb{K})\;\xi_{i,j}=0}$.
\end{example}

\begin{example}
\label{ex:MMSE}
Let $X$ and $W$ be mutually independent $\RR^N$-valued random
vectors. Assume that $X$ is absolutely continuous and 
square-integrable, and that its probability density function is
log-concave. Further, assume that $W$ is Gaussian with zero-mean 
and covariance $\sigma^2\mathrm{I}_{N}$, where $\sigma\in\RPP$. 
Let $Y=X+W$. For every $y\in\RR^N$, $Qy=\EE(X\mid Y=y)$ is the
minimum mean square error (MMSE) denoiser for $X$ given the
observation $y$. The properties of $Q$ have been investigated in
\cite{Gribo13}. It can be shown that $Q$ is the proximity operator
of the conjugate of $h=(-\sigma^2\log p)^*-\|\cdot\|^2/2\in
\Gamma_{0}(\RR^N)$, where $p$ is the density of $Y$. Let
$g\colon\RR^N\to\RR$ be a differentiable convex function with a
$\delta$-Lipschitzian gradient for some $\delta\in\RPP$, and let
$\gamma\in\left]0,2/\delta\right[$. The iteration 
\begin{equation}
(\forall n\in\NN)\quad 
x_{n+1}=Q\big(x_{n}-\gamma\nabla g(x_n)\big)
\end{equation}
therefore turns out to be a special case of the forward-backward
algorithm \eqref{e:FB2}, where $f=h^*/\gamma$ and 
$(\forall n\in\NN)$ $\lambda_n=1$. This algorithm is studied in
\cite{Xu20} from a different perspective. 
\end{example}

The projection-gradient method goes back to the classical papers
\cite{Gold64,Lev66a}. A version can be obtained by setting
$f=\iota_C$ in Proposition~\ref{p:fb17}, where $C$ is the
constraint set. Below, we describe the simpler formulation
resulting from the application of Theorem~\ref{t:1} to
$T=\proj_C\circ(\Id-\gamma\nabla g)$.

\begin{example}[projection-gradient]
\label{ex:19}
Let $C$ be a nonempty closed convex subset of $\HH$ and let
$g\colon\HH\to\RR$ be a differentiable convex function, with 
a $\delta$-Lipschitzian gradient for some $\delta\in\RPP$. 
The task is to
\begin{equation}
\label{e:28q}
\minimize{x\in C}{g(x)},
\end{equation}
under the assumption that $\lim_{\|x\|\to\pinf}g(x)=\pinf$ or $C$
is bounded. Let $\gamma\in\left]0,2/\delta\right[$ and set
$\alpha=2/(4-\gamma\delta)$. Furthermore, let 
$(\lambda_n)_{n\in\NN}$ be an
$\alpha$-relaxation sequence and let $x_0\in\HH$. Iterate
\begin{equation}
\label{e:pg1}
\begin{array}{l}
\text{for}\;n=0,1,\ldots\\
\left\lfloor
\begin{array}{l}
y_n=x_n-\gamma\nabla g(x_n)\\
x_{n+1}=x_n+\lambda_n\big(\proj_{C}y_n-x_n\big).
\end{array}
\right.\\[2mm]
\end{array}
\end{equation}
Then $(x_n)_{n\in\NN}$ converges to a solution to 
\eqref{e:28q}.
\end{example}

As a special case of Example~\ref{ex:19}, we obtain the convergence
of the \emph{alternating projections} algorithm
\cite{Che59a,Lev66a}.

\begin{example}[alternating projections]
\label{ex:20}
Let $C_1$ and $C_2$ be nonempty closed convex subsets of $\HH$, one
of which is bounded. Given $x_0\in\HH$, iterate
\begin{equation}
\label{e:u8}
(\forall n\in\NN)\quad 
x_{n+1}=\proj_{C_1}\big(\proj_{C_2}x_n\big).
\end{equation}
Then $(x_n)_{n\in\NN}$ converges to a solution to the 
constrained minimization problem
\begin{equation}
\label{e:28f}
\minimize{x\in C_1}{d_{C_2}(x)}.
\end{equation}
This follows from Example~\ref{ex:19} applied to $g=d_{C_2}^2/2$.
Note that $\nabla g=\Id-\proj_{C_2}$ has Lipschitz constant
$\delta=1$ (see Example~\ref{ex:jjm7}) and hence \eqref{e:u8} is
the instance of \eqref{e:pg1} obtained by setting $\gamma=1$ and
$(\forall n\in\NN)$ $\lambda_n=1$ (see
Example~\ref{ex:relax}\ref{ex:relaxi}).
\end{example}

The following version of Problem~\ref{prob:5} involves $m$ smooth
functions. 

\begin{problem}
\label{prob:33}
Let $(\omega_i)_{1\leq i\leq m}$ be real numbers in $\rzeroun$ such
that $\sum_{i=1}^m\omega_i=1$. Let $f_0\in\Gamma_0(\HH)$ and, for
every $i\in\{1,\ldots,m\}$, let $\delta_i\in\RPP$ and let
$f_i\colon\HH\to\RR$ be a differentiable convex function with a
$\delta_i$-Lipschitzian gradient. Suppose that
\begin{equation}
\label{e:coer1}
\lim_{\|x\|\to\pinf}f_0(x)+\sum_{i=1}^m\omega_i f_i(x)=\pinf.
\end{equation}
The task is to 
\begin{equation}
\label{e:prob33}
\minimize{x\in\HH}{f_0(x)+\sum_{i=1}^m\omega_i f_i(x)}. 
\end{equation}
\end{problem}

To solve Problem~\ref{prob:33}, an option is to apply 
Theorem~\ref{t:3} to obtain a forward-backward algorithm with
block-updates. 

\begin{proposition}[\cite{Upda20}]
\label{p:33}
Consider the setting of Problem~\ref{prob:33}. Let 
$(I_n)_{n\in\NN}$ be a sequence
of nonempty subsets of $\{1,\ldots,m\}$ such that \eqref{e:K} holds
for some $M\in\NN\smallsetminus\{0\}$. 
Let $\gamma\in\left]0,2/\max_{1\leq i\leq m}\delta_i\right[$,
let $x_0\in\HH$, let $(t_{i,-1})_{1\leq i\leq m}\in\HH^m$, and 
iterate
\begin{equation}
\label{e:a4}
\begin{array}{l}
\text{for}\;n=0,1,\ldots\\
\left\lfloor
\begin{array}{l}
\text{for every}\;i\in I_n\\
\left\lfloor
\begin{array}{l}
t_{i,n}=x_n-\gamma\nabla f_i(x_n)\\
\end{array}
\right.\\
\text{for every}\;i\in\{1,\ldots,m\}\smallsetminus I_n\\
\left\lfloor
\begin{array}{l}
t_{i,n}=t_{i,n-1}\\
\end{array}
\right.\\[1mm]
x_{n+1}=\prox_{\gamma f_0}\big(\sum_{i=1}^m\omega_it_{i,n}\big).
\end{array}
\right.\\
\end{array}
\end{equation}
Then the following hold:
\begin{enumerate}
\setlength{\itemsep}{0pt}
\item 
\label{p:33i}
Let $x$ be a solution to Problem~\ref{prob:33} and let
$i\in\{1,\ldots,m\}$. Then $\nabla f_i(x_n)\to\nabla f_i(x)$.
\item 
\label{p:33ii}
$(x_n)_{n\in\NN}$ converges to a solution to
Problem~\ref{prob:33}. 
\item
\label{p:33iv}
Suppose that, for some $i\in\{0,\ldots,m\}$, $f_i$ is strongly
convex. Then $(x_n)_{n\in\NN}$ converges linearly to
the unique solution to Problem~\ref{prob:33}.
\end{enumerate}
\end{proposition}

A method related to \eqref{e:a4} is proposed in \cite{Mish20};
see also \cite{Mokh18} for a special case. Here is a data analysis
application. 

\begin{example}
\label{ex:23}
Let $(e_k)_{1\leq k\leq N}$ be an orthonormal basis of $\HH$ and,
for every $k\in\{1,\ldots,N\}$, let 
$\psi_k\in\Gamma_0(\RR)$. For
every $i\in\{1,\ldots,m\}$, let $0\neq a_i\in\HH$, let 
$\mu_i\in\RPP$, and let
$\phi_i\colon\RR\to\RP$ be a differentiable convex function such
that $\phi_i'$ is $\mu_i$-Lipschitzian. The task is to 
\begin{equation}
\label{e:prob9}
\minimize{x\in\HH}{\sum_{k=1}^N\psi_k(\scal{x}{e_k})+\dfrac{1}{m}
\sum_{i=1}^m\phi_i(\scal{x}{a_i})}. 
\end{equation}
As shown in \cite{Upda20}, \eqref{e:prob9} is an instantiation of 
\eqref{e:prob33} and, given
$\gamma\in\left]0,2/(\max_{1\leq i\leq m}\mu_i\|a_i\|^2)\right[$ 
and subsets $(I_n)_{n\in\NN}$ of $\{1,\ldots,m\}$ such that 
\eqref{e:K} holds, it can be solved by \eqref{e:a4}, which becomes
\begin{equation}
\label{e:a41}
\begin{array}{l}
\text{for}\;n=0,1,\ldots\\
\left\lfloor
\begin{array}{l}
\text{for every}\;i\in I_n\\
\left\lfloor
\begin{array}{l}
t_{i,n}=x_n-\gamma\phi_i'(\scal{x_n}{a_i})a_i\\
\end{array}
\right.\\
\text{for every}\;i\in\{1,\ldots,m\}\smallsetminus I_n\\
\left\lfloor
\begin{array}{l}
t_{i,n}=t_{i,n-1}\\
\end{array}
\right.\\[1mm]
y_n=\sum_{i=1}^m\omega_it_{i,n}\\
x_{n+1}=\sum_{k=1}^N\big(
\prox_{\gamma\psi_k}\scal{y_n}{e_k}\big)e_k.
\end{array}
\right.\\
\end{array}
\end{equation}
A popular setting is obtained by choosing $\HH=\RR^N$ and 
$(e_k)_{1\leq k\leq N}$ as the canonical basis, $\alpha\in\RPP$,
and, for every $k\in\{1,\ldots,K\}$, $\psi_k=\alpha|\cdot|$. This
reduces \eqref{e:prob9} to
\begin{equation}
\label{e:prob91}
\minimize{x\in\RR^N}{\alpha\|x\|_1+\sum_{i=1}^m
\phi_i(\scal{x}{a_i})}. 
\end{equation}
Choosing, for every $i\in\{1,\ldots,m\}$,
$\phi_i\colon t\mapsto|t-\eta_i|^2$
where $\eta_i\in\RR$ models an observation, yields the lasso 
formulation, whereas choosing 
$\phi_i\colon t\mapsto\ln(1+\exp(t))-\eta_it$, where
$\eta_i\in\{0,1\}$ models a label, yields the penalized 
logistic regression framework \cite{Hast09}. 
\end{example}

Next, we extend Problem~\ref{prob:5} to a flexible composite
minimization problem. See 
\cite{Botr14,Chie15,Chie14,Chou19,Chou15,Icip14,Siim19,%
Ejst20,Moer15,Papa14,Pham14,Repe19}
for concrete instantiations of this model in data science.

\begin{problem}
\label{prob:888}
Let $\delta\in\RPP$ and let $f\in\Gamma_0(\HH)$.
For every $k\in\{1,\ldots,q\}$, let $g_k\in\Gamma_0(\GG_k)$, let 
$0\neq L_k\colon\HH\to\GG_k$ be linear, and let
$h_k\colon\GG_k\to\RR$ be a differentiable convex function, with 
a $\delta$-Lipschitzian gradient. Suppose that 
$\lim_{\|x\|\to\pinf}f(x)+\sum_{k=1}^q (g_k(L_kx)+
h_k(L_kx))=\pinf$ and that
\begin{equation}
\label{e:z3}
(\exi z\in\reli\dom f)(\forall k\in\{1,\ldots,q\})\quad 
L_kz\in\reli\dom g_k.
\end{equation}
The task is to 
\begin{equation}
\minimize{x\in\HH}{f(x)+\sum_{k=1}^q\big(g_k(L_kx)+h_k(L_kx)\big)}.
\end{equation}
\end{problem}

Thanks to the qualification condition \eqref{e:z3},
Problem~\ref{prob:888} is an instance of Problem~\ref{prob:88}
where $A=\partial f$ and, for every $k\in\{1,\ldots,q\}$,
$B_k=\partial g_k$ and $C_k=\nabla g_k$. Since the operators
$(C_k)_{1\leq k\leq q}$ are $1/\delta$-cocoercive, the iterative
algorithms from Propositions~\ref{p:71b}, \ref{prop:MLFB},
and~\ref{p:fb13b} are applicable. For example,
Proposition~\ref{p:fb13b} with the substitution 
$J_{\sigma^{-1} B_k}=\prox_{\sigma^{-1} g_k}$ 
(see Example~\ref{ex:mono2}) allows us to solve the problem.
In particular, the resulting algorithm was proposed in
\cite{Chen13,Lori11} in the case when $W=\tau\Id$ with
$\tau\in\RPP$. See also
\cite{Cham11,Optim14,Cond13,Cond20,Esse10,Xiao12,Komo15,Bang13} 
for related work. 

\begin{example}
Let $o\in\RR^N$ and let $M\in\RR^{K\times N}$ be such that
$\mathrm{I}_N-M^\top M$ is positive semidefinite. 
Let $\varphi\in\Gamma_0(\RR^N)$ and let
$C$ be a nonempty closed convex subset of $\RR^N$.
The denoising problem of \cite{Sele20} is cast as
\begin{equation}
\label{e:debnonconv}
\minimize{x\in C}{\psi(x)+\frac12\|x-o\|^2},
\end{equation}
where the function
\begin{equation}
\psi\colon x\mapsto \varphi(x)-\inf_{y\in\HH}\Big(
\varphi(y)+\frac12\|M(x-y)\|^2\Big)
\end{equation}
is generally nonconvex. However,
\eqref{e:debnonconv} is a convex problem. 
Further developments can be found in \cite{Abe20}.
Note that \eqref{e:debnonconv} is actually
equivalent to Problem~\ref{prob:888} with $q=2$,
$\HH=\RR^N\times\RR^N$,
$\GG_1=\HH$, $\GG_2=\RR^{N}$,
$f\colon (x,y)\mapsto\varphi(x)$,
$h_1\colon (x,y)\mapsto\iota_{C}(x)$,
$g_1\colon (x,y)\mapsto x^\top(\mathrm{I}_N-M^\top M) 
x/2-\scal{x}{o}+\|My\|^2/2$, $g_2=\varphi^*$,
$L_1=\Id$, $L_2\colon (x,y)\mapsto M^\top M (x-y)$, and
$h_2=0$. 
\end{example}

\begin{remark}[ADMM]
\label{r:ADMM}
Let us revisit the composite minimization problem of
Proposition~\ref{p:17} and Example~\ref{ex:m+s}.
Let $f\in\Gamma_0(\HH)$, let $g\in\Gamma_0(\GG)$, and let
$L\colon\HH\to\GG$ be linear. Suppose that
$\lim_{\|x\|\to\pinf}f(x)+g(Lx)=\pinf$ and 
$\reli(L(\dom f))\cap\reli(\dom g)\neq\emp$. Then the problem
\begin{equation}
\label{e:r24p}
\minimize{x\in\HH}{f(x)+g(L x)}
\end{equation}
is a special case of Problem~\ref{prob:888} and it can therefore
be solved by any of the methods discussed above. Now let
$\gamma\in\RPP$ and let us make
the following additional assumptions:
\begin{enumerate}
\item
$L^*\circ L$ is invertible.
\item
The operator 
\[
\prox_{\gamma f}^L\colon\GG\to\HH\colon 
y\mapsto\underset{x\in\HH}{\operatorname{argmin}}
\bigg(f(x)+\dfrac{\|Lx-y\|^2}{2}\bigg)
\]
is easy to implement. 
\end{enumerate}
Then, given $y_0\in\GG$ and $z_0\in\GG$, the 
alternating-direction method of multipliers (ADMM) 
constructs a sequence $(x_n)_{n\in\NN}$ that converges to a
solution to \eqref{e:r24p} via the iterations
\cite{Boyd10,Ecks92,Gaba83,Glow89}
\begin{equation}
\begin{array}{l}
\text{for}\;n=0,1,\ldots\\
\left\lfloor
\begin{array}{l}
x_n=\prox_{\gamma f}^L(y_n-z_n)\\
d_n=Lx_n\\
y_{n+1}=\prox_{\gamma g}(d_n+z_n)\\
z_{n+1}=z_n+d_n-y_{n+1}.
\end{array}
\right.\\[2mm]
\end{array}
\end{equation}
This iteration process can be viewed as an application of the
Douglas-Rachford algorithm \eqref{e:DRv} to the Fenchel dual of
\eqref{e:r24p} \cite{Gaba83,Ecks92}. Variants of this algorithm are
discussed in \cite{Bane21,Banf11,Ecks94}, and applications to image
recovery in \cite{Afon10,Afon11,Figu10,Giov05,Gold09,Setz10}.
\end{remark}

\subsection{Inconsistent feasibility problems}
\label{e:inc1}

We consider a more structured variant of Problem~\ref{prob:1} which
can also be considered as an extension of Problem~\ref{prob:9}.

\begin{problem}
\label{prob:10}
Let $C$ be a nonempty closed convex subset of $\HH$ and, for every
$i\in\{1,\ldots,m\}$, let $L_i\colon\HH\to\GG_i$ be a nonzero
linear operator and let $D_i$ be a nonempty closed convex subset
of $\GG_i$. The task is to 
\begin{multline}
\label{e:split3}
\text{find}\;\;x\in C\;\;\text{such that}\;
(\forall i\in\{1,\ldots,m\})\;L_ix\in D_i.
\end{multline}
\end{problem}

To address the possibility that this problem has no solution due to
modeling errors \cite{Cens18,Sign94,Youl86}, we fix 
weights $(\omega_i)_{1\leq i\leq m}$ in $\rzeroun$ such that 
$\sum_{i=1}^m\omega_i=1$ and consider the surrogate problem
\begin{equation}
\label{e:hig1}
\minimize{x\in C}{\frac12\sum_{i=1}^m\omega_id_{D_i}^2(L_ix)},
\end{equation}
where $C$ acts as a hard constraint. This is a valid relaxation of
\eqref{e:split3} in the sense that, if \eqref{e:split3} does have
solutions, then those are the only solutions to \eqref{e:hig1}. 
Now set $f_0=\iota_C$. In addition, for every $i\in\{1,\ldots,m\}$,
set $f_i\colon x\mapsto (1/2)d_{D_i}^2(L_ix)$
and notice that $f_i$ is differentiable and that its gradient 
$\nabla f_i=L_i^*\circ(\Id-\proj_{D_i})\circ L_i$ has Lipschitz
constant $\delta_i=\|L_i\|^2$. Furthermore, \eqref{e:coer1} holds
as long as $C$ is bounded or, for some $i\in\{1,\ldots,m\}$, $D_i$
is bounded and $L_i$ is invertible. We have thus cast
\eqref{e:hig1} as an instance of Problem~\ref{prob:33}
\cite{Upda20}. In view of \eqref{e:a4}, a solution is found as the
limit of the sequence $(x_n)_{n\in\NN}$ produced by the
block-update algorithm 
\begin{equation}
\label{e:a34}
\begin{array}{l}
\text{for}\;n=0,1,\ldots\\
\left\lfloor
\begin{array}{l}
\text{for every}\;i\in I_n\\
\left\lfloor
\begin{array}{l}
t_{i,n}=x_n+\gamma L_i^*\big(\proj_{D_i}(L_ix_n)-L_ix_n\big)\\
\end{array}
\right.\\
\text{for every}\;i\in\{1,\ldots,m\}\smallsetminus I_n\\
\left\lfloor
\begin{array}{l}
t_{i,n}=t_{i,n-1}\\
\end{array}
\right.\\[1mm]
x_{n+1}=\proj_{C}\big(\sum_{i=1}^m\omega_it_{i,n}\big),
\end{array}
\right.\\
\end{array}
\end{equation}
where $\gamma$ and $(I_n)_{n\in\NN}$ are as in
Proposition~\ref{p:33}.

\subsection{Stochastic forward-backward method}
Consider the minimization of $f+g$,
where $f\in\Gamma_0(\HH)$ and 
$g\colon\HH\to\RR$ is a differentiable convex function.
In certain applications, it may happen that only stochastic
approximations to $f$ or $g$ are available.
A generic stochastic form of the forward-backward
algorithm for such instances is \cite{Comb16}
\begin{equation}
\label{e:FBstoch}
(\forall n\in\NN)\quad
x_{n+1}=
x_n+\lambda_n\big(\prox_{\gamma_n f_n}
(x_n-\gamma_nu_n)+a_n-x_n\big),
\end{equation}
where $\gamma_n\in\RPP$, $\lambda_n\in\rzeroun$, 
$f_n\in\Gamma_0(\HH)$ is an approximation to $f$,
$u_n$ is a random variable approximating 
$\nabla g(x_n)$, and $a_n$ is a random variable modeling a
possible additive error. When $f=f_n=0$,
$\lambda_n=1$, and $a_n=0$, we
recover the standard stochastic gradient method for minimizing
$g$, which was pioneered in \cite{Ermo69,Ermo66}.

\begin{example}
As in Problem~\ref{prob:33},
let $f\in\Gamma_0(\HH)$ and let $g=m^{-1}\sum_{i=1}^m g_i$, where
each $g_i\colon\HH\to\RR$ is a
differentiable convex function. The following specialization of 
\eqref{e:FBstoch} is obtained by setting, for every $n\in\NN$,
$f_n=f$ and $u_n=\nabla g_{\mathrm{i}(n)}(x_n)$,
where $\mathrm{i}(n)$ is a $\{1,\ldots,m\}$-valued random
variable. This leads to the incremental proximal stochastic
gradient algorithm described by the update equation 
\begin{equation}
x_{n+1}=x_n+\lambda_n\Big(\prox_{\gamma_n f}
\big(x_n-\gamma_n \nabla g_{\mathrm{i}(n)}(x_n)\big)
-x_n\Big).
\end{equation}
For related algorithms, see
\cite{Bert11,Def14a,Defa14,John13,Schm17}.
\end{example} 

Various convergence results have been established for
algorithm~\eqref{e:FBstoch}. If $\nabla g$ is Lipschitzian,
\eqref{e:FBstoch} is closely related to the fixed point
iteration in Theorem~\ref{t:1stoch}. The almost sure
convergence of $(x_n)_{n\in\NN}$ to a minimizer of $f+g$
can be guaranteed in several scenarios
\cite{Atch17,Comb16,Rosa20}. Fixed point strategies allow us to
derive convergence results such as the following.

\begin{theorem}[\cite{Comb16}]
\label{t:23}
Let $f\in\Gamma_0(\HH)$, let $\delta\in\RPP$, and let 
$g\colon\HH\to\RR$ be a differentiable convex function 
such that $\nabla g$ is $\delta$-Lipschitzian and
$S=\Argmin(f+g)\neq\emp$.
Let $\gamma\in\left]0,2/\delta\right[$ and
let $(\lambda_n)_{n\in\NN}$ be a sequence in 
$\left]0,1\right]$ such that 
$\sum_{n\in\NN}\lambda_n=\pinf$.
Let $x_0$, $(u_n)_{n\in\NN}$, and $(a_n)_{n\in\NN}$ 
be $\HH$-valued random variables with finite second-order moments. 
Let $(x_n)_{n\in\NN}$ be a sequence produced by 
\eqref{e:FBstoch} with $\gamma_n=\gamma$ and $f_n=f$.
For every $n\in\NN$, let $\XX_n$ be
the $\sigma$-algebra generated by $(x_0,\ldots,x_n)$ and set
$\zeta_n=\EC{\|u_n-\EC{u_n}{\XX_n}\|^2}{\XX_n}$.
Assume that the following are satisfied \as:
\begin{enumerate}
\item 
\label{a:t2i}
$\sum_{n\in\NN}\lambda_n\sqrt{\EC{\|a_n\|^2}{\XX_n}}<\pinf$.
\item 
\label{a:t2ii}
$\sum_{n\in\NN}\sqrt{\lambda_n}
\|\EC{u_n}{\XX_n}-\nabla g(x_n)\|<\pinf$.
\item 
\label{a:t2iii} 
$\sup_{n\in\NN} \zeta_n<\pinf$ and $\sum_{n\in\NN}
\sqrt{\lambda_n\zeta_n}<\pinf$.
\end{enumerate}
Then $(x_n)_{n\in\NN}$ converges $\as$ to 
an $S$-valued random variable.
\end{theorem}

Extensions of these stochastic optimization approaches can be
designed by introducing an inertial parameter \cite{Rosa16b} or 
by bringing into play primal-dual formulations \cite{Comb16}. 

\subsection{Random block-coordinate optimization algorithms}
We design block-coordinate versions of optimization algorithms 
presented in Section~\ref{sec:convmin}, in which blocks 
of variables are updated randomly.

\begin{problem}
\label{prob:2}
For every $i\in\{1,\ldots,m\}$ and $k\in\{1,\ldots,q\}$, 
let $f_i\in\Gamma_0(\HH_i)$, let $g_k\in\Gamma_0(\GG_k)$, and let
$0\neq L_{k,i}\colon\HH_i\to\GG_k$ be linear. Suppose that 
\begin{multline}
(\exi\boldsymbol{z}\in\HHH)(\exi\boldsymbol{w}\in\GGG)
(\forall i\in\{1,\ldots,m\})
(\forall k\in\{1,\ldots,q\})\\
-\sum_{j=1}^qL_{j,i}^*w_j\in\partial f_i(z_i)
\;\:\text{and}\:\;
\sum_{j=1}^mL_{k,j}z_j\in \partial g_k^*(w_{k}).
\end{multline}
The task is to 
\begin{equation}
\label{p:probstocopt}
\minimize{\boldsymbol{x}\in\HHH}
{\sum_{i=1}^m f_i(x_i)+\sum_{k=1}^q
g_k\bigg(\sum_{i=1}^m L_{k,i}x_i\bigg)}.
\end{equation}
\end{problem}

Let $\gamma\in\RPP$, let $(\lambda_n)_{n\in\NN}$ be a sequence
in $\left]0,2\right[$, and set
\begin{multline}
\boldsymbol{V}=\bigg\{(x_1,\ldots,x_m,y_1,\ldots,y_q)\in
\HHH\times\GGG\\
\bigg|~ (\forall k\in\{1,\ldots,q\})\;y_k=
\sum_{i=1}^m L_{k,i}x_i\bigg\}
\end{multline}
Let us decompose $\proj_{\boldsymbol{V}}$ as
$\proj_{\boldsymbol{V}}\colon\boldsymbol{x}\mapsto
({Q}_j\boldsymbol{x})_{1\leq j\leq m+q}$.
A random block-coordinate form of the Douglas-Rachford algorithm
for solving Problem~\ref{prob:2} is \cite{Comb15}
\begin{equation}
\label{e:DRbcr}
\begin{array}{l}
\text{for}\;n=0,1,\ldots\\
\left\lfloor
\begin{array}{l}
\text{for}\;i=1,\ldots,m\\
\left\lfloor
\begin{array}{l}
z_{i,n+1}=z_{i,n}+\varepsilon_{i,n}
\big({Q}_i(\boldsymbol{x}_n,\boldsymbol{y}_n)
-z_{i,n}\big)\\[1mm]
x_{i,n+1}=x_{i,n}\\
\qquad\qquad+\varepsilon_{i,n}\lambda_n
\big(\prox_{\gamma f_i}(2z_{i,n+1}-x_{i,n})-z_{i,n+1}\big)
\end{array}
\right.\\
\text{for}\;k=1,\ldots,q\\
\left\lfloor
\begin{array}{l}
w_{k,n+1}=w_{k,n}+\varepsilon_{m+k,n}
\big({Q}_{m+k}(\boldsymbol{x}_n,\boldsymbol{y}_n)
-w_{k,n}\big)\\[1mm]
y_{k,n+1}=y_{k,n}\\
\;\;+\varepsilon_{m+k,n}\lambda_n
\big(\prox_{\gamma g_k}
(2w_{k,n+1}-y_{k,n})-w_{k,n+1}\big),
\end{array}
\right.
\end{array}
\right.\\
\end{array}
\end{equation}
where $\boldsymbol{x}_n=(x_{i,n})_{1\leq i\leq m}$ and 
$\boldsymbol{y}_n=(y_{k,n})_{1\leq k\leq q}$. Moreover,
$(\varepsilon_{j,n})_{1\leq j\leq m+q,n\in\NN}$ are binary random
variables signaling the activated components. 

\begin{proposition}[\cite{Comb15}]
\label{p:l4}
Let $\boldsymbol{S}$ be the set of solutions to 
Problem~\ref{prob:2} and set 
$D=\{0,1\}^{m+q}\smallsetminus\{\boldsymbol{0}\}$.
Let $\gamma\in\RPP$, let $\epsilon\in\zeroun$,
let $(\lambda_n)_{n\in\NN}$ 
be in $\left[\epsilon,2-\epsilon\right]$,
let $\boldsymbol{x}_0$ and $\boldsymbol{z}_0$ be $\HHH$-valued
random variables, let $\boldsymbol{y}_0$ and $\boldsymbol{w}_0$ be
$\GGG$-valued random variables, and let
$(\boldsymbol{\varepsilon}_n)_{n\in\NN}$ be 
identically distributed $D$-valued random variables. 
In addition, suppose that the following hold:
\begin{enumerate}
\item
\label{c:2014-04-09iiv-}
For every $n\in\NN$, $\boldsymbol{\varepsilon}_n$ and
$(\boldsymbol{x}_0,\ldots,\boldsymbol{x}_n,\boldsymbol{y}_0,
\ldots,\boldsymbol{y}_n)$ are mutually independent.
\item
\label{c:2014-04-09iiv}
$(\forall j\in\{1,\ldots,m+q\})$ 
$\PP[\varepsilon_{j,0}=1]>0$.
\end{enumerate}
Then the sequence $(\boldsymbol{z}_n)_{n\in\NN}$ generated by
\eqref{e:DRbcr} converges $\as$ to an 
$\boldsymbol{S}$-valued random variable.
\end{proposition}

Applications based on Proposition~\ref{p:l4} appear in the areas
of machine learning \cite{Nume19} and binary logistic regression
\cite{Bric19}.

If the functions $(g_k)_{1\leq k\leq q}$ are differentiable in
Problem~\ref{prob:2}, a block-coordinate version of the
forward-backward algorithm can also be employed, namely,
\begin{equation}
\label{e:bcFB}
\begin{array}{l}
\text{for}\;n=0,1,\ldots\\
\left\lfloor
\begin{array}{l}
\text{for}\;i=1,\ldots,m\\
\left\lfloor
\begin{array}{l}
r_{i,n}=\varepsilon_{i,n}\Big(x_{i,n}-\\
\hskip 18mm \gamma_{i,n}\sum_{k=1}^q
L_{k,i}^*\Big(\nabla g_k\big(\sum_{j=1}^mL_{k,j}
x_{j,n}\big)\Big)\Big)\\[2mm]
x_{i,n+1}=x_{i,n}+\varepsilon_{i,n}\lambda_n
\big(\prox_{\gamma_{i,n}f_i}r_{i,n}-x_{i,n}\big),
\end{array}
\right.
\end{array}
\right.\\
\end{array}
\end{equation}
where $\gamma_{i,n}\in\RPP$ and $\lambda_n\in\rzeroun$. The
convergence of \eqref{e:bcFB} has been investigated in various
settings in terms of the expected value of
the cost function \cite{Neco16,Rich14,Rich16,Salz19}, the
mean square convergence of the iterates
\cite{Comb19,Rich14,Rich16}, or the almost sure convergence of the
iterates \cite{Comb15,Salz19}. It is shown in \cite{Salz19} that
algorithms such as the so-called \emph{random Kaczmarz method} 
to solve standard linear systems are special cases of
\eqref{e:bcFB}.

A noteworthy feature of the block-coordinate forward-backward
algorithm \eqref{e:bcFB} is that, at iteration $n$, it allows for
the use of distinct parameters $(\gamma_{i,n})_{1\leq i\leq m}$ to
update each component. This was observed to be beneficial to the
convergence profile in several applications \cite{Chou16,Rich14}.
See also \cite{Salz19} for further developments along these lines.

\subsection{Block-iterative multivariate minimization algorithms}

We investigate a specialization of a primal-dual version of the
multivariate inclusion Problem~\ref{prob:82} in the context of
Problem~\ref{prob:2}.

\begin{problem}
\label{prob:3}
Consider the setting of Problem~\ref{prob:2}.
The task is to solve the primal minimization problem
\begin{equation}
\label{e:12p}
\minimize{\boldsymbol{x}\in\HHH}
{\sum_{i=1}^mf_i(x_i)+\sum_{k=1}^q 
g_k\bigg(\sum_{i=1}^mL_{k,i}x_i\bigg)},
\end{equation}
along with its dual problem
\begin{equation}
\label{e:12d}
\minimize{\boldsymbol{v}^*\in\GGG}
{\sum_{i=1}^mf_i^*\bigg(-\sum_{k=1}^qL_{k,i}^*v^*_k\bigg)
+\sum_{k=1}^qg^*_k(v^*_k)}.
\end{equation}
\end{problem}

We solve Problem~\ref{prob:3} with algorithm \eqref{e:n03a} by
replacing $J_{\gamma_{i,n}A_i}$ by 
$\prox_{\gamma_{i,n}f_i}$ and $J_{\mu_{k,n}B_k}$ by 
$\prox_{\mu_{k,n}g_k}$. This block-iterative
method then produces a 
sequence $(\boldsymbol{x}_n)_{n\in\NN}$ which converges to a
solution to \eqref{e:12p} and a sequence 
$(\boldsymbol{v}^*_n)_{n\in\NN}$ which converges to a solution to
\eqref{e:12d} \cite{MaPr18}.

Examples of problems that conform to the format of 
Problems~\ref{prob:2} or \ref{prob:3} are 
encountered in image processing \cite{Berg16,Nmtm09,Jmiv11} as well
as in machine learning 
\cite{Argy12,Bach12,Nume19,Jaco09,Jena11,McDo16,Vill14,Yuan06}.

\subsection{Splitting based on Bregman distances}
The notion of a Bregman distance goes back to \cite{Breg67} and it
has been used since the 1980s in signal recovery; see
\cite{Byrn01,Cens97}. Let $\varphi\in\Gamma_0(\HH)$ be strictly
convex, and differentiable on $\intdom\varphi\neq\emp$ (more
precisely, we require a Legendre function, see \cite{Baus97,Sico03}
for the technical details). The associated \emph{Bregman distance}
between two points $x$ and $y$ in $\HH$ is 
\begin{equation}
\label{e:Df}
D_\varphi(x,y)=
\begin{cases}
\varphi(x)-\varphi(y)-\scal{x-y}{\nabla \varphi(y)},\\
\hskip 21mm \text{if}\;\;y\in\intdom \varphi;\\
\pinf,\hskip 13mm \text{otherwise}.
\end{cases}
\end{equation}
This construction captures many interesting discrepancy measures in
data analysis such as the Kullback-Leibler divergence. Another
noteworthy instance is when $\varphi=\|\cdot\|^2/2$, which yields
$D_\varphi(x,y)=\|x-y\|^2/2$ and suggests extending standard 
tools such as projection and proximity operators (see 
Theorems~\ref{t:11} and \ref{t:12}) by replacing the quadratic 
kernel by a Bregman distance 
\cite{Baus97,Sico03,Breg67,Cens92,Ecks93,Tebo92}. 
For instance, under mild conditions
on $f\in\Gamma_0(\HH)$ \cite{Sico03}, the 
\emph{Bregman proximal point} of
$y\in\intdom\varphi$ relative to $f$ is the unique point
$\prox^\varphi_fy$ which solves
\begin{equation}
\label{e:Dprox}
\minimize{p\in\intdom\varphi}{f(p)+D_\varphi(p,y)}.
\end{equation}
The \emph{Bregman projection} $\proj_C^\varphi y$ of $y$ onto 
a nonempty closed convex 
set $C$ in $\HH$ is obtained by setting $f=\iota_C$ above. Various
algorithms such as the POCS algorithm \eqref{e:g3} or the 
proximal point algorithm \eqref{e:sd2} have been extended in the
context of Bregman distances \cite{Baus97,Sico03}. For instance
\cite{Baus97} establishes the convergence to a solution to
Problem~\ref{prob:1} of a notable extension of POCS in which 
the sets are Bregman-projected onto in arbitrary order, namely
\begin{equation}
\label{e:h}
(\forall n\in\NN)\quad
x_{n+1}=\proj^\varphi_{C_{\mathrm{i}(n)}}x_n,
\end{equation}
where $\mathrm{i}\colon\NN\to\{1,\ldots,m\}$ is such that, for
every $p\in\NN$ and every $j\in\{1,\ldots,m\}$, there exists 
$n\geq p$ such that $\mathrm{i}(n)=j$.

A motivation for such extensions is that, for certain functions,
proximal points are easier to compute in the Bregman sense than in
the standard quadratic sense \cite{Baus17,Joca16,Nguy17}. Some work
has also focused on monotone operator splitting using Bregman
distances as an extension of standard methods \cite{Joca16}. The
Bregman version of the basic forward-backward minimization method
of Proposition~\ref{p:fb17}, namely,
\begin{equation}
\label{e:FB4}
\begin{array}{l}
\text{for}\;n=0,1,\ldots\\
\left\lfloor
\begin{array}{l}
u_n=\nabla\varphi(x_n)-\gamma_n\nabla g(x_n)\\
x_{n+1}=\big(\nabla\varphi+\gamma_n\partial f\big)^{-1}u_n
\end{array}
\right.\\[2mm]
\end{array}
\end{equation}
has also been investigated in \cite{Baus17,Buim20,Nguy17}
(note that the standard quadratic kernel corresponds 
to $\nabla\varphi=\Id$).
In these papers, it was shown to converge in instances when
\eqref{e:FB2} cannot be used because $\nabla g$ is not 
Lipschitzian.

\section{Fixed point modeling of Nash equilibria}
\label{sec:6}
In addition to the notation of Section~\ref{sec:not}, given
$i\in\{1,\ldots,m\}$, $x_i\in\HH_i$, and $\boldsymbol{y}\in\HHH$,
we set
\begin{equation}
\begin{cases}
\HHH_{\smallsetminus i}=
\HH_1\times\cdots\times\HH_{i-1}\times\HH_{i+1}\times
\cdots\times\HH_{m}\\
\boldsymbol{y}_{\smallsetminus i}=(y_j)_{1\leq j\leq m, j\neq i}\\
(x_i;\boldsymbol{y}_{\smallsetminus i})=(y_1,\ldots,y_{i-1},x_i,
y_{i+1},\ldots,y_m).
\end{cases}
\end{equation}
In various problems arising in signal recovery 
\cite{Aujo04,Aujo06,Berg16,Nmtm09,Jmiv11,Dani12,Darb20,Demo04}, 
telecommunications \cite{Lasa11,Scut10}, 
machine learning \cite{Brav18,Dasg19}, 
network science \cite{Yip19a,Yinh11}, 
and control \cite{Belg19,Borz13,Zhan19}, the solution is not
a single vector but a collections of vectors
$\boldsymbol{x}=(x_1,\ldots,x_m)\in\HHH$ representing the
actions of $m$ competing players. Oftentimes, such solutions 
cannot be modeled via a standard minimization problem of the form 
\begin{equation}
\label{e:s4}
\minimize{\boldsymbol{x}\in\HHH}{\boldsymbol{h}(\boldsymbol{x})}
\end{equation}
for some function $\boldsymbol{h}\colon\HHH\to\RX$, but rather as 
a Nash equilibrium \cite{Nash51}. 
In this game-theoretic setting \cite{Lara19}, 
player $i$ aims at minimizing his individual loss (or negative 
payoff) function $\boldsymbol{h}_i\colon\HHH\to\RX$, that
incorporates the actions of the other players. An action
profile $\overline{\boldsymbol{x}}\in\HHH$ is called a \emph{Nash
equilibrium} if unilateral deviations from it are not
profitable, i.e., 
\begin{equation}
\label{e:u4}
(\forall i\in\{1,\ldots,m\})\quad\boldsymbol{h}_i
(\overline{x}_i;\overline{\boldsymbol{x}}_{\smallsetminus i})
=\min_{x_i\in\HH_i}{\boldsymbol{h}_i
(x_i;\overline{\boldsymbol{x}}_{\smallsetminus i})}.
\end{equation}
In other words, if 
\begin{multline}
\label{e:best}
\!\!\best_i\colon\HHH_{\smallsetminus i}\to 2^{\HH_i}\colon
\boldsymbol{x}_{\smallsetminus i}\mapsto\\
\quad\menge{x_i\in\HH_i}{(\forall y_i\in\HH_i)\:
\boldsymbol{h}_i(y_i;\boldsymbol{x}_{\smallsetminus i})\geq 
\boldsymbol{h}_i(x_i;\boldsymbol{x}_{\smallsetminus i})}
\end{multline}
denotes the \emph{best response operator} of player $i$, 
$\overline{\boldsymbol{x}}\in\HHH$ is a Nash equilibrium if and
only if
\begin{equation}
\label{e:u2}
(\forall i\in\{1,\ldots,m\})\quad\overline{x}_i\in
\best_i(\overline{\boldsymbol{x}}_{\smallsetminus i}).
\end{equation}
This property can also be expressed in terms of the set-valued
operator
\begin{equation}
\label{ma719}
\boldsymbol{B}\colon\HHH\to 2^{\HHH}\colon \boldsymbol{x}\mapsto
\best_1(\boldsymbol{x}_{\smallsetminus 1})\times\cdots\times
\best_m(\boldsymbol{x}_{\smallsetminus m}).
\end{equation}
Thus, a point $\boldsymbol{\overline{x}}\in\HHH$ is a Nash
equilibrium if and only if it is a fixed point of $\boldsymbol{B}$
in the sense that $\overline{\boldsymbol{x}}\in\boldsymbol{B}
\overline{\boldsymbol{x}}$.

\subsection{Cycles in the POCS algorithm}
\label{sec:cycles}

\begin{figure}[h!tb]
\begin{center}
\scalebox{0.43}{
\begin{pspicture}(1,-4.85)(18.24,6.1)
\definecolor{color96b}{rgb}{0.9,0.90,1.0}
\definecolor{red1}{rgb}{0.0,0.60,0.6}
\definecolor{red2}{rgb}{0.90,0.0,0.0}
\rput{18.0}(-0.18953674,-2.8316424)%
{\psellipse[linewidth=0.06,dimen=outer,fillstyle=solid,%
fillcolor=color96b](8.844375,-2.0141652)(5.0,1.5)}
\psellipse[linewidth=0.06,dimen=outer,fillstyle=solid,%
fillcolor=color96b](11.444375,3.185835)(5.8,1.8)
\psellipse[linewidth=0.06,dimen=outer,fillstyle=solid,%
fillcolor=color96b](17.4,-2.18)(1.8,2.0)
\psline[linewidth=0.06cm,linestyle=solid,linecolor=red1,%
arrowsize=0.18cm 2.0,arrowlength=1.4,arrowinset=0.4]{->}%
(0.204375,0.20583488)(5.7243,2.926)
\psline[linewidth=0.06cm,linestyle=solid,linecolor=red1,%
arrowsize=0.18cm 2.0,arrowlength=1.4,arrowinset=0.4]{->}%
(5.7243,2.926)(7.304375,-0.994165)
\psline[linewidth=0.06cm,linestyle=solid,linecolor=red1,%
arrowsize=0.18cm 2.0,arrowlength=1.4,%
arrowinset=0.4]{->}(7.304375,-0.994165)(15.63,-1.9)
\psline[linewidth=0.06cm,linestyle=solid,linecolor=red1,%
arrowsize=0.18cm 2.0,arrowlength=1.4,%
arrowinset=0.4]{->}(15.63,-1.9)(14.70,1.74)
\psline[linewidth=0.06cm,linestyle=solid,linecolor=red1,%
arrowsize=0.18cm 2.0,arrowlength=1.4,%
arrowinset=0.4]{->}(14.70,1.74)(13.25,-0.11)
\psline[linewidth=0.06cm,arrowsize=0.18cm 2.0,arrowlength=1.4,%
arrowinset=0.4,linecolor=red2]{->}(15.1,1.83)(13.54,-0.29)
\psline[linewidth=0.06cm,arrowsize=0.18cm 2.0,arrowlength=1.4,%
arrowinset=0.4,linecolor=red2]{->}(13.5,-0.36)(15.75,-1.36)
\psline[linewidth=0.06cm,arrowsize=0.18cm 2.0,arrowlength=1.4,%
arrowinset=0.4,linecolor=red2]{->}(15.8,-1.36)(15.1,1.78)
\psdots[dotsize=0.18,linecolor=red1](0.24,0.22)
\rput(8.0,-2.8){\LARGE $C_2$}
\rput(10.8,3.3){\LARGE $C_3$}
\rput(18.0,-2.0){\LARGE $C_1$}
\rput(0.22,-0.20){\LARGE $x_0$}
\rput(7.5,-1.4){\LARGE $x_2$}
\rput(6.2,3.1){\LARGE $x_1$}
\rput(16.0,-2.2){\LARGE $x_3$}
\rput(14.3,2.04){\LARGE $x_4$}
\rput(12.9,0.3){\LARGE $x_5$}
\rput(13.1,-0.7){\LARGE $\overline{x}_2$}
\rput(15.4,2.37){\LARGE $\overline{x}_3$}
\rput(16.4,-1.33){\LARGE $\overline{x}_1$}
\psdots[dotsize=0.18](15.1,1.83)
\psdots[dotsize=0.18](13.5,-0.36)
\psdots[dotsize=0.18](15.8,-1.36)
\end{pspicture} 
}
\caption{The POCS algorithm with $m=3$ sets and 
initialized at $x_0$ produces the cycle 
$(\overline{x}_1,\overline{x}_2,\overline{x}_3)$.}
\label{fig:2}
\end{center}
\end{figure}
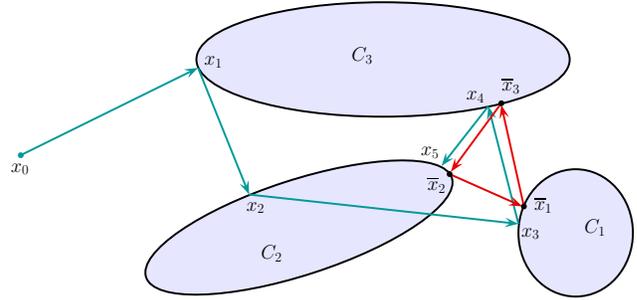

Let us go back to feasibility and Problem~\ref{prob:1}.
The POCS algorithm \eqref{e:g3} converges to a solution to the
feasibility problem \eqref{e:cfp1} when one exists. Now suppose
that Problem~\ref{prob:1} is inconsistent, with $C_1$ bounded. 
Then, as seen in Example~\ref{ex:20}, in the case of $m=2$ sets,
the sequence $(x_{2n})_{n\in\NN}$ produced by the alternating
projection algorithm \eqref{e:u8}, written as
\begin{equation}
\label{e:23}
\begin{array}{l}
\text{for}\;n=0,1,\ldots\\
\left\lfloor
\begin{array}{l}
x_{2n+1}=\proj_{C_2}x_{2n}\\
x_{2n+2}=\proj_{C_1}x_{2n+1},
\end{array}
\right.\\[2mm]
\end{array}
\end{equation}
converges to a point
$\overline{x}_1\in\Fix(\proj_{C_1}\circ\proj_{C_2})$, i.e.,
to a minimizer of $d_{C_2}$ over $C_1$. 
More precisely \cite{Che59a}, if we
set $\overline{x}_2=\proj_{C_2}\overline{x}_1$, then 
$\overline{x}_1=\proj_{C_1}\overline{x}_2$ and 
$(\overline{x}_1,\overline{x}_2)$ solves
\begin{equation}
\label{e:jeux1}
\minimize{x_1\in C_1,\,x_2\in C_2}{\|x_1-x_2\|}.
\end{equation}
An extension of the alternating projection method 
\eqref{e:23} to $m$ sets is the POCS algorithm \eqref{e:g3}, which
we write as
\begin{equation}
\label{e:pocs2}
\begin{array}{l}
\text{for}\;n=0,1,\ldots\\
\left\lfloor
\begin{array}{ll}
x_{mn+1}&\hskip -3mm=\proj_{C_m}x_{mn}\\
x_{mn+2}&\hskip -3mm=\proj_{C_{m-1}}x_{mn+1}\\
&\hskip -2mm\vdots\\
x_{mn+m}&\hskip -3mm=\proj_{C_1}x_{mn+m-1}.
\end{array}
\right.\\[2mm]
\end{array}
\end{equation}
As first shown in \cite{Gubi67} (this is also a
consequence of Theorem~\ref{t:5}),
for every $i\in\{1,\ldots,m\}$, $(x_{mn+i})_{n\in\NN}$
converges to a point $\overline{x}_{m+1-i}\in C_{m+1-i}$; in
addition $(\overline{x}_i)_{1\leq i\leq m}$ forms a \emph{cycle} in
the sense that (see Fig.~\ref{fig:2})
\begin{multline}
\label{e:21c}
\overline{x}_1=\proj_{C_1}\overline{x}_2,\;\ldots,\;
\overline{x}_{m-1}=\proj_{C_{m-1}}\overline{x}_m,\\
\text{and}\quad
\overline{x}_m=\proj_{C_m}\overline{x}_1.
\end{multline}
As shown in \cite{Baill12}, in stark contrast with the case of
$m=2$ sets and \eqref{e:jeux1}, there exists no function
$\Phi\colon\HH^m\to\RR$ such that cycles solve the minimization
problem 
\begin{equation} 
\label{e:best3} 
\minimize{x_1\in C_1,\ldots,\,x_m\in C_m}{\Phi(x_1,\ldots,x_m)},
\end{equation} 
which deprives cycles of a minimization interpretation. 
Nonetheless, cycles are equilibria in a more general sense, which
can be described from three different perspectives.
\begin{itemize}
\item
Fixed point theory:
Define two operators $\boldsymbol{P}$ and $\boldsymbol{L}$ from
$\HH^m$ to $\HH^m$ by
\begin{equation}
\begin{cases}
\boldsymbol{P}\colon\boldsymbol{x}\mapsto
(\proj_{C_1}x_1,\ldots,\proj_{C_m}x_m)\\
\boldsymbol{L}\colon\boldsymbol{x}\mapsto
(x_2,\ldots,x_m,x_1).
\end{cases}
\end{equation}
Then, in view of \eqref{e:21c}, the set of cycles is precisely the
set of fixed points of $\boldsymbol{P\circ L}$, which is also the 
set of fixed points of $\boldsymbol{T}=\boldsymbol{P\circ F}$,
where $\boldsymbol{F}=(\ID+\boldsymbol{L})/2$ (see
\cite[Corollary~26.3]{Livre1}). Since 
Example~\ref{ex:1} implies that $\boldsymbol{P}$ is firmly
nonexpansive and since $\boldsymbol{L}$ is nonexpansive,
$\boldsymbol{F}$ is firmly nonexpansive as well. It thus follows
from Example~\ref{ex:roma}, that the cycles are the fixed points of
the $2/3$-averaged operator $\boldsymbol{T}$. 
\item
Game theory: Consider a game in $\HH^m$ in which the goal of player
$i$ is to minimize the loss 
\begin{equation}
\label{e:u3}
\boldsymbol{h}_i
\colon(x_i;\boldsymbol{x}_{\smallsetminus i})\mapsto
\iota_{C_i}(x_i)+\frac{1}{2}\|x_i-x_{i+1}\|^2,
\end{equation}
i.e., to be in $C_i$ and as close as possible to the action 
of player $i+1$ (with the convention $x_{m+1}=x_1$). Then a cycle
$(\overline{x}_1,\ldots,\overline{x}_m)$ is a solution to
\eqref{e:u4} and therefore a Nash equilibrium. 
Let us note that the best response operator of player $i$ is
$\best_i\colon\boldsymbol{x}_{\smallsetminus i}\mapsto
\proj_{C_i}x_{i+1}$.
\item
Monotone inclusion: Applying Fermat's rule to each line of
\eqref{e:u4} in the setting of \eqref{e:u3}, and using
\eqref{e:normalcone}, we obtain
\begin{equation}
\label{e:as}
\begin{cases}
0\in N_{C_1}\overline{x}_1+\overline{x}_1-\overline{x}_2\\
\hskip 4mm\vdots\\
0\in N_{C_{m-1}}\overline{x}_{m-1}+\overline{x}_{m-1}
-\overline{x}_m\\
0\in N_{C_m}\overline{x}_m+\overline{x}_m
-\overline{x}_1.
\end{cases}
\end{equation}
In terms of the maximally monotone operator 
$\boldsymbol{A}=N_{C_1\times\cdots\times C_m}$ and the cocoercive
operator 
\begin{equation}
\boldsymbol{B}\colon\boldsymbol{x}\mapsto
(x_1-x_2,\ldots,{x}_{m-1}-{x}_m,x_m-x_1),
\end{equation}
\eqref{e:as} can be rewritten as an instance of
Problem~\ref{prob:7} in $\HH^m$, namely, 
$\boldsymbol{0}\in\boldsymbol{A}\overline{\boldsymbol{x}}+
\boldsymbol{B}\overline{\boldsymbol{x}}$.
\end{itemize}

\subsection{Proximal cycles}
\label{sec:IIC}
We have seen in Section~\ref{sec:cycles} a first example of a Nash
equilibrium. This setting can be extended by replacing the
indicator function $\iota_{C_i}$ in \eqref{e:u3} by a
general function $\varphi_i\in\Gamma_0(\HH)$ modeling the 
self-loss of player $i$, i.e., 
\begin{equation}
\label{e:u5}
\boldsymbol{h}_i
\colon(x_i;\boldsymbol{x}_{\smallsetminus i})\mapsto
\varphi_i(x_i)+\frac{1}{2}\|x_i-x_{i+1}\|^2.
\end{equation}
The solutions to the resulting problem \eqref{e:u4} are 
\emph{proximal cycles}, i.e., 
$m$-tuples $(\overline{x}_i)_{1\leq i\leq m}\in\HH^m$ such that 
\begin{multline}
\label{e:21d}
\overline{x}_1=\prox_{\varphi_1}\overline{x}_2,\;\ldots,\;
\overline{x}_{m-1}=\prox_{\varphi_{m-1}}\overline{x}_m,\\
\text{and}\quad
\overline{x}_m=\prox_{\varphi_m}\overline{x}_1.
\end{multline}
Furthermore, the equivalent monotone inclusion and fixed point 
representations of the cycles in Section~\ref{sec:cycles} 
remain true with
\begin{equation}
\boldsymbol{P}\colon\HHH\to\HHH\colon\boldsymbol{x}\mapsto
\big(\prox_{\varphi_1}x_1,\ldots,\prox_{\varphi_m}x_m\big)
\end{equation}
and $\boldsymbol{A}=\partial\boldsymbol{f}$, where
$\boldsymbol{f}\colon\boldsymbol{x}\mapsto
\sum_{i=1}^m\varphi_i(x_i)$.
Here, the best response operator of player $i$ is
$\best_i\colon\boldsymbol{x}_{\smallsetminus i}
\mapsto\prox_{\varphi_i}x_{i+1}$.
Examples of such cycles appear in \cite{Nmtm09,Smms05}.

\subsection{Construction of Nash equilibria}

A more structured version of the Nash equilibrium formulation
\eqref{e:u4}, which captures \eqref{e:u5} and therefore 
\eqref{e:u3}, is provided next.

\begin{problem}
\label{prob:j}
For every $i\in\{1,\ldots,m\}$, let $\psi_i\in\Gamma_0(\HH_i)$,
let $\boldsymbol{f}_i\colon\HHH\to\RX$, let
$\boldsymbol{g}_i\colon\HHH\to\RX$
be such that, for every $\boldsymbol{x}\in\HHH$, 
$\boldsymbol{f}_i
(\cdot;\boldsymbol{x}_{\smallsetminus i})\in\Gamma_0(\HH_i)$ and
$\boldsymbol{g}_i
(\cdot;\boldsymbol{x}_{\smallsetminus i})\in\Gamma_0(\HH_i)$.
The task is to 
\begin{multline}
\label{e:u7}
\text{find}\;\;\overline{\boldsymbol{x}}\in\HHH
\quad\text{such that}\quad
(\forall i\in\{1,\ldots,m\})\\
\overline{x}_i\in\Argmind{x_i\in\HH_i}{\psi_i(x_i)+
\boldsymbol{f}_i(x_i;\overline{\boldsymbol{x}}_{\smallsetminus i})+
\boldsymbol{g}_i
(x_i;\overline{\boldsymbol{x}}_{\smallsetminus i})}.
\end{multline}
\end{problem}

Under suitable assumptions on 
$(\boldsymbol{f}_i)_{1\leq i\leq m}$ and
$(\boldsymbol{g}_i)_{1\leq i\leq m}$, monotone operator
splitting strategies can be contemplated to solve 
Problem~\ref{prob:j}. This approach was
initiated in \cite{Cohe87} in a special case of the following
setting, which reduces to that investigated in \cite{Bric13} when
$(\forall i\in\{1,\ldots,m\})$ $\psi_i=0$.

\begin{assumption}
\label{a:j1}
In Problem~\ref{prob:j},
the functions $(\boldsymbol{f}_i)_{1\leq i\leq m}$ coincide with
a function $\boldsymbol{f}\in\Gamma_0(\HHH)$. For every
$i\in\{1,\ldots,m\}$ and every $\boldsymbol{x}\in\HHH$, 
$\boldsymbol{g}_i(\cdot;\boldsymbol{x}_{\smallsetminus i})$ is
differentiable on $\HH_i$ and
$\nabli{i}\boldsymbol{g}_i(\boldsymbol{x})$ denotes its derivative
relative to $x_i$. Moreover, 
\begin{multline}
\label{e:Cmon}
(\forall\boldsymbol{x}\in\HHH)(\forall\boldsymbol{y}\in\HHH)\\
\sum_{i=1}^m\scal{\nabli{i}\boldsymbol{g}_i
(\boldsymbol{x})-\nabli{i}\boldsymbol{g}_i
(\boldsymbol{y})}{x_i-y_i}\geq 0,
\end{multline} 
and 
\begin{multline}
(\exi\boldsymbol{z}\in\HHH)\quad -\big(\nabli{1}\boldsymbol{g}_1
(\boldsymbol{z}),\ldots,\nabli{m}\boldsymbol{g}_m
(\boldsymbol{z})\big)\\
\in\partial\boldsymbol{f}(\boldsymbol{z})+
\overset{m}{\underset{i=1}{\cart}}\partial\psi_i(z_i).
\end{multline}
\end{assumption}

In the context of Assumption~\ref{a:j1}, let us introduce the
maximally monotone operators on $\HHH$
\begin{equation}
\label{e:B}
\begin{cases}
\boldsymbol{A}=\partial\boldsymbol{f}\\
\boldsymbol{B}\colon\boldsymbol{x}\mapsto\cart_{\!i=1}^{\!m}
\partial\psi_i(x_i)\\
\boldsymbol{C}\colon
\boldsymbol{x}\mapsto\big(\nabli{1}\boldsymbol{g}_1
(\boldsymbol{x}),\ldots,\nabli{m}\boldsymbol{g}_m
(\boldsymbol{x})\big).
\end{cases}
\end{equation}
Then the solutions to the inclusion problem 
(see Problem~\ref{prob:8}) $\boldsymbol{0}\in\boldsymbol{A}
\boldsymbol{x}+\boldsymbol{B}\boldsymbol{x}+
\boldsymbol{C}\boldsymbol{x}$
solve Problem~\ref{prob:j} \cite{Bric13}. In turn, applying the
splitting scheme of Proposition~\ref{p:71} leads
to the following implementation.

\begin{proposition}
\label{p:n1}
Consider the setting of Assumption~\ref{a:j1} with the additional
requirement that, for some $\delta\in\RPP$,
\begin{multline}
\label{e:CLip}
(\forall\boldsymbol{x}\in\HHH)
(\forall\boldsymbol{y}\in\HHH)\quad
\sum_{i=1}^m\|\nabli{i}\boldsymbol{g}_i
(\boldsymbol{x})-\nabli{i}\boldsymbol{g}_i
(\boldsymbol{y})\|^2\\
\leq\delta^2\sum_{i=1}^m\|x_i-y_i\|^2.
\end{multline}
Let $\varepsilon\in\left]0,1/(2+\delta)\right[$,
let $(\gamma_n)_{n\in\NN}$ be in 
$\left[\varepsilon,(1-\varepsilon)/(1+\delta)\right]$,
let $\boldsymbol{x}_0\in\HHH$, and let $\boldsymbol{v}_0\in\HHH$.
Iterate
\begin{equation}
\label{e:fz79}
\begin{array}{l}
\text{for}\;n=0,1,\ldots\\
\left\lfloor 
\begin{array}{l}
\text{for}\:\:i=1,\ldots,m\\
\lfloor\:y_{i,n}=x_{i,n}-\gamma_n\big(\nabli{i}\boldsymbol{g}_i
(\boldsymbol{x}_n)+v_{i,n}\big)\\
\boldsymbol{p}_n=\prox_{\gamma_n\boldsymbol{f}}\:
\boldsymbol{y}_n\\
\text{for}\:\:i=1,\ldots,m\\
\left\lfloor
\begin{array}{l}
q_{i,n}=v_{i,n}+\gamma_n\big(x_{i,n}-\prox_{\psi_i/\gamma_n}
(v_{i,n}/\gamma_n+x_{i,n})\big)\\
x_{i,n+1}=x_{i,n}-y_{i,n}+p_{i,n}-\gamma_n
\big(\nabli{i}\boldsymbol{g}_i(\boldsymbol{p}_n)+q_{i,n}\big)\\
v_{i,n+1}=q_{i,n}+\gamma_n(p_{i,n}-x_{i,n}).
\end{array}
\right.
\end{array}
\right.
\end{array}
\end{equation}
Then there exists a solution $\overline{\boldsymbol{x}}$ to
Problem~\ref{prob:j} such that, for every $i\in\{1,\ldots,m\}$,
$x_{i,n}\to\overline{x}_i$.
\end{proposition}

\begin{example}
\label{ex:juil1}
Let $\varphi_1\colon\HH_1\to\RR$ be convex and differentiable with
a $\delta_1$-Lipschitzian gradient, 
let $\varphi_2\colon\HH_2\to\RR$ be convex and differentiable with
a $\delta_2$-Lipschitzian gradient, let
$L\colon\HH_1\to\HH_2$ be linear, and let $C_1\subset\HH_1$,
$C_2\subset\HH_2$, and $\boldsymbol{D}\subset\HH_1\times\HH_2$ be
nonempty closed convex sets. Suppose that there exists 
$\boldsymbol{z}\in\HH_1\times\HH_2$ such that
$-(\nabla\varphi_1(z_1)+L^*z_2,\nabla\varphi_2(z_2)-Lz_1)
\in N_{\boldsymbol{D}}(z_1,z_2)+N_{C_1}z_1\times N_{C_2}z_2$.
Then the 2-player game 
\begin{equation}
\label{e:j23}
\begin{cases}
\overline{x}_1\in\Argmind{x_1\in C_1}
{\iota_{\boldsymbol{D}}(x_1,\overline{x}_2)+
\varphi_1(x_1)+\scal{Lx_1}{\overline{x}_2}}\\
\overline{x}_2\in\Argmind{x_2\in C_2}
{\iota_{\boldsymbol{D}}(\overline{x}_1,x_2)+
\varphi_2(x_2)-\scal{L\overline{x}_1}{x_2}}
\end{cases}
\end{equation}
is an instance of Problem~\ref{prob:j} with
$\boldsymbol{f}_1=\boldsymbol{f}_2=\iota_{\boldsymbol{D}}$, 
$\psi_1=\iota_{C_1}$, $\psi_2=\iota_{C_2}$, and 
\begin{equation}
\label{e:juil2}
\begin{cases}
\boldsymbol{g}_1\colon(x_1,x_2)\mapsto
\varphi_1(x_1)+\scal{Lx_1}{x_2}\\
\boldsymbol{g}_2\colon(x_1,x_2)\mapsto
\varphi_2(x_2)-\scal{Lx_1}{x_2}.
\end{cases}
\end{equation}
In addition, Assumption~\ref{a:j1} is satisfied, as well as 
\eqref{e:CLip} with $\delta=\max\{\delta_1,\delta_2\}+\|L\|$.
Moreover, in view of \eqref{e:98}, algorithm \eqref{e:fz79}
becomes
\begin{equation}
\label{e:fz80}
\hskip -3mm
\begin{array}{l}
\text{for}\;n=0,1,\ldots\\
\left\lfloor 
\begin{array}{l}
y_{1,n}=x_{1,n}-\gamma_n\big(\nabla\varphi_1(x_{1,n})+L^*x_{2,n}
+v_{1,n}\big)\\
y_{2,n}=x_{2,n}-\gamma_n\big(\nabla\varphi_2(x_{2,n})-Lx_{1,n}
+v_{2,n}\big)\\
\boldsymbol{p}_n=\proj_{\boldsymbol{D}}\:
\boldsymbol{y}_n\\
q_{1,n}=v_{1,n}+\gamma_n\big(x_{1,n}-\proj_{C_1}
(v_{1,n}/\gamma_n+x_{1,n})\big)\\
q_{2,n}=v_{2,n}+\gamma_n\big(x_{2,n}-\proj_{C_2}
(v_{2,n}/\gamma_n+x_{2,n})\big)\\
x_{1,n+1}=x_{1,n}-y_{1,n}+p_{1,n}\\
\hskip 15mm -\gamma_n
\big(\nabla\varphi_1(p_{1,n})+L^*p_{2,n}+q_{1,n}\big)\\
x_{2,n+1}=x_{2,n}-y_{2,n}+p_{2,n}\\
\hskip 15mm -\gamma_n
\big(\nabla\varphi_2(p_{2,n})-Lp_{1,n}+q_{2,n}\big)\\
v_{1,n+1}=q_{1,n}+\gamma_n(p_{1,n}-x_{1,n})\\
v_{2,n+1}=q_{2,n}+\gamma_n(p_{2,n}-x_{2,n}).
\end{array}
\right.
\end{array}
\end{equation}
\end{example}

Condition \eqref{e:CLip} means that the operator $\boldsymbol{C}$
of \eqref{e:B} is $\delta$-Lipschitzian. The stronger assumption
that it is cocoercive, allows us to bring into play the
three-operator splitting algorithm of Proposition~\ref{p:72} to
solve Problem~\ref{prob:j}.

\begin{proposition}
\label{p:n2}
Consider the setting of Assumption~\ref{a:j1} with the additional
requirement that, for some $\beta\in\RPP$,
\begin{multline}
\label{e:Ccoco}
(\forall\boldsymbol{x}\in\HHH)
(\forall\boldsymbol{y}\in\HHH)\quad\sum_{i=1}^m
\scal{x_i-y_i}{\nabli{i}\boldsymbol{g}_i(\boldsymbol{x})-
\nabli{i}\boldsymbol{g}_i(\boldsymbol{y})}\\
\geq\beta\sum_{i=1}^m\|\nabli{i}\boldsymbol{g}_i
(\boldsymbol{x})-\nabli{i}\boldsymbol{g}_i
(\boldsymbol{y})\|^2.
\end{multline}
Let $\gamma\in\left]0,2\beta\right[$ and set 
$\alpha=2\beta/(4\beta-\gamma)$.
Furthermore, let $(\lambda_n)_{n\in\NN}$ be an
$\alpha$-relaxation sequence and let $\boldsymbol{y}_0\in\HHH$. 
Iterate
\begin{equation}
\label{e:j9}
\begin{array}{l}
\text{for}\;n=0,1,\ldots\\
\left\lfloor
\begin{array}{l}
\text{for}\:\:i=1,\ldots,m\\
\left\lfloor 
\begin{array}{l}
x_{i,n}=\prox_{\gamma\psi_i}\,y_{i,n}\\
r_{i,n}=y_{i,n}+\gamma \nabli{i}\boldsymbol{g}_i
(\boldsymbol{x}_n)\\
\end{array}
\right.\\[2mm]
\boldsymbol{z}_n=\prox_{\gamma\boldsymbol{f}}
(2\boldsymbol{x_n}-\boldsymbol{r}_n)\\
\boldsymbol{y}_{n+1}=\boldsymbol{y}_n+
\lambda_n(\boldsymbol{z}_n-\boldsymbol{x}_n).
\end{array}
\right.\\[2mm]
\end{array}
\end{equation}
Then there exists a solution $\overline{\boldsymbol{x}}$ to
Problem~\ref{prob:j} such that, for every $i\in\{1,\ldots,m\}$,
$x_{i,n}\to\overline{x}_i$.
\end{proposition}

\begin{example}
\label{ex:n45}
For every $i\in\{1,\ldots,m\}$,
let $C_i\subset\HH_i$ be a nonempty closed convex set, let 
$L_i\colon\HH_i\to\GG$ be linear, and let $o_i\in\GG$. The 
task is to solve the Nash equilibrium (with the convention
$L_{m+1}\overline{x}_{m+1}=L_1\overline{x}_1$)
\begin{multline}
\label{e:j8}
\text{find}\;\;\overline{\boldsymbol{x}}\in\HHH
\quad\text{such that}\quad
(\forall i\in\{1,\ldots,m\})\\
\overline{x}_i\in\Argmind{x_i\in C_i}
{\psi_i(x_i)+\dfrac{\|L_ix_i+L_{i+1}\overline{x}_{i+1}
-o_i\|^2}{2}}.
\end{multline}
Here, the action of player $i$ must lie in $C_i$, and it is further
penalized by $\psi_i$ and the proximity of the linear mixture
$L_ix_i+L_{i+1}\overline{x}_{i+1}$ to some vector $o_i$. For
instance if, for every $i\in\{1,\ldots,m\}$, $C_i=\HH_i$, $o_i=0$,
and $L_i=(-1)^i\Id$, we recover the setting of
Section~\ref{sec:IIC}. The equilibrium \eqref{e:j8} is an
instantiation of Problem~\ref{prob:j} with
$\boldsymbol{f}_1=\boldsymbol{f}_2\colon\boldsymbol{x}
\mapsto\sum_{i=1}^m\iota_{C_i}(x_i)$ and, for every
$i\in\{1,\ldots,m\}$, $\boldsymbol{g}_i\colon\boldsymbol{x}\mapsto
\|L_ix_i+L_{i+1}x_{i+1}-o_i\|^2/2$. In addition, as in
\cite[Section~9.4.3]{Bric13}, \eqref{e:Ccoco} holds with
$\beta=(2\max_{1\leq i\leq m}\|L_i\|^2)^{-1}$. Finally,
\eqref{e:j9} reduces to (with the convention
$L_{m+1}x_{m+1,n}=L_1x_{1,n}$)
\begin{equation}
\label{e:j6}
\hskip -2mm
\begin{array}{l}
\text{for}\;n=0,1,\ldots\\
\left\lfloor
\hskip -1mm
\begin{array}{l}
\text{for}\:\:i=1,\ldots,m\\
\left\lfloor 
\begin{array}{l}
\hskip -1mm 
x_{i,n}=\prox_{\gamma\psi_i}\,y_{i,n}\\
\hskip -1mm 
r_{i,n}=y_{i,n}+\gamma L_i^*(L_ix_{i,n}+L_{i+1}x_{i+1,n}-o_i)\\
\hskip -1mm 
z_{i,n}=\proj_{C_i}(2x_{i,n}-r_{i,n})\\
y_{i,n+1}=y_{i,n}+\lambda_n(z_{i,n}-x_{i,n}).
\end{array}
\right.\\[2mm]
\end{array}
\right.\\[2mm]
\end{array}
\end{equation}
\end{example}

\begin{remark}\
\label{r:j1}
\begin{enumerate}
\item 
As seen in Example~\ref{ex:juil1}, the functions of 
\eqref{e:juil2} satisfy the
Lipschitz condition \eqref{e:CLip}. However the cocoercivity
condition \eqref{e:Ccoco} does not hold. For instance, if
$\varphi_1=0$ and $\varphi_2=0$ then, for every $\boldsymbol{x}$
and $\boldsymbol{y}$ in $\HH_1\times\HH_2$,
\begin{multline}
\label{e:211}
\scal{\nabli{1}\boldsymbol{g}_1
(\boldsymbol{x})-\nabli{1}\boldsymbol{g}_1
(\boldsymbol{y})}{x_1-y_1}\hspace{2cm}\\+
\scal{\nabli{2}\boldsymbol{g}_2
(\boldsymbol{x})-\nabli{2}\boldsymbol{g}_2
(\boldsymbol{y})}{x_2-y_2}=0.
\end{multline}
\item 
Distributed splitting algorithms for finding Nash equilibria
are discussed in \cite{Belg18,Belg19,Yip19a,Yip19b}.
\item 
An asynchronous block-iterative decomposition algorithm to solve
Nash equilibrium problems involving a mix of nonsmooth and smooth
functions acting on linear mixtures of actions is proposed in
\cite{Nash21}.
\end{enumerate}
\end{remark}

\section{Fixed point modeling of other non-minimization problems}
\label{sec:7}

\subsection{Neural network structures}

\begin{figure}
\scalebox{0.50} 
{
\begin{pspicture}(-0.05,-1.7)(15.9,2.1)
\psline[linewidth=0.04cm,arrowsize=2.2mm]{->}(0.35,0.0)(1.0,0.0)
\psline[linewidth=0.04cm,arrowsize=2.2mm]{->}(3.0,0.0)(3.65,0.0)
\psline[linewidth=0.04cm,arrowsize=2.2mm]{->}(4.35,0.0)(5.0,0.0)
\psline[linewidth=0.04cm,arrowsize=2.2mm]{->}(4.0,1.2)(4.0,0.36)
\psframe[linewidth=0.04,dimen=outer](1.0,-1.0)(3.0,1.0)
\pscircle[linewidth=0.04,dimen=outer](4.0,0.0){0.35}
\rput(0.15,0.0){\Large$\boldsymbol{{x}}$}
\rput(2.0,0.0){\Large$\boldsymbol{{W_1}}$}
\rput(4.0,0.0){\Large$\boldsymbol{+}$}
\rput(4.0,1.5){\Large$\boldsymbol{{b_1}}$}
\rput(6.0,0.0){\Large$\boldsymbol{{R_1}}$}
\psframe[linewidth=0.04,dimen=outer](5.0,-1.0)(7.0,1.0)
\psline[linewidth=0.04cm,arrowsize=2.2mm]{->}(7.00,0.0)(7.65,0.0)
\rput(8.5,0.0){\Large$\boldsymbol{{\cdots}}$}
\psline[linewidth=0.04cm,arrowsize=2.2mm]{->}(9.35,0.0)(10.0,0.0)
\psline[linewidth=0.04cm,arrowsize=2.2mm]{->}(12.0,0.0)(12.65,0.0)
\psline[linewidth=0.04cm,arrowsize=2.2mm]{->}(13.35,0.0)(14.0,0.0)
\psline[linewidth=0.04cm,arrowsize=2.2mm]{->}(13.0,1.2)(13.0,0.36)
\psframe[linewidth=0.04,dimen=outer](10.0,-1.0)(12.0,1.0)
\pscircle[linewidth=0.04,dimen=outer](13.0,0.0){0.35}
\rput(11.0,0.0){\Large$\boldsymbol{{W_m}}$}
\rput(13.0,0.0){\Large$\boldsymbol{+}$}
\rput(13.0,1.5){\Large$\boldsymbol{{b_m}}$}
\rput(15.0,0.0){\Large$\boldsymbol{{R_m}}$}
\rput(17.0,0.05){\Large$\boldsymbol{{Tx}}$}
\psframe[linewidth=0.04,dimen=outer](14.0,-1.0)(16.0,1.0)
\psline[linewidth=0.04cm,arrowsize=2.2mm]{->}(16.00,0.0)(16.65,0.0)
\end{pspicture} 
}
\vskip -1mm
\caption{Feedforward neural network: the $i$th layer involves a
linear weight operator $W_i$, a bias vector $b_i$, and an 
activation operator $R_i$, which is assumed to be an
averaged nonexpansive operator.}
\label{fig:4}
\end{figure}
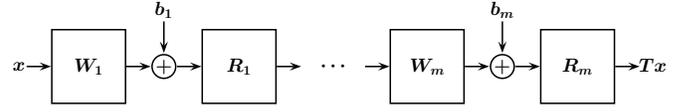

A feedforward neural network (see Fig.~\ref{fig:4}) consists 
of the composition of nonlinear activation operators and affine 
operators. More precisely, such an $m$-layer network can be
modeled as
\begin{equation}
\label{e:zib}
T=T_m\circ\cdots\circ T_1,
\end{equation}
where $T_i=R_i\circ(W_i\cdot +\,b_i)$, 
with $W_i\in\RR^{N_i\times N_{i-1}}$, $b_i\in\RR^{N_i}$, and
$R_i\colon\RR^{N_i}\to\RR^{N_i}$ (see Fig.~\ref{fig:4}).
If the $i$-th layer is convolutional, then the corresponding
weight matrix $W_i$ has a Toeplitz (or block-Toeplitz) structure.
Many common activation operators are separable, i.e.,
\begin{equation}
\label{e:actsep}
R_i\colon(\xi_k)_{1\leq k\leq N_i}\mapsto
\big(\varrho_{i,k}(\xi_k)\big)_{1\leq k\leq N_i},
\end{equation}
where $\varrho_{i,k}\colon\RR\to\RR$.
For example, the ReLU activation function is given by
\begin{equation}
\label{e:ReLU}
\varrho_{i,k}\colon\xi\mapsto
\begin{cases}
\xi,&\text{if}\;\;\xi>0;\\
0,&\text{if}\;\;\xi\leq 0,
\end{cases}
\end{equation}
and the unimodal sigmoid activation function is 
\begin{equation}
\label{e:sig}
\varrho_{i,k}\colon\xi\mapsto
\frac{1}{1+e^{-\xi}}-\frac{1}{2}.
\end{equation}
An example of a nonseparable operator is the softmax activator
\begin{equation}
R_i\colon (\xi_k)_{1\leq k\leq N_i}\mapsto
\left(e^{\xi_k}\left/ \displaystyle\sum_{j=1}^{N_i} 
e^{\xi_j}\right. \right)_{1\leq k\leq N_i}.
\end{equation}
It was observed in \cite{Smds20} that almost all standard 
activators are actually averaged operators in the sense of 
\eqref{e:105-}. In particular, as discussed in \cite{Svva20}, many 
activators are proximity operators in the sense of 
Theorem~\ref{t:12}. In this case, in \eqref{e:actsep}, there exist
functions $(\phi_k)_{1\leq k\leq N_i}$ in $\Gamma_0(\RR)$ such that
\begin{equation}
R_i\colon(\xi_k)_{1\leq k\leq N_i}\mapsto
\big(\prox_{\phi_k}\xi_k\big)_{1\leq k\leq N_i}.
\end{equation}
For ReLU, $\phi_k$ reduces to $\iota_{[0,+\infty[}$
whereas, for the unimodal sigmoid, it is the function 
\begin{equation}
\xi\mapsto
\begin{cases}
(\xi+1/2)\ln(\xi+1/2)+(1/2-\xi)\ln(1/2-\xi)\\
\hskip 13mm -(|\xi|^2+1/4)/2,
\hskip 7mm\text{if}\;\;|\xi|<1/2;\\
-1/4,\hskip 36mm\text{if}\;\;|\xi|=1/2;\\
\pinf,\hskip 37.8mm\text{if}\;\;|\xi|>1/2.
\end{cases}
\end{equation}
For softmax, we have $R_i=\prox_{\varphi_i}$ where
\begin{equation}
\varphi_i\colon(\xi_k)_{1\leq k\leq N_i}\mapsto
\begin{cases}
\sum_{i=1}^{N_i}\big(\xi_k\ln\xi_k-
{|\xi_k|^2}/{2}\big),\\
\qquad\text{if}\;\;
\displaystyle{\min_{1\leq k\leq N_i}}\xi_k\geq 0\;\text{and}\;
\sum_{k=1}^{N_i}\xi_k=1;\\
\pinf,\;\text{otherwise.}
\end{cases}
\end{equation}

The weight matrices $(W_i)_{1\leq i\leq m}$ play a
crucial role in the overall nonexpansiveness of the network. 
Indeed, under suitable conditions on these matrices, the network
$T$ is averaged. For example, let $W=W_m\cdots W_1$ and let
\begin{multline}
\label{e:defthetam}
\theta_m=\|W\|+\sum_{\ell=1}^{m-1}\sum_{0\leq
j_1<\cdots<j_\ell\leq m-1}\|W_m\cdots W_{j_\ell+1}\|\\
\times\|W_{j_\ell}\cdots
W_{j_{\ell-1}+1}\|\cdots\|W_{j_1}\cdots W_0\|.
\end{multline}
Then, if there exists $\alpha\in [1/2,1]$ such that
\begin{equation}
\label{e:alphaNN}
\|W-2^m(1-\alpha)\Id\|-\|W\|+2\theta_m\leq 2^m\alpha, 
\end{equation}
$T$ is $\alpha$-averaged. 
Other sufficient conditions have been established in \cite{Svva20}.
These results pave the way to a theoretical analysis of neural
networks from the standpoint of fixed point methods. In particular,
assume that $N_m=N_0$ and consider a recurrent network of the form
\begin{equation}
(\forall n\in\NN)\quad
x_{n+1}=(1-\lambda_n)x_n+\lambda_n T x_n,
\end{equation}
where $\lambda_n\in\RPP$ models
a skip connection. Then, according to
Theorem~\ref{t:1}, the convergence of $(x_n)_{n\in\NN}$ to a fixed
point of $T$ is guaranteed under condition~\eqref{e:alphaNN}
provided that $(\lambda_n)_{n\in\NN}$ is an $\alpha$-relaxation
sequence.
As shown in \cite{Svva20}, when for every $i\in\{1,\ldots,m\}$,
$R_i$ is the proximity operator of some function 
$\varphi_i\in\Gamma_0(\RR^{N_i})$, the recurrent network 
delivers asymptotically a solution to the system of inclusions
\begin{equation}
\label{e:vi1}
\begin{cases}
b_1\in&\hskip -2mm
\overline{x}_1-W_1\overline{x}_m+\partial
\varphi_1(\overline{x}_1)\\
b_2\in&\hskip -2mm\overline{x}_2-W_2\overline{x}_1+
\partial\varphi_2(\overline{x}_2)\\
\hskip 6mm\vdots\\
b_m\in&\hskip -2mm\overline{x}_m-W_m\overline{x}_{m-1}+
\partial\varphi_m(\overline{x}_m),
\end{cases}
\end{equation}
where $\overline{x}_m\in\Fix T$ and, for every 
$i\in\{2,\ldots,m\}$, $\overline{x}_i=T_i\overline{x}_{i-1}$.
Alternatively, \eqref{e:vi1} is a Nash equilibrium of
the form \eqref{e:u4} where
(we set $\overline{x}_0=\overline{x}_m$)
\begin{equation}
\label{e:u9}
\boldsymbol{h}_i\colon(x_i;
\overline{\boldsymbol{x}}_{\smallsetminus i})\mapsto
\varphi_i(x_i)+\dfrac{1}{2}\|x_i-b_i-W_i\overline{x}_{i-1}\|^2.
\end{equation}

Fixed point theory also allows us to provide conditions for $T$ to
be Lipschitzian and to calculate an associated Lipschitz constant.
Such results are useful to evaluate the robustness of the network
to adversarial perturbations of its input \cite{Szeg13}.
As shown in \cite{Smds20}, if $\theta_m$ is given by 
\eqref{e:defthetam},
$\theta_m/2^{m-1}$ is a Lipschitz constant of $T$ and
\begin{equation}
\label{e:sandLip}
\|W\|\leq\frac{\theta_m}{2^{m-1}}\leq\|W_1\|\cdots\|W_m\|.
\end{equation}
This bound is thus more accurate than the product of the individual
bounds corresponding to each layer used in \cite{Szeg13}. Tighter
estimations can also be derived, especially when the activation
operators are separable \cite{Smds20,Lato20,Scam18}. Note that the
lower bound in \eqref{e:sandLip} would correspond to a linear
network where all the nonlinear activation operators would be
removed. Interestingly, when all the weight matrices have
components in $\RP$, $\|W\|$ is a Lipschitz constant of the
network \cite{Smds20}.

Special cases of the neural network model of \cite{Svva20} are
investigated in \cite{Hasa20,Tang20}. Another special case of
interest is when the operator $T$ in \eqref{e:zib} corresponds
to the \emph{unrolling} (or \emph{unfolding}) of a fixed point
algorithm \cite{Mong21}, that is, each operator $T_i$ 
corresponds to one iteration of such an algorithm
\cite{Bane20,Lecu10,Yang16,Zhan18}. The algorithm
parameters, as well as possible hyperparameters of the problem, can
then be optimized from a training set by using differentiable
programming. Let us note that the results of \cite{Svva20,Smds20}
can be used to characterize the nonexpansiveness properties of the
resulting neural network \cite{Bert20}.

\subsection{Plug-and-play methods}
\label{sec:7b}
The principle of the so-called \emph{plug-and-play} (PnP) methods 
\cite{Buzz18,Onos17,Rond16,Ryue19,Suny19,Venka13} is to replace a
proximity operator appearing in some proximal minimization
algorithm by another operator $Q$. The rationale is
that, since a proximity operator can be interpreted as a denoiser
\cite{Smms05}, one can consider replacing this proximity operator
by a more sophisticated denoiser $Q$, or even learning it in a
supervised manner from a database of examples.
Example~\ref{ex:MMSE} described implicitly a PnP algorithm that
can be interpreted as a minimization problem. Here are some
techniques that go beyond the optimization setting.

\begin{algorithm}[PnP forward-backward]
Let $f\colon\HH\to\RR$ be a differentiable convex function, let
$Q\colon\HH\to\HH$, let $\gamma\in\RPP$, let 
$(\lambda_n)_{n\in\NN}$ be a sequence
in $\RPP$, and let $x_0\in\HH$. Iterate
\begin{equation}
\label{e:fbpnp}
\begin{array}{l}
\text{for}\;n=0,1,\ldots\\
\left\lfloor\begin{array}{l}
y_n=x_n-\gamma \nabla f(x_n)\\
x_{n+1}=x_n+\lambda_n(Qy_n-x_n).
\end{array}\right. 
\end{array}
\end{equation}
\end{algorithm}

The convergence of $(x_n)_{n\in\NN}$ in
\eqref{e:fbpnp} is related to the properties of 
$T=Q\circ(\Id-\gamma\nabla f)$.
Suppose that $T$ is $\alpha$-averaged with 
$\alpha\in\rzeroun$, and that $S=\Fix T\neq\emp$. 
Then it follows from Theorem~\ref{t:1}
that, if $(\lambda_n)_{n\in\NN}$ is an $\alpha$-relaxation 
sequence, then $(x_n)_{n\in\NN}$ converges to a point in $S$. 

\begin{algorithm}[PnP Douglas-Rachford]
\label{a:DR}
Let $f\in\Gamma_0(\HH)$, let $Q\colon\HH\to\HH$, let
$\gamma\in\RPP$, let $(\lambda_n)_{n\in \NN}$ be a sequence in
$\RPP$, and let $x_0\in\HH$. Iterate
\begin{equation}
\label{e:m1}
\begin{array}{l}
\text{for}\;n=0,1,\ldots\\
\left\lfloor\begin{array}{l}
x_n=\prox_{\gamma f} y_n\\
y_{n+1}=y_n+\lambda_n \big(Q(2x_n-y_n)-x_n\big).
\end{array}\right. 
\end{array}
\end{equation}
\end{algorithm}

In view of \eqref{e:m1},
\begin{equation}
\label{e:dr}
(\forall n\in\NN)\quad y_{n+1}=
\Big(1-\frac{\lambda_n}{2}\Big)y_n+\frac{\lambda_n}{2}Ty_n,
\end{equation}
where $T=(2Q-\Id)\circ(2\prox_{\gamma f}-\Id)$. Now assume that
$Q$ is such that $T$ is $\alpha$-averaged for some
$\alpha\in\rzeroun$
and $\Fix T\neq\emp$. Then it follows from Theorem~\ref{t:1}
that, if $(\lambda_n/2)_{n\in\NN}$ is an $\alpha$-relaxation 
sequence, then $(y_n)_{n\in\NN}$ converges to a point in 
$\Fix T$ and we deduce that $(x_n)_{n\in\NN}$
converges to a point in $S=\prox_{\gamma f}(\Fix T)$.
Conditions for $T$ to be a Banach contraction in the two previous
algorithms are given in \cite{Ryue19}.

Applying the Douglas-Rachford algorithm to the dual of Problem
\ref{prob:5} leads to a simple form of the alternating direction
method of multipliers. Thus, consider
algorithm~\ref{a:DR}, where $f$, $\gamma$, and $Q$ are 
replaced by $f^*$, $1/\gamma$ and $\Id+\gamma^{-1}Q(-\gamma\cdot)$,
respectively, and $(\forall n\in\NN)$ $\lambda_n=1$. Then we obtain
the following algorithm \cite{Chan17}, which is applied to image
fusion in \cite{Teod19}.

\begin{algorithm}[PnP ADMM]
\label{a:ADMM}
Let $f\in\Gamma_0(\HH)$, let $Q\colon\HH\to\HH$, let 
$\gamma\in\RPP$, let $y_0\in\HH$, let $z_0\in\HH$, and let 
$\gamma\in\RPP$. Iterate
\begin{equation}
\begin{array}{l}
\text{for}\;n=0,1,\ldots\\
\left\lfloor\begin{array}{l}
x_n=Q(y_n-z_n)\\
y_{n+1}=\prox_{\gamma f}(x_n+z_n)\\
z_{n+1}=z_n+x_n-y_{n+1}.
\end{array}
\right.
\end{array}
\end{equation}
\end{algorithm}

Note that, beyond the above fixed point descriptions of $S$, the 
properties of the solutions in plug-and-play methods are elusive 
in general. 

\subsection{Adjoint mismatch problem}
\label{sec:7c}
A common inverse problem formulation is to 
\begin{equation}
\minimize{x\in\HH}{f(x)+\frac{1}{2}\|Hx-y\|^2+
\frac{\kappa}{2}\|x\|^{2}},
\end{equation}
where $f\in\Gamma_{0}(\HH)$, $y\in\GG$ models the observation,
$H\colon\HH\to\GG$ is a linear operator, and $\kappa\in\RP$.
This is a particular case of Problem~\ref{prob:5} where 
\begin{equation} 
g=\frac12\|H\cdot-y\|^2+\frac{\kappa}{2}\|\cdot\|^{2}, 
\end{equation}
has Lipschitzian gradient 
$\nabla g\colon x\mapsto H^*(Hx-y)+\kappa x$.
It can therefore be solved via 
Proposition~\ref{p:fb17}, which therefore 
requires the application of the
adjoint operator $H^*$ at each iteration.
Due to both physical and computational limitations in
certain applications, this adjoint may be 
hard to implement and it is replaced by a linear approximation 
$K\colon\GG\to\HH$ \cite{Lore18,Zeng20}. This leads to a
surrogate of the proximal-gradient scheme \eqref{e:FB2} of the
form 
\begin{multline}
\label{e:FBt}
(\forall n\in\NN)\quad 
x_{n+1}=x_{n}+\\
\lambda_n\Big(\prox_{\gamma f}\big((1-\gamma\kappa)x_{n}-
\gamma K(Hx_n-y)\big)-x_n\Big),
\end{multline}
with $\gamma\in\RPP$ and $\{\lambda_n\}_{n\in\NN}\subset\rzeroun$.
Let us assume that $L=K\circ H+\kappa \Id$ is a cocoercive
operator. Then the above algorithm is an instance of the
forward-backward splitting algorithm introduced in
Proposition~\ref{p:fb13} to solve Problem~\ref{prob:7} where
$A=\partial f$ and $B=L\cdot -K y$. This means that a solution
produced by algorithm~\eqref{e:FBt} no longer solves a minimization
problem since $L$ is not a gradient in general
\cite[Proposition~2.58]{Livre1}. However, suppose that $g$ is
$\nu$-strongly convex with $\nu\in\RPP$, 
let $\zeta_{\rm min}$ be
the minimum eigenvalue of $L+L^*$,
set $\chi=1/(\nu+\zeta_{\rm min})$,
let $\widehat{x}$ be the
solution to Problem~\ref{prob:5}, and let $\widetilde{x}$ be the
solution to Problem~\ref{prob:7}. Then, as shown in 
\cite{Chou21},
\begin{equation}
\label{e:biasgen}
\|\widetilde{x}-\widehat{x}\|\leq
\chi \,\|(H^{*}-K)(H\widehat{x}-y)\|.
\end{equation}
A sufficient condition ensuring that $L$ is cocoercive is
that $\zeta_{\rm min}>0$. The problem of adjoint mismatch when
$f=0$ is studied in \cite{Dong19}.

\subsection{Problems with nonlinear observations}

We describe the framework presented in \cite{Eusi20,Ibap20} to
address the problem of recovering an ideal object
$\overline{x}\in\HH$ from linear and nonlinear transformations
$(r_k)_{1\leq k\leq q}$ of it.

\begin{problem}
\label{prob:z}
For every $k\in\{1,\ldots,q\}$, let $R_k\colon\HH\to\GG_k$ and let
$r_k\in\GG_k$. The task is to 
\begin{equation}
\label{e:z1}
\text{find}\;x\in\HH\;\:\text{such that}\;\:
(\forall k\in \{1,\ldots,q\})\;\:R_kx=r_k.
\end{equation}
\end{problem}

In the case when $q=2$, $\GG_1=\GG_2=\HH$, and $R_1$ and $R_2$ 
are projectors onto vector subspaces, Problem~\ref{prob:z} reduces
to the classical linear recovery framework of \cite{Youl78} which
can be solved by projection methods. We can also express
Problem~\ref{prob:1} as a special case of Problem~\ref{prob:z} 
by setting $m=q$ and 
\begin{equation}
\label{e:z2}
(\forall k\in\{1,\ldots,q\})\quad
r_k=0\quad\text{and}\quad R_k=\Id-\proj_{C_k}.
\end{equation}
In the presence of more general nonlinear operators, however, 
projection techniques are not applicable to solve \eqref{e:z1}.
Furthermore, standard minimization approaches such as minimizing
the least-squares residual $\sum_{k=1}^q\|R_kx-r_k\|^2$ typically
lead to an intractable nonconvex problem. Yet, we can employ 
fixed point arguments to approach the problem and design a provably
convergent method to solve it. To this 
end, assume that \eqref{e:z1} has a solution and that each 
operator $R_k$ is \emph{proxifiable} in the sense that 
there exists $S_k\colon\GG_k\to\HH$ such that 
\begin{equation}
\label{e:pr}
\begin{cases}
S_k\circ R_k\;\text{is firmly nonexpansive}\\
(\forall x\in\HH)\quad
S_k(R_kx)=S_kr_k\;\;\Rightarrow\;\; R_kx=r_k. 
\end{cases}
\end{equation}
Clearly, if $R_k$ is firmly nonexpansive, e.g., a projection or 
proximity operator (see Fig.~\ref{fig:1}), then it is proxifiable 
with $S_k=\Id$. Beyond that, many transformations found in data 
analysis, including discontinuous operations such as wavelet 
coefficients hard-thresholding, are proxifiable
\cite{Eusi20,Ibap20}. Now set
\begin{equation}
\label{e:21}
(\forall k\in\{1,\ldots,q\})\quad T_k=S_kr_k+\Id-S_k\circ R_k.
\end{equation}
Then the operators $(T_k)_{1\leq k\leq q}$ are firmly nonexpansive
and Problem~\ref{prob:z} reduces finding one of their
common fixed points. In view of Propositions~\ref{p:roma} and 
\ref{p:f3}, this can be achieved by applying Theorem~\ref{t:1} 
with $T=T_1\circ\cdots\circ T_q$. The more sophisticated
block-iterative methods of \cite{Nume06,Ibap20} are also
applicable.

Let us observe that the above model is based purely on a fixed 
point formalism which does not involve monotone inclusions or
optimization concepts. See \cite{Eusi20,Ibap20}
for data science applications.

\section{Concluding remarks}
\label{sec:8}
We have shown that fixed point theory provides an essential set of
tools to efficiently model, analyze, and solve a broad range of
problems in data science, be they formulated as traditional
minimization problems or in more general forms such as Nash
equilibria, monotone inclusions, or nonlinear operator equations.
Thus, as illustrated in Section~\ref{sec:7}, nonlinear models that
would appear to be predestined to nonconvex minimization methods
can be effectively solved with the fixed point machinery. The
prominent role played by averaged operators in the construction of
provably convergent fixed point iterative methods has been
highlighted. Also emphasized is the fact that monotone operators
are the backbone of many powerful modeling approaches. We believe
that fixed point strategies are bound to play an increasing role in
future advances in data science.

\medskip
{\bfseries Acknowledgment.} The authors thank Minh N. B\`ui and 
Zev C. Woodstock for their careful proofreading of the paper.

\newpage

\begin{IEEEbiography}
[{\includegraphics[width=1in,height=1.25in,clip,%
keepaspectratio]{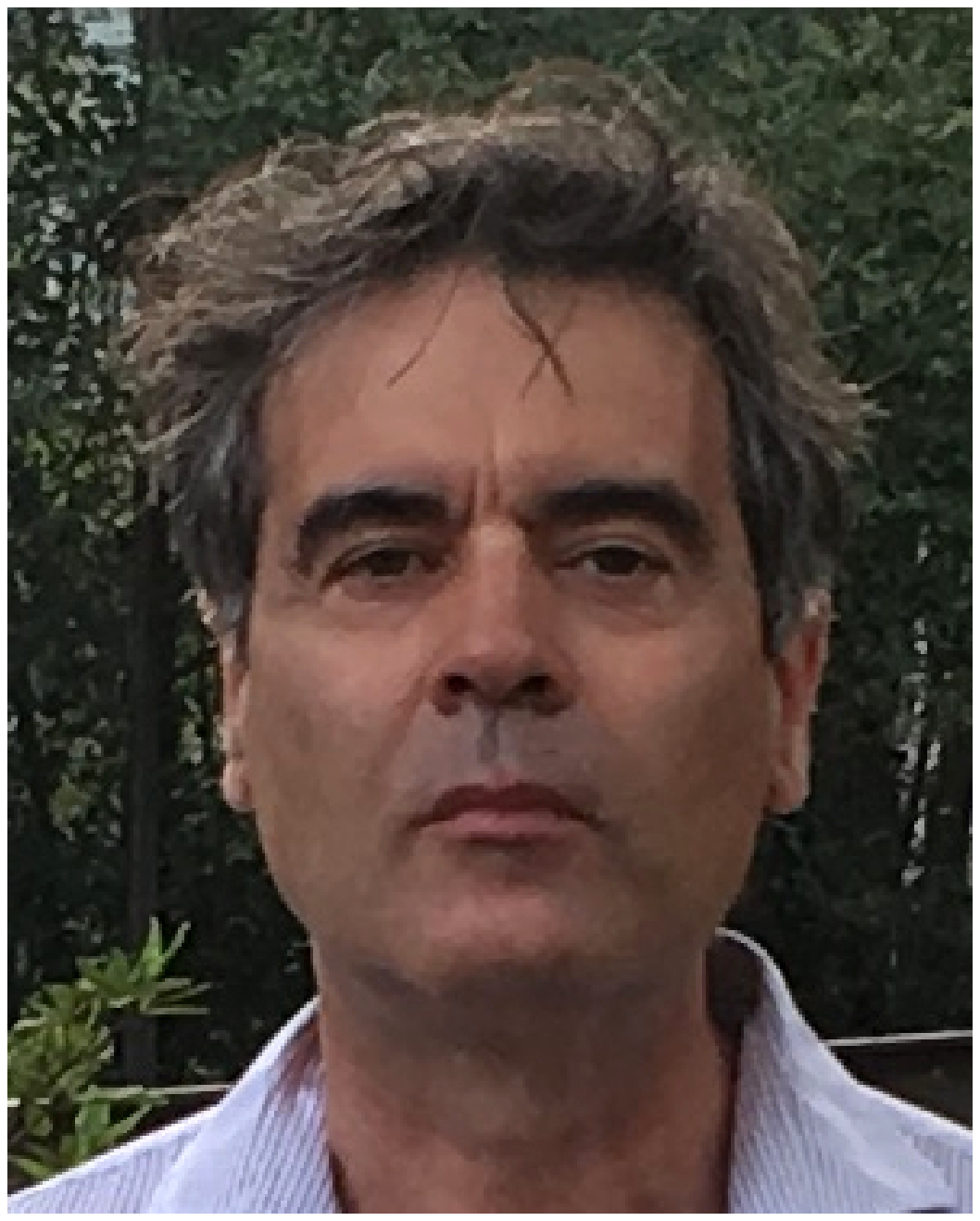}}]
{Patrick L. Combettes}
(Fellow, IEEE 2006) joined the faculty of the City University of
New York (City College and Graduate Center) in 1990 and the
Laboratoire Jacques-Louis Lions of Université Pierre et
Marie Curie--Paris 6 (now Sorbonne Universit\'e) in 1999.
He has been a Distinguished Professor of Mathematics at
North Carolina State University since 2016. He is the founding
director (2009-2013) of the CNRS research consortium 
MOA on mathematical optimization and its applications.
\end{IEEEbiography}

\vskip -153mm

\begin{IEEEbiography}
[{\includegraphics[width=1.2in,height=1.35in,clip,%
keepaspectratio]{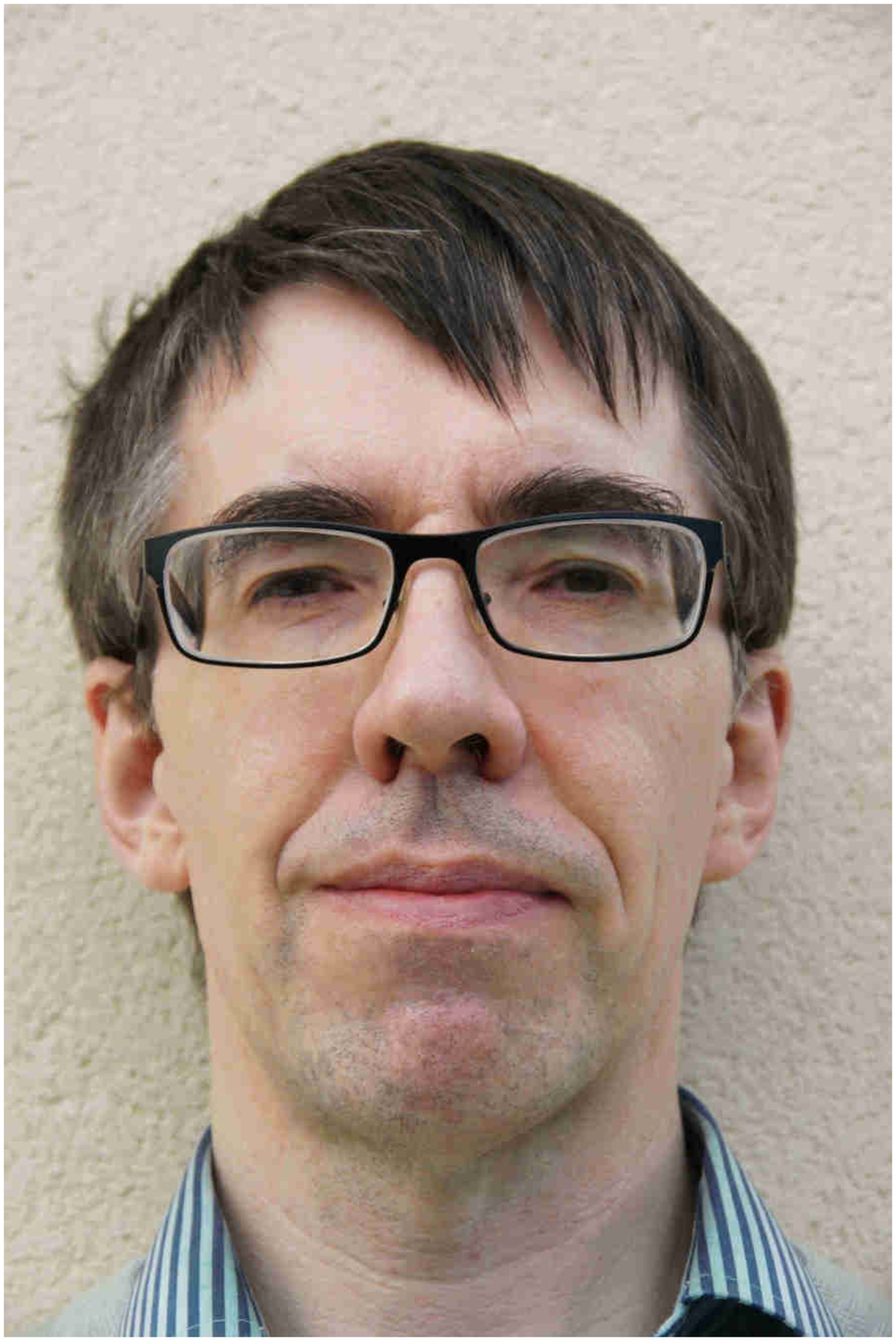}}]
{Jean-Christophe Pesquet}
(Fellow, IEEE 2012) received the engineering degree from Supélec,
Gif-sur-Yvette, France, in 1987, the Ph.D. and HDR degrees
from Universit\'e Paris-Sud in 1990 and 1999,
respectively. From 1991 to 1999, he was an Assistant
Professor at Universit\'e Paris-Sud, and a Research
Scientist at the Laboratoire des Signaux et Syst\`emes
(CNRS). From 1999 to 2016, he was a Full Professor at
Universit\'e Paris-Est and from 2012 to 2016, he was the
Deputy Director of the Laboratoire d’Informatique of the
university (CNRS). He is currently a Distinguished
Professor at CentraleSupélec, Universit\'e Paris-Saclay,
and the Director of the Center for Visual Computing and
OPIS Inria group. His research interests include
statistical signal/image processing and optimization
methods with applications to data science. He has also been
a Senior Member of the Institut Universitaire de France
since 2016.
\end{IEEEbiography}

\end{document}